\mag=\magstep1
\documentstyle{amsppt}
\input amsppt1
\pageheight{23.5true cm}
\pagewidth{15.5true cm}
\parindent=4mm
\baselineskip=13pt plus.1pt
\parskip=3pt plus1pt minus.5pt
\nologo
\NoRunningHeads
\NoBlackBoxes
\topmatter

\title Flips from 4-folds with isolated complete intersection 
singularities \\
whose downstairs have rational bi-elephants 
\endtitle
\author Yasuyuki Kachi 
\endauthor
\affil 
Department of Mathematics \\
University of Utah \\
Salt Lake City, UT 84112, U.S.A. \\
kachi \! \! \@ \! math.utah.edu
\endaffil

\abstract
{We shall investigate a flipping contraction 
$g : X \to Y$ from a 4-fold $X$ with at most 
isolated complete intersection singularities. 
If $Y$ has an anti-bi-canonical divisor (=bi-elephant) 
with only rational singularities, then $g$ 
carries an inductive structure chained up by blow-ups 
({\it La Torre Pendente\/}), and in particular 
the flip exists. This naturally contains Miles Reid's 
`{\it Pagoda\/}' as an anti-canonical divisor 
(=elephant) and its proper transforms.} 
\endtopmatter

\document
\head Contents 
\endhead
\vskip -6mm
$$
\align
&\S 0.\, \, \text{Introduction.} \\
&\S 1.\, \, \widetilde{E}_i \simeq \Bbb P^2. \\ 
&\S 2.\, \, \text{Freeness of $|-K_X|$ and its consequences.} \\
&\S 3.\, \, \text{The invariant } 
\varepsilon_P (Z \supset S). \\
&\S 4.\, \, \text{The normal bundle } N_{E/X}. \\
&\S 5.\, \, \text{Widths (after M.~Reid).} \\
&\S 6.\, \, \text{Deformations of contractions } g : X \to Y. \\
&\S 7.\, \, \operatorname{Sing} X = \{P\}, \, \, \,  
\varepsilon_P (X \supset E) = 1. \\
&\S 8.\, \, \text{Existence of flips \, 
--- \, La Torre Pendente.} \\
&\S \text{Appendix.} \, \, \text{On terminal complete 
intersection singularities} \\
& \qquad \text{--- \, S.~Ishii's theorem.}
\endalign
$$

We will work over $\Bbb C$, the complex number field.

\newpage

\head \S 0. \ \ Introduction.
\endhead

\proclaim{0.1. Minimal Model Conjecture}\ \ 

Let $X$ be a smooth projective algebraic 
variety of dimension $n$. Then $X$ is 
birationally equivalent to 
a projective variety $X'$ which has only 
{\it terminal singularities\/} 
(Reid [R1]), such that either one of the following 
holds: 

\flushpar
(0.1.1) \, \, \, $X'$ is a 
{\it minimal model, i.e.\/} $K_{X'}$ is 
nef, or

\flushpar
(0.1.2) \, \, \, there exists a surjective 
morphism $\varphi : X' \to Y$ with connected 
fibers such that $-K_{X'}$ is $\varphi$-ample, 
and $\dim Y \leq n-1$.  
\endproclaim

For dimension 2, this was classically 
known to be true by the works of the Italian 
school in the 19-th century. For dimension 
3 or more, however, this has been left to be 
unknown for a hundred of years. 

The attempt overcoming this has been made 
by the idea of introducing 
{\it extremal rays\/} (Mori [Mo2]) in early 80's, 
and a program has been raised by several people 
such as Reid [Re1], Kawamata [Kaw1,2], Shokurov [Sh1] 
and others, how to approach minimal models in terms of 
extremal rays, which is nowadays known as 
the so-called {\it Minimal Model Program\/} 
(see for {\it e.g.\/} [KaMaMa] for introduction). 
Especially, it has been made it clear that the 
only obstacle toward the whole problem 
is concentrated on the appearance of 
{\it small contractions, i.e.\/} 
those birational contractions of extremal rays 
with higher codimensional exceptional loci. 

\definition{0.2. Definition}\ \ 
Let $X$ be a projective algebraic 
variety of dimension $n$ 

\flushpar
(0.2.1) \, \, \, with 
only $\Bbb Q$-factorial terminal singularities, 

\flushpar
and let $g : X \to Y$ be the contraction 
of an extremal ray. Then a priori 

\flushpar
(0.2.2) \, \, \, the relative Picard number 
$\rho(X/Y)=1$, and 

\flushpar
(0.2.3) \, \, \, $-K_X$ is $g$-ample. 

\flushpar
Then we say that $g$ is a {\it small contraction\/}, 
or alternatively a {\it flipping contraction\/}, if 

\flushpar
(0.2.4) \, \, \, $\dim \text{Exc }g \leq n-2$. 

If assume that there exists a projective variety 
$X^+$ and a birational morphism 
$g^+ : X^+ \to Y$ such that 

\flushpar
(0.2.1)$^+$ \, \, \, $X^+$ has only 
$\Bbb Q$-factorial terminal singularities, 

\flushpar
(0.2.2)$^+$ \, \, \, $\rho(X^+/Y) =1$, 

\flushpar
(0.2.3)$^+$ \, \, \, $K_{X^+}$ is $g^+$-ample, 
and 

\flushpar
(0.2.4)$^+$ \, \, \, $\dim \text{Exc }g^+ \leq n-2$, 

\flushpar
then we call this $g^+$ (and also the composite 
birational map 
$\psi := g^{+ \, -1} \circ g: 
X \dashrightarrow X^+$,) the {\it flip\/} 
of $g$. 
\enddefinition

A concrete example of such kind of transformation 
is first constructed by P.~Francia [Fra]. 

Now we state the main conjecture, 
which implies the Minimal Model Conjecture: 

\proclaim{0.3. Flip Conjecture}\ \ 

\flushpar
(1) \, \, (Existence) \, \, \, 
For any $g : X \to Y$ of Definition 0.1, 
the flip $g^+$ exists. 

\flushpar
(2) \, \, (Termination) \, \, \, 
There is no infinite sequence of flips: 
$$
\not\exists \qquad 
X \dashrightarrow X^+ \dashrightarrow 
X^{++} \dashrightarrow \dots. 
$$
\endproclaim

\vskip 2mm

\flushpar
{\bf Remark.}\ \  
This formulation does work for varieties 
with a certain class of singularities which 
is worse than terminal, or for 
varieties together with divisors, 
called {\it log-flips\/} 
(Shokurov [Sh2], Kawamata [Kaw7]). 
In fact these are harder and more complicated. So 
it is sometimes necessary to distinguish 
the above from these, and in such a case we call 
the above flip specifically {\it terminal flips\/}. 

\vskip 2mm

The termination part was first proved in dimension 3 
by Shokurov [Sh1]. Here it should be noted that 
he discovered a numerical invariant 
coming from the discrepancy, which he called 
the {\it difficulty\/} [loc.cit], played an essential 
role. Then the termination was generalized to 
dimension 4 by Kawamata--Matsuda--Matsuki [KaMaMa]. 
On the other hand, the existence part 
was also investigated by several people, 
first for the semi-stable 3-fold case by 
Kawamata [Kaw3], Tsunoda, and Shokurov, 
and later on for the general 3-fold case by Mori 
[Mo4], applying Kawamata's method. 
In this way the Minimal Model Conjecture 
has been found to be true also in dimension 3. 

There are also further developments 
and generalizations in dimension 3, 
such as Reid [Re5], Koll\'ar--Mori [KoMo], 
Shokurov [Sh2] (the existence of log-flips), 
Kawamata [Kaw7] (the termination of log-flips), 
[Kaw8], and Corti [Co] (see also [Utah]). 

So we are interested in the existence problem of flips 
in dimension 4. Because the problem is now local on the 
downstairs $Y$, we may discuss everything from now on 
under the following setting (as [Kaw3], [Mo4], [KoMo] 
did in dimension 3):

\definition{0.4. Notation--Assumption}\ \ 
Let $g : X \to Y$ be a proper bimeromorphic morphism 
from a 4-dimensional analytic space ({\it 4-fold\/} in short) 
$X$ with only terminal singularities to a normal 4-fold $Y$. 
Assume that the exceptional locus $E := \text{Exc }g$ 
is compact, connected, and that $X$ is a sufficiently small 
analytic neighborhood of $E$. Then $g$ is called an 
{\it (analytic) flipping contraction\/} 
if 

\flushpar
(0.4.1) \, \, \, $\dim E \leq 2$ and $-K_X$ is $g$-ample. 

\flushpar
We simply write this 
$X \supset E \overset{g}\to{\longrightarrow}
Y \supset g(E)$. 

A proper bimeromorphic morphism 
$g^+ : X^+ \to Y$ from a 4-fold $X^+$ with only terminal 
singularities is called the {\it flip\/} of $g$ if 

\flushpar
(0.4.1)$^+$ \, \, \, $\dim \text{Exc }g^+ \leq 2$, 
and $K_{X^+}$ is $g^+$-ample. 
\enddefinition

In this direction, the first achievement has been reached by 
Kawamata [Kaw4] in the case of {\it smooth\/} 4-folds: 

\proclaim{0.5. Theorem}\ \ (Kawamata [Kaw4], Structure theorem 
of flips from smooth 4-folds)

Let $X \supset E \overset{g}\to{\longrightarrow} Y \supset g(E)$ 
be a flipping contraction as in 0.4. Assume that $X$ is a 
smooth 4-fold. Then 
$E \simeq \Bbb P^2$, $\operatorname{Bs} |-K_X| = \emptyset$, 
and the normal bundle $N_{E/X} 
\simeq \Cal O_{\Bbb P^2}(-1)^{\oplus 2}$. 

In particular, the flip $g^+ : X^+ \to Y$ of $g$ exists. 
\endproclaim

This should be considered as the first step of generalizing the 
Minimal Model Conjecture to dimension 4. 

\definition{0.6. Remark}\ \ 
In this case the downstairs $Y$ has 
{\it rational bi-elephants, i.e.\/} $|-2K_Y|$ has a member with 
only rational singularities. ({\it cf.\/} [Nak].)
\enddefinition

Although the assumption $X$ being smooth seems at a first glance 
quite casual and collects less attention, the conclusion 
of the theorem says that it is rather strong in this context, 
and it in fact determines the structure of $g$ essentially in 
a unique way. Also in this case the flip can be composed just by 
a single blow up-and-down, and so we may say this is the 
simplest flip in dimension 4. 
Toward the Minimal Model Conjecture, it is necessary 
however to generalize this to the singular case, {\it i.e.\/} 
the case that $X$ has terminal singularities. 
 
In his previous paper [Kac2], the author investigated flipping 
contractions from {\it semi-stable\/} 4-folds whose degenerate 
fibers satisfy a certain special assumption on singularities. 
Especially it is found that even we naively generalize 
Kawamata [Kaw3,7]'s definition of 3-fold semi-stable degenerations 
to the case of 4-folds, flips may in general destroy this 
condition, contrary to dimension 3. So we have to 
look for an appropriate re-definition of semi-stability which is 
preserved under any flips and divisorial contractions. 
This might be a problem to confront with in a future. 

So we turn back to {\it non-semi-stable\/} flipping 
contractions. In this paper, we deal with flipping 
contractions from 4-folds with isolated complete intersection 
singularities. 
The following is the object we are going to investigate: 

\definition{0.7. Assumption A}\ \ 
Let $g : X \to Y$ be as above. 
Assume 

\flushpar
(A-1)\, \, \, (Singularities on $X$)

\flushpar
$X$ has only isolated terminal complete 
intersection singularities, and 

\flushpar
(A-2)\, \, \, ({\it Existence of a good bi-elephant on\/} 
$Y$) ({\it cf.\/} 0.6 above)

\flushpar
$|-2K_Y|$ has a member with only rational singularities. 

A flipping contraction $g : X \to Y$ is said to be of 
{\it Type (R)\/} if $g$ satisfies both (A-1) and (A-2). 
$g$ is said to be of {\it Type (I)\/} if $g$ satisfies 
(A-1) and does {\it not\/} satisfy (A-2). 
\enddefinition

Then our main result is stated as follows: 

\proclaim{0.8. Main Theorem} \ \ 
Let $X \supset E \overset{g}\to{\longrightarrow} Y \supset g(E)$ 
be a 4-fold flipping contraction satisfying the 
Assumption A above, {\it i.e.\/} of Type (R). 
Assume $\operatorname{Sing} X (\cap E) \not= \emptyset$. 
Then 
$$
E \simeq \Bbb P^2, \quad \operatorname{Bs} |-K_X| = \emptyset, 
\quad N_{E/X} \simeq \Cal O_{\Bbb P^2} \oplus 
\Cal O_{\Bbb P^2}(-2). 
$$
Moreover, the flip $g^+$ exists. 

(As for the description of the flip $g^+ : X^+ \to Y$, 
see Corollary 8.11.)
\endproclaim

To be more precise, $g$ is classified in terms of the 
discrete invariant $\operatorname{width} g$ (see \S 5), 
and $g$ carries an inductive structure chained up by blow-ups 
with respect to this number ({\it La Torre Pendente\/}), 
and this gives a complete understanding of the mechanism 
of flips (see \S 8).

After the completion of this paper, the author succeeded in 
proving the implication ``(A-1) $\Longrightarrow$ (A-2)'' so that 
we obtained the following stronger theorem: 

\proclaim{0.9. Theorem}\ \ (Existence of flips from 4-folds  
with isolated complete intersection singularities) 

Let $X$ be a projective 4-fold, with at most 
isolated complete intersection singularities. 
Then for any flipping contraction $g : X \to Y$, 
the flip $g^+$ exists. 
\endproclaim

The proof will be given in the forthcoming paper. 

What should be coming next is to investigate $g : X \to Y$ 
where $X$ has a {\it non-isolated\/} singular locus. 
For example it is unknown whether the condition 

\flushpar
(A-1)$'$ \, \, \, $X$ has terminal complete intersection 
singularities with $\dim \text{Sing }X = 1$, 

\flushpar
implies (A-2), while on the other hand there is an example 
of a flipping contraction $g : X \to Y$ 
satisfying both (A-1)$'$ and (A-2) by M.~Gross 
(which will be appeared in the forthcoming paper).

This paper is organized as follows: 

The first couple of sections are preliminary, 
and the arguments are mostly parallel to Kawamata [Kaw4]'s.
({\it cf.\/} [AW1], [ABW].)

First in \S 1, we shall prove that the normalization 
$\widetilde{E}_i$ of each irreducible component $E_i$ 
of $E$ is isomorphic to $\Bbb P^2$. 
We apply Koll\'ar's theorem [Kol2] on Hilbert schemes 
parametrizing rational curves on a variety 
with complete intersection singularities, 
which is a generalization of Mori [Mo1] in the smooth case.
This part we are indebted to Koll\'ar's idea. All the rest 
is the reproduction of arguments in [Kaw4]. 
At the same time, it is proved that $-K_X$ pulled-back 
to $\widetilde{E}_i\text{'s}$ are all $\Cal O_{\Bbb P^2}(1)$. 

In \S 2, we shall discuss the freeness of $|-K_X|$ 
(Theorem 2.1). Kawamata [Kaw4] first proved this for the 
smooth case by developping the `{\it base-point-free 
technique\/}' initiated by himself, and later on 
Andreatta--Wi\'sniewski [AW1] extended this to a certain 
singular case, which is sufficient for what we need.
We shall give the proof for readers' convenience, 
following Kawamata [Kaw4]'s original argument.  
Also we shall collect several immediate consequences:
The normality of $E_i\text{'s}$ (Corollary 2.8), 
a cohomological interpretation of the condition (A-2) 
of 0.7 (Corollary 2.9), and the irreducibility of $E$
: $E \simeq \Bbb P^2$ for Type (R) contraction 
(Corollary 2.10). 

In \S 3, we shall introduce an invariant
$\varepsilon_P = \varepsilon_P (Z \! \supset \! S)$ for 
a pair of a germ of an analytic space $Z = (Z,P)$ and its 
closed subspace $S$ passing through $P$. 
We give an explicit formula (Theorem 3.1) for computing 
this $\varepsilon_P$ in the case that $Z$ is a complete 
intersection singularity and $S$ is a smooth subspace. 
Especially when $Z$ is a 3-dimensional terminal Gorenstein 
singularity (= $cDV${\it-singularity\/}) and $S$ is a smooth curve, 
then this coincides exactly with Mori [Mo4]'s invariant 
$i_P(1)$ (Corollary 3.2), and thus our $\varepsilon_P$ can be 
thought as a sort of a higher dimensional analogue of 
$i_P(1)$. 

So far is the preparation. 

From \S 4 on, we concentrate on Type (R) contractions. 

In \S 4, we shall determine the normal bundle  
$N_{E/X} := (I_E/I_E^2)^{\vee}$ of $E$ in $X$ 
for Type (R) contractions. We shall prove that 
$N_{E/X} \simeq \Cal O_{\Bbb P^2} \oplus \Cal O_{\Bbb P^2}(-2)$
(Theorem 4.1). The strategy is first to limit the types of 
its restrictions to lines (Proposition 4.2), where the 
assumption (A-2) is essentially used, and then to 
deduce the uniformity of $N_{E/X}$ (Proposition 4.4) 
by applying the generalized Namikawa [Nam3]'s local moduli 
(Theorem 4.5), where we develop the deformation theory for 
contraction morphisms. (See also Remark 4.9). Then Van de Ven [V]'s 
characterization of uniform bundles immediately gives us 
the conclusion. 

In \S 5, we shall introduce an invariant called {\it widths\/} 
for our contractions $g$ (Definition 5.5), after Reid [Re1] 
for 3-fold {\it flopping contractions\/}.  
Roughly, this measures the infinitesimal length of the 
embedding $E \subset X$ along a general locus of $E$, while 
it actually specifies the singularities of $X$ as 
we will see in \S 7 (Corollary 7.6). This is thus so to speak 
a local-to-global invariant, and in fact plays an essential 
role toward the existence of flips in \S 8.  

In \S 6, we discuss deformations of a contraction 
$g : X \to Y$. Since $R^i g_* \Cal O_X =0$ $(i=1,2)$, 
we can talk about deformations of the contraction 
$g : X \to Y$, rather than that of $X$ itself 
(see Koll\'ar--Mori [KoMo] \S 11). 
We shall prove that any set of local deformations 
of $X$ along the singular points $(X,P_i)$ $(i=1, \dots , r)$ 
$\big( \text{where } \{P_1, \dots , P_r\} 
:= \text{Sing }X \cap E \big)$
can be patched together to produce a global deformation 
of $X$, and of $g:X \to Y$ (Theorem 6.1). Especially $X$ is 
smoothable by a flat deformation. 
What we actually prove are the following couple of 
statements: 

\flushpar
(a) \, \, \, The {\it unobstructedness\/} of deformations 
(see Ran [Ra2], Kawamata [Kaw6], Nami-
\allowlinebreak
kawa [Nam1,2,3], Gross [G1,2]), 
and

\flushpar
(b) \, \, \, The surjectivity of \, \, 
$\text{Ext}^1_X (\Omega^1_X, \Cal O_X) 
\overset{\alpha}\to{\longrightarrow} 
g_* \Cal E \! xt^1_X (\Omega^1_X, \Cal O_X)$. 

\flushpar
These imply the desired surjectivity of the natural holomorphic map

\qquad \qquad \qquad \qquad \qquad \quad \, \, \, \, 
$\text{Def }X \longrightarrow 
\prod\limits_{i=1}^r \text{Def } (X,P_i)$ 

\flushpar
between the global and the local {\it versal deformation spaces\/} 
(={\it Kuranishi spaces\/}). 
In dimension 3, the corresponding result was proved by Mori 
[Mo4] \S 1b in rather complex analytic languages 
({\it L-deformations\/}), including an essential use of the 
implicit function theorem. 
Saying in algebro-geometric words, his argument might correspond 
to the fact that the obstruction lies in $R^2 g_* \Cal T_X$, 
which automatically vanishes in dimension 3. In dimension 4, 
however, the proof of $R^2 g_* \Cal T_X=0$ 
essentially requires the structure of the normal bundle 
$N_{E/X}$, which we have already determined in \S 4. 
This can be viewed as a sort of a relative 4-dimensional 
analogue of Namikawa [Nam1]'s smoothing theory of Fano 
varieties. ({\it cf.\/} [Nam2,3], [G1,2].)

Also, for such given family 
$\Cal X \to \Cal Y = \{X_t \overset{g_t}\to{\longrightarrow}
Y_t\}_{t \in \Delta}$ 
we especially focus on the ideal structure 
of $E$, which is defined naturally as the degeneration 
of $E_t := \text{Exc }g_t$ $(t \not= 0)$ (6.7). 
We shall prove that this subscheme supported on $E$ 
has no embedded components, 
and that this has the multiplicity just equal to the number 
of connected components of $E_t$ $(t \not= 0)$ 
(Proposition 6.8). 
This observation is crucial toward the subsequent arguments. 
As opposed to Mori's L-deformation ([Mo4], [KoMo]), the basic 
idea of which was to restrict attention to one irreducible 
component after a deformation and forget everything else, ours 
is rather to look at all irreducible components simultaneously, 
and this in fact makes us possible to get desirable 
informations on the local structure of $X \supset E$ 
along singularities. For instance, at the final stage of 
the local classification (\S 7), we will show that the singularity 
$(X, P)$ can be completely determined by the width (Corollary 7.6), 
which is a primary consequence of this observation. 
This may thus shed a new light toward the problem of flips.

Based on these preparations, in \S 7 we now are going to 
classify $g$ of Type (R) in \S 7. We shall prove that 
$X$ has exactly one singular point $P$, and 
that $\varepsilon_P (X \supset E) = 1$ (Theorem 7.1). 
Hence $g$ is completely classified 
up to the invariant $\text{width }g$, introduced in \S 5. 
We shall make full use of the deformation theory 
established in \S 6. 
Technically, the hardest part is to show that every 
singular point of $X$ is a hypersurface 
singularity (Proposition 7.2). This can be done 
by giving several different local deformations and 
looking at the associated subscheme structures of $E$ 
for every corresponding global deformations. 

Finally in \S 8 we shall prove the existence of 
flips for Type (R) cxontractions. We shall prove 
this by induction on $m := \text{width }g$. 
The case $m = 1$ is nothing but Kawamata [Kaw4]'s 
result (Theorem 0.5). 
The construction of the flip is as follows: 
First blow $X$ up with the center $E$: 
$\varphi: \overline{X} \to X$. 
Then as we will see in Proposition 8.7, which is the hardest part 
of this section, there exists a unique flipping contraction 
$\overline{g}$ from $\overline{X}$ of Type (R), 
with $\text{width }\overline{g} = m-1$. (This is also a 
consequence of the deformation theory.) So by the 
induction hypothesis we are able to flip it: 
$\overline{X} \dashrightarrow \overline{X}^+$. 
Then the proper transform of $\text{Exc }\varphi$ 
in $\overline{X}^+$ can be contracted down; $\overline{X}^+ 
\to X^+$, to produce the desired flip $X^+ \to Y$. 
Here we essentially use the structure of Reid [Re1]'s 
`{\it Pagoda\/}' (see 5.1). 

To see the mechanism of this flip operation more concretely, 
let us follow all the induction steps successively upstreams, 
then we will eventually arrive $m = 1$, and will get a sequence 
of blow-ups and blow-downs:
$$
X \leftarrow X^{(1)} \leftarrow ... \leftarrow X^{(m-1)} 
\leftarrow X^{(m)} \to X^{(m-1) \, +} \to ... \to 
X^{(1) \, +} \to X^+,
$$
where $X^{(i)} \dashrightarrow X^{(i) \, +}$ gives the 
flip, with the width $m-i$. This should be compared to 
Reid's Pagoda, not only since the patterns of 
the constructions of the flip or the flop look similar, 
but also since ours in fact contains Pagoda as an 
anti-canonical divisor (and its proper transforms). 
On the other hand, though Pagoda was symmetric with respect to 
the flop (see 5.1), ours is not. 
So with a great esteem for Reid's humor of this lovely naming, 
we name by a special grace our 4-fold flip 
`{\it La Torre Pendente\/}'.

\definition{Acknowledgements}\ \ 

I would like to express my gratitude to Professor J.~Koll\'ar. 
This work has been started when I stayed University of Utah in 
April 1995, where he suggested me Theorem 1.2, which was 
the practical starting point of this work. 
Also I was invited by him to Utah in the spring semester, 1996. 
I appreciate his hospitality, and also I was much stimulated by 
his mathematical ingenuity.  

I am grateful to Professor M.~Andreatta, who invited me to 
Universit\`a di Trento in winter semester of 1996, and gave me 
much hospitality and entertainment. 

I am also grateful to Professor Yo.~Namikawa, who communicated 
me all through, and gave me great many of important suggestions. 
Especially he told me Theorem 5.4 and Example 4.8. Also I got 
Theorem 4.5 out of the conversation with him. These 
actually form the key of this paper.

I am grateful to Professor M.~Reid who showed me special generosity 
during my stay at University of Warwick in the fall term of 1995. 

During the preparation of this work, I was also invited by 
the following universities: 
Universit\"at Bayreuth, 
Universit\`a di Genova, 
Universit\`a di Torino, 
The Johns Hopkins University, 
Cornell University, 
University of Georgia, and 
Purdue University. 
I would like to thank Professors 
Th.~Peternell, 
M.~Beltrametti, 
A.~Conte, 
Y.~Kawamata, 
V.~Shokurov, 
M.~Gross, 
V.~Alexeev, and 
K.~Matsuki 
for their hospitalities. 

I would also like to thank Professors 
S.~Ishii, Yu.~Prokhorov and K.-i.~Watanabe who 
informed me of the structure of terminal complete 
intersection singularities (\S Appendix), 
and Professor A.~Degtyarev who indicated me 
a topological point of view.

Finally, I again thank Professor Y.~Kawamata, my teacher. 
He suggested me the necessity to formulate a deformation theory 
of a pair $(X, E)$ consisting of a variety $X$ together with its 
subspace $E$ meeting the singular locus of $X$, 
which led me later to inspire the observation in \S 6 
(Definition 6.7 and Proposition 6.8). 
Also, he so generously allowed these journeys and gave me 
consistent encouragement. I appreciate him very much for these 
kindness. 
\enddefinition

\newpage

\head \S 1.\  $\widetilde{E}_i \simeq \Bbb P^2$ 
\endhead

\definition{Notation 1.0}\ \
Let $X \supset E \overset{g}\to{\longrightarrow} Y$ 
be a flipping contraction satisfying 
the Assumption (A-1), {\it i.e.\/} 
of Type (R) or (I). Throughout this section, 
we fix the following notation: 

\flushpar
(1.0.1) \, \, \, 
Let $E = \bigcup\limits_{i=1}^n E_i$
be the irreducible decomposition of $E$, and 
$\nu_i : \widetilde{E}_i \to E_i$ the 
normalization of $E_i$ $(i = 1, \dots , n)$. 
\enddefinition

The aim of this section is to prove the following: 

\proclaim{Theorem 1.1} \qquad 
$\widetilde{E}_i \simeq \Bbb P^2$, \, and 
\, $\nu_i^* \Cal O_{E_i}(-K_X) \simeq \Cal O_{\Bbb P^2}(1)$ 
 \, \, $(i = 1, \dots , n)$. 
\endproclaim

First we shall give the following theorem due to J.~Koll\'ar [Kol2], 
which is a generalization of Mori [Mo1], and is one of the key to 
our whole argument:

\proclaim{Theorem 1.2}\ \ (Koll\'ar [Kol2] II 1.3 Theorem)

Let $Z$ be an $n\text{-dimensional}$ analytic space 
and $C$ an irreducible rational curve on it. 
Assume that $Z$ has at most complete intersection singularities, 
and that $C \not\subset \text{Sing } Z$. 
 
Let 
$\alpha : \widetilde{C} \simeq \Bbb P^1 \to C \subset Z$ 
be the normalization of $C$. Then 

\flushpar
\qquad \qquad \qquad \qquad \qquad \, \, 
$\dim \text{Hom }_{[\alpha]} (\Bbb P^1, Z) 
\geq n + (-K_Z \, . \, C)$. \quad \qed
\endproclaim

\definition{Remark 1.3}\ \ 
Koll\'ar's Theorem 1.2 has many far reaching consequences. 
For instance, Mori [Mo4]'s result  
``there is no flipping contraction from a 
terminal Gorenstein 3-fold" is now an 
immediate corollary of Theorem 1.2. 
\enddefinition

Thus we can proceed completely the same argument as 
in Kawamata [Kaw4] in the smooth case, which is an application 
of Mori [Mo1], Ionescu [Io] ({\it cf.\/} Wi\'sniewski [W]) 
to deduce: 

\proclaim{Proposition 1.4}\ \ (following Kawamata [Kaw4]. 
See also Andreatta--Ballico--Wi\'sniewski [ABW].) 

\flushpar
(1) \, \, \, $E$ is purely 2-dimensional, and 

\flushpar
(2) \, \, \, Each $E_i$ is covered by rational curves $C$ such that 
$(-K_X \, . \, C) = 1$. 
\endproclaim

\demo{Proof}\ \ 
First by Kawamata [Kaw5], each $E_i$ is 
covered by rational curves. Thus by the assumption 
(0.4.4), (1) is immediate from Theorem 1.2. 

Next, let $C$ be a rational curve in $E_i$ 
which attains the minimum number:

\flushpar
(1.4.0) \qquad \qquad \, \, \, 
$r := \text{Min } \{ (-K_X \, . \, C) \, | \, 
C: \text{ rational curves in } E_i \}$. 

\flushpar
Let $\alpha : \widetilde{C} \to C$ be the normalization of $C$. 
Then 

\flushpar
(1.4.1) \qquad \qquad \qquad \qquad \quad \, \, 
$\dim \text{Hom}_{[\alpha]} (\Bbb P^1, X) \geq 4+r$,  

\flushpar
by Theorem 1.2. Hence there exists a family of curves 
$p : \Cal U \to H$, with the projection 
$q : \Cal U \to X$ satisfying $q(\Cal U) = E_i$, 
which gives the universal deformation 
of $C$ (with a certain generically reduced closed subscheme 
structure) inside $X$. 
By the choice of $C$, 

\flushpar
(1.4.2) \, \, \, All fibers of $p$ are irreducible rational curves, 

\flushpar
and 
$$
\align 
(1.4.3) \qquad \qquad \qquad \,
\dim H &= \dim \text{Hom}_{[\alpha]} (\Bbb P^1, X) 
- \dim \text{Aut } \Bbb P^1 
\qquad \qquad \qquad \qquad \quad \, \\ 
&\geq 4 + r - 3 = r + 1. 
\endalign
$$
From now on we shall prove $\dim H = 2$. Let $x$ be a 
general point of $E_i$, and let 
$$
\cases
H_x := p(q^{-1} (x)), \, \, \Cal U_x := p^{-1} (H_x), \\
p_x := p|_{\Cal U_x} : \Cal U_x \to H_x, \text{ and } 
q_x := q|_{\Cal U_x} : \Cal U_x \to X.
\endcases
\tag 1.4.4
$$
Then 

\flushpar
(1.4.5) \qquad \qquad \qquad \qquad \qquad \qquad \, \, 
$q_x(\Cal U_x) = E_i$, 

\flushpar
and 
$$
\align
(1.4.6) \qquad \qquad \qquad \, \, \, \, 
\dim H_x = \dim q^{-1}(x) &= \dim \Cal U - \dim E_i \\
&= \dim H -1. 
\qquad \qquad \qquad \qquad \qquad \quad \, \, \, 
\endalign
$$

Assume $\dim H \geq 3$ to derive a contradiction, till (1.4.10). 
Let $y$ be another general point of $E_i$, 
and let 

$$
\cases
H_{x,y} := p_x(q_x^{-1}(y)), \, \, \Cal U_{x,y} := p_x^{-1} (H_{x,y}), \\
p_{x,y} := p_x|_{\Cal U_{x,y}} : \Cal U_{x,y} \to H_{x,y}, \text{ and } 
q_{x,y} := q_x|_{\Cal U_{x,y}} : \Cal U_{x,y} \to X.
\endcases
\tag 1.4.7
$$
Then similarly we have 
$$
\cases 
q_{x,y} (\Cal U_{x,y}) = E_i, \, \text{ and} \\
\dim H_{x,y} \geq \dim H - 2 \, \, \, (1.4.6). 
\endcases
\tag 1.4.8
$$
If $\dim H \geq 3$, then this means that 

\flushpar
(1.4.9) \, \, \, $p_{x,y} : \Cal U_{x,y} \to H_{x,y}$ 
admits two disjoint sections 
$$
s_x := q_{x,y}^{-1} (x) \, \, \, \text{ and } \, \, 
s_y := q_{x,y}^{-1} (y). 
$$

Let $H_1 \subset H_{x,y}$ be 
any 1-dimensional irreducible closed subset, let 
$\widetilde{H}_1 \to H_1$ be the normalization, 
and $p_1 : \Cal U_{x,y} \times_{H_{x,y}} \widetilde{H}_1 
\to \widetilde{H}_1$ the induced morphism. 
Then 

\flushpar 
(1.4.10) \, \, \, $p_1$ still admits two disjoint sections 
$s_{x,1}$, $s_{y,1}$ induced from $s_x$, $s_y$, respectively, 
both of which are contractible. 

This is however impossible, since $p_1$ is a 
$\Bbb P^1\text{-bundle}$ over a smooth complete curve 
$\widetilde{H}_1$ (1.4.2). 
Hence we must have 

\flushpar
(1.4.11) \qquad \qquad \qquad \qquad \qquad \qquad \, 
$\dim H = 2$. 

\flushpar
In particular, the inequality (1.4.3) must be the equality: 

\flushpar
(1.4.12)  
\qquad \qquad \qquad \qquad \qquad \qquad \quad \, \,  
$r=1$. 
\quad \qed 
\enddemo

\proclaim{1.5}\ \ Proof of Theorem 1.1. 
\endproclaim

By the previous proposition, it is enough to prove 
$\widetilde{E}_i \simeq \Bbb P^2$, so let us return back 
to the situation (1.4.4) through (1.4.6). 
Note that $\dim H_x = 1$ by (1.4.6) and (1.4.11). Let 
$H_0$ be any 1-dimensional irreducible component 
of $H_x$, and $\widetilde{H}_0 \to H_0$ 
the normalization. Then the induced 
$$
p_0 : \widetilde{\Cal U}_0 := \Cal U_x \times_{H_x} \widetilde{H}_0
\to \widetilde{H}_0 
$$
is a $\Bbb P^1\text{-bundle}$ over a smooth 
complete curve $\widetilde{H}_0$ dominating $E_i$ 
through $q_0 : \widetilde{\Cal U}_0 \to X$. Moreover, 
$p_0$ has a section $s_x  := q_0^{-1} (x)$, by the construction. 
Hence the normalization of $E_i$ must be $\Bbb P^2$. 
\quad \qed

\definition{Remark 1.6}\ \ 
The conclusion of Theorem 1.1 is no longer true if we assume $X$ 
merely to be Gorenstein, instead of assuming $X$ to have complete 
intersection singularities (see \S Appendix A.1).
In fact there is an example (Example 8.13) due to Mukai, 
of a flipping contraction $g : X \to Y$ from a Gorenstein 
4-fold which contracts a singular quadric surface. 

A ruling $l$ satisfies $(-K_X \, . \, l) = 1$, while $l$ 
deforms inside $X$ exactly with 1-dimensional parameter space. 
So Kollar's formula (Theorem 1.2) also fails in this case. 
(For details see 8.13.)
\enddefinition

\vskip 5mm

\head \S 2.\  Freeness of $|-K_X|$ and its consequences. 
\endhead

In this section, we shall overview the proof of the freeness 
of the anti-canonical linear system $|-K_X|$ on $X$ (Theorem 2.1) 
(without the assumption (A-2)), following Kawamata [Kaw4] and 
Andreatta--Wi\'sniewski [AW1]. In the case $X$ is smooth, 
this was proved by Kawamata [Kaw4], and was extended later 
on to a more general situation, {\it i.e.\/} the case that 
$X$ has rational Gorenstein singularities, by 
Andreatta--Wi\'sniewski [AW1]. 

This result provides us three corollaries. The first one is 
the smoothness of each irreducible components of $E$ 
(Corollary 2.8). The second is an interpretation of 
the condition (A-2) in terms of a certain cohomology 
vanishing (Corollary 2.9). The final one is the irreducibility 
of $E$ under the asumption (A-2), {\it i.e.\/} 
for Type (R) (Corollary 2.10).

\proclaim{Theorem 2.1} \ \ (Kawamata [Kaw4], 
Andreatta--Wi\'sniewski [AW1].)

Let $X \supset E \overset{g}\to{\longrightarrow} 
Y \ni Q$ be a flipping contraction satisfying the 
assumption 
(A-1), {\it i.e.\/} of Type (R) or (I). 
Then $|-K_X|$ is free, {\it i.e.\/} the homomorphism 
$$
\rho : g^* g_* \Cal O_X (-K_X) \to \Cal O_X (-K_X) 
$$
is surjective. 
\endproclaim

Here we shall adopt Kawamata [Kaw4]'s original proof. 
We completely follow the argument of [loc.cit] from 2.2 through 
2.7 ({\it cf.\/} [Kac1] \S 4). 

Most of the notations fixed therein (such as $A$, $B$, 
$C$, $F_i$, $G_i$, $L$ {\it etc.\/}) are valid only 
for this section.

\definition{2.2}\ \ 
Assume $\text{Bs }|-K_X| 
:= \text{Supp } \text{Coker } \rho \not= \emptyset$ 
to derive a contradiction, till 2.7. 

Let $\varphi : X' \to X$ be a projective bimeromorphic 
morphism from a smooth 4-fold $X'$, 
together with a simple normal crossing divisor 
$\sum\limits_i G_i$ on it, such that 
$$ 
\cases
\varphi^* |-K_X| = |D'| + \sum\limits_i r_i G_i 
\, \text{ with } \, \text{Bs } |D'| = \emptyset, \\
K_{X'} = \varphi^*K_X + \sum\limits_i a_i G_i, \, \text{ and} \\
- \varphi^*K_X - \sum\limits_i \delta_i G_i \text{ is } 
(g \circ \varphi)\text{-ample}, 
\endcases
\tag 2.2.0
$$
for some $r_i \in \Bbb Z_{\geq 0}$, $a_i \in \Bbb Z$, 
and $\delta_i \in \Bbb Q_{\geq 0}$, $0 <\delta_i <1$.  
Since $X$ is assumed to have at most terminal singularities, 

\flushpar
(2.2.1) \qquad \qquad \qquad \qquad \qquad 
\qquad \quad \, \, \, \, 
$a_i \geq 0$. 

\flushpar
Note that 

\flushpar
(2.2.2) \qquad \qquad \qquad \qquad \, 
$\varphi(\bigcup\limits_i G_i) = \text{Bs } |-K_X| \cup \text{Sing }X$. 

Let $c := \min \dfrac{a_i + 1 - \delta_i}{r_i}$. 
By shrinking $\delta_i\text{'s}$ if necessary, we may assume 
that the minimum $c$ is attained exactly for a single $i$, 
say $i=1$. Let

\flushpar
(2.2.3) \qquad \qquad \qquad 
$A := \sum\limits_{i \geq 2}(-cr_i + a_i - \delta_i) G_i$, 
\, and \, $B := G_1$. 

\flushpar
Then 

\flushpar
(2.2.4) \qquad \qquad \qquad \quad \, \, \, \, \, 
$\ulcorner A \urcorner \geq 0$ 
\, and \, 
$\text{Supp } \ulcorner A \urcorner \subset 
\text{Exc } \varphi$. 
\enddefinition

\proclaim{Lemma 2.3} \ \ 
Under the notations and the assumption of 2.2, 
$$
a_i \geq r_i \, \, \, (\forall i). 
$$
\endproclaim 

\demo{Proof}\ \ 
Assume to the contrary, then 
$c < \dfrac{a_i + 1}{r_i} \leq 1$, and hence 
$$
\align
C := & -\varphi^* K_X - K_{X'} + (A - B) \\
= & cD' - (2-c) \varphi^* K_X - \sum\limits_{i \geq 1} \delta_i G_i 
\endalign
$$
is $(g \circ \varphi)\text{-ample}$. 
By the Kawamata-Viehweg vanishing theorem 
(for {\it e.g.\/} [KaMaMa]), 
$R^1 (g \circ \varphi)_* \Cal O_{X'} (-\varphi^* K_X 
+ \ulcorner A \urcorner - B) =0$, and thus 
$$
s : (g \circ \varphi)_* \Cal O_{X'} (-\varphi^* K_X 
+ \ulcorner A \urcorner) \to H^0 (B, \Cal O_B (-\varphi^* K_X 
+ \ulcorner A \urcorner))
$$
is surjective. Here 

\flushpar
(2.3.2) \qquad \qquad \quad \, \, \, 
$(g \circ \varphi)_* \Cal O_{X'} (-\varphi^* K_X 
+ \ulcorner A \urcorner) \simeq g_* \Cal O_X(-K_X)$ 

\flushpar
(2.2.4). On the other hand, 

\flushpar
(2.3.3) \qquad \qquad \qquad \qquad 
$H^0 (B, \Cal O_B (-\varphi^* K_X 
+ \ulcorner A \urcorner)) \not= 0$, 

\flushpar
by the fact that $\nu_i^* \Cal O_{E_i} (-K_X) \simeq 
\Cal O_{\Bbb P^2}(1)$ is globally generated $(\forall i)$. 
These imply $\varphi(B) \not\subset \text{Bs }|-K_X|$, 
while from (2.2.2) we must have 
$\varphi(B)\subset \text{Bs }|-K_X|$, 
a contradiction. Hence 

\flushpar
\qquad \qquad \qquad \qquad \qquad \qquad \qquad \quad \, \, \, \, 
$a_i \geq r_i \, \, \, (\forall i)$. \quad \qed 
\enddemo

\proclaim{Lemma 2.4}\ \ 
Let $D$ be a general member of $|-K_X|$. 
Then $D$ has at most canonical singularities. 
\endproclaim

\demo{Proof}\ \ 
Take a general smooth $D' \in |D'|$ (2.2.0) 
and consider $\varphi|_{D'} : D' \to D$. 
This gives a resolution of $D$. 
By adjunction and (2.2.0), 
$$ 
\align
K_{D'} = K_{X'} + D' \big|_{D'} &= \sum\limits_i (a_i - r_i) (G_i|_{D'}) \\
&= (\varphi|_{D'})^* K_D + \sum\limits_i (a_i - r_i) (G_i|_{D'})
\endalign
$$
(recall $K_D \sim 0$). Since $a_i - r_i \geq 0$ (Lemma 2.3), 
this means that $D$ has at most canonical singularities. 
\quad \qed
\enddemo

\definition{2.5}\ \ 
Let $D \in |-K_X|$ be a general member as in Lemma 2.4, and let 
$$
h := g|_D : D \to V := g(D). 
$$
$h$ is a projective bimeromorphic morphism 
such that 

\flushpar
(2.5.1) \qquad \qquad \qquad \qquad \qquad \quad \, \, \, 
$F := \text{Exc } h \subset E$. 

\flushpar
Let $F = \sum\limits_j F_{1, j} + \sum\limits_k F_{2,k}$ 
be the irreducible decomposition, where 
$$
\dim F_{1,j} =1 \, \text{ and } \, \dim F_{2,k} =2. 
$$
Moreover let $F_1 := \sum\limits_j F_{1, j}$ and 
$F_2 := \sum\limits_k F_{2,k}$: 
$$
F = F_1 + F_2. 
$$
Note that each $F_{2,k}$ coincides with some $E_i$. 
\enddefinition

\proclaim{Lemma 2.6} 
\qquad \qquad \qquad \qquad \quad 
$\text{Bs } |D| = F_2$. 
\endproclaim

\demo{Proof}\ \ 
Let $L := -K_X|_D$ and consider 
the exact sequence 
$$
0 \longrightarrow \Cal O_X \longrightarrow 
\Cal O_X(-K_X) \longrightarrow \Cal O_D(L) 
\longrightarrow 0. 
$$
Since $R^1 g_*\Cal O_X =0$ [loc.cit], 

\flushpar
(2.6.1) \qquad \qquad \qquad \qquad \qquad \, \, \, \, 
$\text{Bs } |L| = \text{Bs } |-K_X|$. 

\flushpar
Let $x \in F_1 - F_2$ be an arbitrary point, 
then we can find an effective Cartier divisor $M_x$ 
on $D$ such that 

\flushpar
(2.6.2) \qquad \qquad \qquad \qquad \qquad \, \, \, \, 
$M_x \cap F = \{x\}$, \, \, and 

\flushpar
(2.6.3) \, \, \, 
$(M_x \, . \, F_{1,j}) =1$ \, for any $F_{1,j}$ containing $x$. 

\flushpar
In particular, $L - M_x$ is $h\text{-nef}$. 
Again by [loc.cit], 
$R^1 h_* \Cal O_D (L - M_x) =0$, and 
$$ 
h_* \Cal O_D (L) \longrightarrow h_* \Cal O_{M_x} (L)
$$
is surjective. 
Thus $x \not\in \text{Bs }|L|$. 
From this and (2.6.1), we have 
$x \not\in \text{Bs }|-K_X|$,  
and we conclude 
$\text{Bs } |-K_X| = F_2$. 
\quad \qed
\enddemo

\proclaim{2.7}\ \ Proof of Theorem 2.1. 
\endproclaim 

\flushpar
(2.7.0) \, \, \, Let $H$ be a general very ample 
divisor on $D$. By shrinking $V$ if necessary, 
$H$ is decomposed as $H = H_1 + H_2$ so that 
$H_1 \cap F_1 = H_2 \cap F_2 = \emptyset$, and 
$|H_2|$ gives a projective bimeromorphic morphism 
$h_1 : D \to \widetilde{V}$ such that 
$$
\text{Exc } h_1 = F_2. 
$$
Let $F^{\circ}$ be a connected component 
of $F_2$, $D^{\circ}$ an analytic neighborhood of 
$F^{\circ}$ in $D$, $V^{\circ} := h_1(D^{\circ})$, and 
$$
h^{\circ} := h_1|_{D^{\circ}} : D^{\circ} \to V^{\circ}. 
$$
Then 

\flushpar
(2.7.1) \, \, \, $h^{\circ}$ is a projective bimeromorphic 
morphism which contracts a connected divisor $F^{\circ}$ 
to a point, and each irreducible component 
of $F^{\circ}$ coincides with some $E_i$. 

We denote $L|_{D^{\circ}}$ again by $L$. Then by (2.6.1) 
and Bertini's theorem, 

\flushpar
(2.7.2) \qquad \qquad \qquad \qquad \qquad \quad \, 
$\text{Sing } D^{\circ} \subset \text{Bs }|L|$. 

\flushpar
Let $h^{\circ}(F^{\circ}) =: P^{\circ} \in V^{\circ}$, let 
$L_0 \in |L|$ be a general member, $r$ a sufficiently large 
integer, and $L'$ a general hyperplane section of 
$(V^{\circ}, P^{\circ})$. Let 

\flushpar
(2.7.3) \qquad \qquad \qquad \qquad \qquad 
$L_r := L_0 + r h_0^*L' \in |L|$. 

\flushpar
Take a projective bimeromorphic morphism 
$\psi : D' \to D^{\circ}$ from a smooth 3-fold 
$D'$, with a simple normal crossing divisor 
$\sum\limits_i G_i'$, such that 
$$
\cases
\psi^* L_r = \sum\limits_i r_i' G_i', \\ 
K_{D'} = \psi^* K_{D^{\circ}} + \sum\limits_i a_i' G_i' 
\sim \sum\limits_i a_i' G_i', \, \text{ and} \\
- \psi^* L - \sum\limits_i \delta_i' G_i' \text{ is } 
(h^{\circ} \circ \psi)\text{-ample}
\endcases
\tag 2.7.4
$$
for some $r_i', a_i' \in \Bbb Z$ 
with $r_i \geq 0$, and $\delta_i' \in \Bbb Q_{>0}$ 
($0 < \delta_i' < 1)$.  

Note the following three things: 

\flushpar
(2.7.5) \qquad \qquad \quad \, \, \, \, \, 
$\psi (\bigcup\limits_i G_i') \subset 
\text{Bs }|L| \cup \text{Sing } D^{\circ} 
= \text{Bs }|L| = F^{\circ}$, 

\flushpar
(2.7.6) \qquad \qquad \qquad \qquad \qquad \qquad \, \, \, \, 
$a_i' \geq 0 \, \, \, (\forall i)$, 

\flushpar
and 

\flushpar
(2.7.7) \qquad \qquad \qquad \qquad \qquad \qquad \, \, \, \, 
$r_i' \geq r \, \, \, (\forall i)$. 

\flushpar
In fact, (2.7.5) comes from Lemma 2.6 (with (2.7.2)), 
(2.7.6) from Lemma 2.4, and (2.7.7) from the definition of 
$L_r \in |L|$ (2.7.3). 

Since $a_i'\text{'s}$ do not depend on $r$, we can shrink 
the minimum 
$$
c' = c'(r) := \min \dfrac{a_i' + 1 - \delta_i'}{r_i'}
$$
as little as we need: 
 
\flushpar
(2.7.8) \qquad \qquad \qquad \qquad \quad \, \, \, \, \, 
$0 < c' \ll 1$ \quad (as $r \gg 0$). 

As before, we may assume that the minimum $c'$ is attained 
exactly for a single $i=1$, say, 
and let 
$$
A' := \sum\limits_{i \geq 2}(-c'r_i' + a_i' - \delta_i') G_i', 
\, \text{ and } \, B' := G_1'. 
$$
Then 
$$
C' := \psi^* L - K_{D'} + (A' - B') 
= (1-c') \biggl( \psi^* L - 
\sum\limits_i \dfrac{\delta_i'}{1-c'} G_i' \biggr)
$$
is $(h^{\circ} \circ \psi)\text{-ample}$ 
(2.7.4), (2.7.8). Thus 
$$
R^1 (h^{\circ} \circ \psi)_* \Cal O_{D'} 
(\psi^*L + \ulcorner A' \urcorner - B') =0 
$$
[loc.cit], and we get the surjection 
$$
(h^{\circ} \circ \psi)_* \Cal O_{D'} 
(\psi^*L + \ulcorner A' \urcorner) 
= h^{\circ}_* \Cal O_{D^{\circ}}(L) 
\longrightarrow 
H^0(B', \Cal O_{B'}(\psi^*L + \ulcorner A' \urcorner)). 
$$
Since $\psi(B') \subset E$ (2.7.5) and 
$\nu_i^* \Cal O_{E_i}(-K_X) \simeq \Cal O_{\Bbb P^2}(1)$ 
$(\forall i)$ (Theorem 1.1), the right-hand side does not vanish. 
Hence $\psi(B') \not\subset E$, a contradiction. 
Now the proof of Theorem 2.1 is completed. 
\quad \qed

\proclaim{Corollary 2.8} 
\qquad \qquad \qquad \quad \, \, \, \, \, 
$E_i \simeq \Bbb P^2 \, \, \, (\forall i)$. 
\endproclaim

\demo{Proof}\ \ 
By Theorem 1.1, it is enough to show that the 
normalization morphism 
$$
\Bbb P^2 \simeq \widetilde{E}_i 
\overset \nu_i \to{\longrightarrow} E_i
$$ 
is an isomorphism. Since $\big|-K_X|_{E_i} \big|$ 
is free (Theorem 2.1), the associated rational map: 
$$
\Phi : E_i \to \Bbb P^{N-1} \quad 
(N = \dim H^0 (E_i, \Cal O_{E_i}(-K_X)))
$$
is actually a morphism. Since $\Cal O_{E_i}(-K_X)$ 
is ample, 

\flushpar
(2.8.1) \qquad \qquad \qquad \qquad \qquad \qquad \quad \, \, \, \, 
$N \geq 3$. 

\flushpar
Consider the natural inclusion: 
$$
H^0 (E_i, \Cal O_{E_i}(-K_X)) \hookrightarrow 
H^0 (\widetilde{E}_i, \nu_i^* \Cal O_{E_i}(-K_X)) 
\simeq H^0 (\Bbb P^2, \Cal O_{\Bbb P^2}(1)). 
$$
Since the right-hand side has dimension just 3, 
this must be an isomorphism: 

\flushpar
(2.8.2) \qquad \qquad \qquad \qquad \qquad \qquad \quad \, \, \, \, 
$N=3$. 

\flushpar
Let us consider the composition 
$\Bbb P^2 \simeq \widetilde{E}_i 
\overset \nu_i \to{\longrightarrow} E_i
\overset \Phi \to{\longrightarrow} \Bbb P^2$. 
Since both arrows are birational morphisms, 
so is the composition, and is hence an isomorphism. 
In particular, $\Phi$ is set theoretically 
a bijection, and by Zariski Main Theorem, 
$\Phi$ is actually an isomorphism: 
$E_i \simeq \Bbb P^2$. \quad \qed 
\enddemo

\proclaim{Corollary 2.9}\ \ 
Let $X \supset E \overset{g}\to{\longrightarrow} 
Y \ni Q$ be a flipping contraction 
with the assumption (A-1). 
Then 
$$
\text{(A-2)} \, \, \, \Longleftrightarrow 
\, \, \, R^2 g_* \Cal O_X(2K_X) = 0. 
$$
\endproclaim

\demo{Proof}\ \ 
Take a general $L \in |-2K_X|$. By Theorem 2.1 we can 
choose such $L$ to be smooth. $g(L) \in |-2K_Y|$, and 
$g|_L : L \to g(L)$ gives a resolution of $g(L)$. 
Consider the exact sequence: 
$$
0 \longrightarrow 
\Cal O_X(2K_X) \longrightarrow 
\Cal O_X \longrightarrow 
\Cal O_L \longrightarrow 0
$$
Then by $R^i g_* \Cal O_X = 0$ $(i = 1,2)$, 
$$
R^2 g_* \Cal O_X(2K_X) \simeq 
R^1 (g|_L)_* \Cal O_L,  
$$
and hence the result. \quad \qed
\enddemo

\proclaim{Corollary 2.10}\ \ 
Let $X \supset E \overset{g}\to{\longrightarrow} 
Y \ni Q$ be satisfying the Assumption A, 
{\it i.e.\/} of Type (R). Then $E$ is irreducible:
$$
E \simeq \Bbb P^2. 
$$
\endproclaim

\demo{Proof}\ \ 
Assume that $E$ has more than one irreducible components. 
Let $L \in |-2K_X|$ be a general member, 
let $k := g|_L : L \to g(L)$, and consider the exact sequence: 
$$
0 \longrightarrow \Cal O_X(2K_X) 
\longrightarrow \Cal O_X 
\longrightarrow \Cal O_L \longrightarrow 0
$$
$R^1 g_* \Cal O_X = 0$ 
by Grauert--Riemenschneider vanishing theorem, 
and 
$R^2 g_* \Cal O_X(2K_X) =0$ 
by Corollary 2.9, hence we have 

\flushpar
(2.10.1) \qquad \qquad \qquad \qquad \qquad \qquad \, 
$R^1 k_* \Cal O_L =0.$

\flushpar
Let $C := L \cap E$, then this implies 

\flushpar
(2.10.2) \qquad \qquad \qquad \qquad \qquad \qquad \,  
$H^1(\Cal O_C)=0$. 

\flushpar
Since $|-2K_X|$ is also free (Theorem 2.1), 
and since $\Cal O_{E_i}(-2K_X) \simeq \Cal O_{\Bbb P^2}(2)$
(Theorem 1.1), if two of $\{E_1, \dots, E_n\}$ 
intersects with each other along a 1-dimensional 
subspace, then $L$ has at least two irreducible components
meeting at distinct two points, which contradicts 
(2.10.2). Thus 

\flushpar
(2.10.3) \, \, \, 
For any two $E_i$ and $E_j$ $(i \not= j)$, 
$E_i \cap E_j$ is at most a discrete set of points. 

\flushpar
(2.10.4) \, \, \, 
In the rest we follow Kawamata [Kaw4]'s argument. 

Take in turn a general $D \in |-K_X|$. Then by Theorem 2.1 
and by the assumption $\# \text{Sing }X < \infty$, 
$D$ is smooth. Consider $g|_D : D \to g(D)$, 
then $-K_D$ is numerically trivial relative to this 
$g|_D$. Since $(Y, Q)$ is a rational singularity 
$R^i g_* \Cal O_X = 0$ $(i = 1,2)$ and is hence 
Cohen--Macaulay, $g(D) \in |-K_Y|$ is Gorenstein 
and is in particular normal. Thus any fiber of 
$g|_D : D \to g(D)$ must be connected. 
This, (2.10.3), and Theorem 2.1 imply the irreducibility 
of $E$. 
\quad \qed
\enddemo

\vskip 5mm

\head \S 3.\ The invariant $\varepsilon_P(Z \! \supset \! S)$. 
\endhead

In this section we introduce a numerical invariant
$\varepsilon_P = \varepsilon_P(Z \! \supset \! S)$ 
for an analytic germ of a singularity $(Z,P)$ together with 
its closed subspace $S$. This section is logically independent 
from the former sections. 

\definition {Definition 3.0}\ \ 
Let $(Z,P)$ be an analytic germ of a singularity, 
and let $S \subset Z$ be an irreducible closed subspace 
passing through $P$. Assume that 

\flushpar
(3.0.1) \qquad \qquad \qquad \qquad \qquad \quad \, 
$S \cap \text{Sing }Z = \{P\}$. 

\flushpar 
Define: 

\flushpar
(3.0.2) \qquad \qquad \quad \, \, 
$\varepsilon_P = \varepsilon_P(Z \! \supset \! S) 
:= \dim_P \Cal E \! xt_{\, S}^1 
(\Omega_Z^1 \otimes \Cal O_S, \Cal O_S)$. 
\enddefinition 

\proclaim{Theorem 3.1}\ \ 
In the above, assume that $Z$ is a complete intersection 
singularity, and that $S$ is smooth. Then we can 
write down the set of defining equations of $Z$ 
inside $(\Bbb C^N,0)$ (where $N$ is the least possible 
embedding dimension of $(Z,P)$) as follows: 

$$
\align
(3.1.1) \qquad \qquad \qquad \quad \, \, \, 
S &= \, \{x_{s+1} = \dots = x_N = 0\} 
\qquad \qquad \qquad \qquad \qquad \qquad \quad \, \, \, \\
\subset Z &= \, \{f_1(x_1, \dots , x_N) = \dots 
= f_r(x_1, \dots , x_N) = 0\} \\
\subset V &:= (\Bbb C^N,0) = \{(x_1, \dots, x_N)\}, \\
f_i(x_1, \dots , x_N) 
&= \, \sum\limits_{j=s+1}^N x_j \cdot g_{i,j}(x_1, \dots , x_s) 
+ h_i (x_1, \dots , x_N), \\
\text{with } \, \, 
g_{i,j} &\in \, (x_1, \dots , x_s) \, \, \text{ and } \, \, 
h_i \in (x_{s+1}, \dots , x_N)^2. 
\endalign
$$
Then 
$$
\Cal E \! xt_{\, S}^1 (\Omega_Z^1 \otimes \Cal O_S, \Cal O_S) 
\simeq \bigoplus\limits_{i=1}^r \Cal O_S \big/ 
(g_{i,s+1}, \dots , g_{i,N})
$$
and thus 
$$
\varepsilon_P (Z \! \supset \! S) = 
\sum\limits_{i=1}^r \operatorname{length} \Cal O_S \big/ 
(g_{i,s+1}, \dots , g_{i,N}). 
$$
\endproclaim

\demo{Proof}\ \ 
Let $J$ be the ideal sheaf of $Z$ inside $V = (\Bbb C^N,0)$. 
Since $J/J^2$ is a free $\Cal O_Z\text{-module}$, 

\flushpar
(3.1.2) \qquad \qquad 
$0 \longrightarrow J/J^2 \otimes \Cal O_S 
\longrightarrow \Omega_V^1 \otimes \Cal O_S 
\longrightarrow \Omega_Z^1 \otimes \Cal O_S 
\longrightarrow 0$ 

\flushpar
is exact. Let 
$\alpha : \Cal T_V \otimes \Cal O_Z \to (J/J^2)^{\vee}$
\allowlinebreak
$:= \Cal H \! om_Z(J/J^2, \Cal O_Z)$  
be the natural $\Cal O_Z\text{-homo-}$
\allowlinebreak
morphism: 
$$
\dfrac{\partial}{\partial x_i} \otimes 1 \mapsto 
\Big( [f] := f \text{ mod } J^2 \mapsto  
\dfrac{\partial f}{\partial x_i}\Big|_Z \Big)
$$
In (3.1.2), $J/J^2 \otimes \Cal O_S$ and 
$\Omega_V^1 \otimes \Cal O_S$ are both free 
$\Cal O_S\text{-modules}$, so 
$$
\Cal T_V \otimes \Cal O_S 
\overset \alpha \otimes \Cal O_S \to{\longrightarrow} 
(J/J^2)^{\vee} \otimes \Cal O_S \longrightarrow
\Cal E \! xt_{ S}^1 (\Omega_Z^1 \otimes \Cal O_S, \Cal O_S) 
\longrightarrow 0
$$
is exact. In other words, 

\flushpar
(3.1.3) \qquad \qquad \qquad \, \, \, 
$\Cal E \! xt_{\, S}^1 (\Omega_Z^1 \otimes \Cal O_S, \Cal O_S) 
\simeq \text{Coker } ( \alpha \otimes \Cal O_S)$. 

\flushpar
Then it is straightforward to see that 

\flushpar
{\bf Fact.}\ \ Let $\varphi : J/J^2 \to \Cal O_Z$ 
be an arbitrary $\Cal O_Z\text{-homomorphism}$. 
Then $\varphi \otimes \Cal O_S 
\in \text{Im } (\alpha \otimes \Cal O_S)$ 
if and only if 
there exist $\xi_1, \dots , \xi_N \in \Cal O_S$ such that 
for all $f \in J$, 
$$
\varphi([f])\big|_S = \sum\limits_{i=1}^N \xi_i 
\dfrac{\partial f}{\partial x_i} \bigg|_S 
\, \, \, \text{ on } S. \quad \text{---}
$$

Since $J/J^2 = \bigoplus\limits_{i=1}^r \Cal O_Z \cdot [f_i]$, 
the condition above is interpreted in terms of 
the Jacobian matrix as: 
$$
\left( \matrix
\dfrac{\partial f_1}{\partial x_1}\Big|_S & 
\ldots & 
\dfrac{\partial f_1}{\partial x_N}\Big|_S \\
\vdots & \ddots & \vdots \\
\dfrac{\partial f_r}{\partial x_1}\Big|_S & 
\ldots & 
\dfrac{\partial f_r}{\partial x_N}\Big|_S 
\endmatrix \right)
\left( \matrix 
\xi_1 \\
\vdots \\
\xi_N
\endmatrix \right)
=
\left( \matrix
\varphi(f_1)|_S \\
\vdots \\
\varphi(f_r)|_S
\endmatrix \right)
\tag 3.1.4
$$
\qquad \qquad \qquad \qquad \qquad \qquad 
\qquad \qquad \qquad \qquad 
$(\exists \xi_1 , \dots , \xi_N \in \Cal O_S)$.

\flushpar
By the expression (3.1.1), 
the left-hand side can be rewritten as 
$$
\left( \matrix
0 & \ldots & 0 & g_{1, s+1} & \ldots & g_{1,N} \\
\vdots &   & \vdots & \vdots &   & \vdots \\
0 & \ldots & 0 & g_{r, s+1} & \ldots & g_{r,N} 
\endmatrix \right)
\left( \matrix 
\xi_1 \\
\vdots \\
\xi_N
\endmatrix \right) , 
\tag 3.1.5
$$
so the condition (3.1.4) is further reduced to 
$$
\left( \matrix
g_{1, s+1} & \ldots & g_{1,N} \\
\vdots &   & \vdots \\
g_{r, s+1} & \ldots & g_{r,N} 
\endmatrix \right)
\left( \matrix 
\xi_{s+1} \\
\vdots \\
\xi_N
\endmatrix \right) 
=
\left( \matrix
\varphi(f_1)|_S \\
\vdots \\
\varphi(f_r)|_S
\endmatrix \right)
\tag 3.1.6
$$
\qquad \qquad \qquad \qquad \qquad \qquad 
\qquad \qquad 
$(\exists (\xi_{s+1} , \dots , \xi_N) 
\in \Cal O_S^{\oplus (N-s)})$.

\flushpar
Hence 
$$
\text{Im } ( \alpha \otimes \Cal O_S) = 
(g_{1, s+1} , \dots , g_{1,N}) \oplus 
\dots \oplus 
(g_{r, s+1} , \dots , g_{r,N}) 
$$
as $\Cal O_S\text{-sub}$ modules of 
$$
(J/J^2)^{\vee} \otimes \Cal O_S \simeq 
\Cal O_S \oplus \dots \oplus \Cal O_S. 
$$
This means 
$$
\Cal E \! xt_{ S}^1 (\Omega_Z^1 \otimes \Cal O_S, \Cal O_S) 
\simeq \bigoplus\limits_{i=1}^r \Cal O_S \big/ 
(g_{i,s+1}, \dots , g_{i,N}), 
$$
as required. 
\quad \qed
\enddemo

\proclaim{Corollary 3.2}\ \ 
Assume furthermore that $Z$ is a hypersurface singularity, 
and that $S$ is a smooth curve. Then by a suitable biholomorphic 
change of coordinates, $S \subset Z \subset V = (\Bbb C^N,0)$ 
is described as: 
$$
\align
(3.2.1) \qquad \qquad \, \, 
S = \{ x_2 = \dots = x_N = 0\} 
\subset Z = \{ f(x_1 , \dots , x_N) =0\} \qquad & 
\qquad \qquad \\
\subset V = (\Bbb C^N,0) = \{(x_1 , \dots , x_N)\},& \\
f(x_1 , \dots , x_N) 
= x_N \cdot x_1^e + h(x_1 , \dots , x_N), 
\, \, \text{ with } \, \, 
h \in (x_2 , \dots , x_N)^2.& 
\endalign
$$
Then 

\flushpar 
\qquad \qquad \qquad \qquad \qquad \qquad \qquad \, \, \, \, \, 
$\varepsilon_P (Z \! \supset \! P) = e$. \quad \qed 
\endproclaim

\definition{Remark 3.3}\ \ 
Mori introduced the numerical invariant 
$i_P(1)$ in [Mo4] \S 2 in quite a different manner. 
The above Corollary 3.2 and [loc.cit] ((2.16) lemma) 
show however that these coincide with each other, 
when $Z$ is a terminal Gorenstein 3-fold singularity. 
\enddefinition

\proclaim{Corollary 3.4}\ \ 
Let $(Z,P) \supset S$ be as in Theorem 3.1. Then 

\flushpar
(1) \qquad \qquad \qquad \qquad \quad \, \, 
$\varepsilon_P(Z \! \supset \! S) \geq 
\operatorname{emb.codim} (Z,P)$. 

\flushpar
In particular, $\varepsilon_P(Z \! \supset \! S)=0$ if and only if 
$Z$ is smooth. 

\flushpar
(2) \, \, \, Let $D$ be a Cartier divisor of $Z$ passing through $P$ 
such that $T := D \cap S$ is again smooth. Then 
$$
\varepsilon_P(Z \! \supset \! S) \geq 
\varepsilon_P(D \! \supset \! T). 
$$
\endproclaim

\demo{Proof} \, \, Clear. 
\quad \qed 
\enddemo

For our purpose the case $(\dim Z, \dim Z) = (4, 2)$ 
is important. Particularly those with 
$\varepsilon_P (Z \supset S) = 1$ can be completely 
determined as follows:  
 
\proclaim{Corollary 3.5}\ \ 
Let $(Z, P) \supset S$ be as in Theorem 3.1. 
Assume that $(Z, P)$ is an isolated singularity, 
and
$$
\dim X = 4, \quad \dim S = 2, 
\quad \varepsilon_P(Z \supset S) = 1.
$$
Then 
$$
(Z, P) \simeq \{x_1 x_3 + x_2 x_4 + x_5^m = 0\} 
\supset S = \{x_3 = x_4 = x_5 = 0\} \quad 
(m \geq 2). 
$$
\endproclaim

\demo{Proof}\ \ 
By Corollary 3.4, the assumption 
$\varepsilon_P (Z \supset S) = 1$ implies that 
$(Z, P)$ should be a hypersurface singularity.
Let 
$$
\align
(Z, P) &\simeq \{f(x_1, ... , x_5) = 0\} 
\supset S = \{x_3 = x_4 = x_5 = 0\}, \\
f(x_1, ... , x_5) &= g_3(x_1, x_2) \cdot x_3 
+ g_4(x_1, x_2) \cdot x_4 + g_5(x_1, x_2) \cdot x_5 
+ h(x_1, ... , x_5), \\
& \qquad \qquad \qquad \qquad \qquad 
\qquad \qquad \qquad \qquad \qquad 
(h \in (x_3, x_4, x_5)^2). 
\endalign
$$
Then again by the assumption 
$\varepsilon_P (Z \supset S) = 1$ and Theorem 3.1, 
$(g_3, g_4, g_5) = (x_1, x_2)$ in $\Bbb C\{x_1, x_2\}$. 
So by a suitable biholomorphic change of the variables 
$\{x_1, x_2\}$, 
we may assume 
$$
g_3 = x_1, \quad g_4 = x_2, \quad g_5 \in (x_1, x_2): 
\, \, \, \text{arbitrary, \quad {\it i.e.\/}}
$$
$$
f(x_1, ... , x_5) = x_1 x_3 
+ x_2 x_4 + g_5(x_1, x_2) \cdot x_5 
+ h(x_1, ... , x_5). 
$$
Then it is easy to see that 
after a biholomorphic change of $\{x_3, x_4, x_5\}$, 
$$
f(x_1, ... , x_5) = x_1 x_3 + x_2 x_4 + c \cdot x_5^m 
\quad (c \in \Bbb C, \, m \geq 2), 
$$
where $c \not= 0$ by the assumption that 
$(Z, P)$ is isolated. Finally, under the above changes 
of coordinates the ideal $(x_3, x_4, x_5)$ is not changed, 
hence $E = \{x_3 = x_4 = x_5 = 0\}$. 
\quad \qed
\enddemo

\vskip 5mm

\head 
\S 4.\ The normal bundle $N_{E/X}$.
\endhead

\definition{Definition 4.0}\ \ 
Assume that $X \supset E 
\overset{g}\to{\longrightarrow} Y \ni Q$ 
satisfies the Assumption (A-1), {\it i.e.\/}
of Type either (R) or (I). 

Let $E_i$ be an irreducible component of $E$. 
$E_i \simeq \Bbb P^2$ (Corollary 2.8). Moreover $E$ is 
irreducible: $E = E_i$ 
for Type (R) (Corollary 2.10). Let $I_{E_i}$ be the ideal 
sheaf of $E_i$ in $X$. Define:
$$
N_{E_i/X} := \Cal H \! om (I_{E_i}/I_{E_i}^2, \Cal O_{E_i}). 
$$
This is a priori a locally free sheaf of rank 2 
on $E_i$. Moreover 
$$
c_1(N_{E_i}/X) = -2. 
\tag 4.0.1
$$
In fact take a general line $l$ (so that 
$\text{Sing }X \cap l = \emptyset$), 
then $(-K_X \, . \, l) = 1$ (Theorem 1.1), 
{\it i.e.\/} $c_1 (N_{l/X}) = -1$, so 
$$
c_1(N_{E_i/X}) = c_1 (N_{E_i/X}|_l) 
= c_1(N_{l/X}) - c_1(N_{l/E_i}) = -2. \quad \text{---}
$$
\enddefinition

The main result of this section is: 

\proclaim{Theorem 4.1}\ \ 
Let 
$$
X \supset E \simeq \Bbb P^2 
\overset{g}\to{\longrightarrow} Y \ni Q
$$ 
be a flipping contraction satisfying the 
Assumptions (A-1), (A-2), {\it i.e.\/} 
of Type (R). (See Corollary 2.10.) Assume 
$\operatorname{Sing} X (\cap E) \not= \emptyset$. 
Then 
$$
N_{E/X} \simeq 
\Cal O_{\Bbb P^2} \oplus \Cal O_{\Bbb P^2}(-2). 
$$
\endproclaim

The proof consists of two steps: Proposition 4.2 and 
Proposition 4.4. 

\proclaim{Proposition 4.2}\ \ (Step 1)

Let $l \subset E$ be any line. 
Then 
$$
N_{E/X} \otimes \Cal O_l \simeq 
\cases
\Cal O_{\Bbb P^1}(-1)^{\oplus 2}, \, \, \, \text{ or} \\
\Cal O_{\Bbb P^1} \oplus \Cal O_{\Bbb P^1}(-2). 
\endcases
$$
\endproclaim

\demo{Proof}\ \ 
Take $D = D_l \in |-K_X|$ 
such that $D \cap E$ coincides with the given $l$: 
$$
D \cap E = l. 
$$
Let $J_l$ be the ideal of $l$ in $D$, 
and let $J_l^{(m)}$ be the saturation of 
$J_l^m$ in $\Cal O_D$ (see Mori [Mo4]). 
$$
\cases
J_l / J_l^{(2)} \simeq (J_l / J_l^2)/ \Cal T \! or
\quad \text{ as } \Cal O_l\text{-modules, and} \\
c_1(J_l/J_l^{(2)}) = 2 - 
\sum\limits_{P \in \, l \, \cap \, \operatorname{Sing} D} 
\varepsilon_P (D \supset l) 
\endcases
\tag 4.2.0
$$
(see [loc.cit] \S 2). Let 
$$
J_l / J_l^{(2)} \simeq \Cal O_{\Bbb P^1}(a_l) \oplus 
\Cal O_{\Bbb P^1}(b_l) \, \, \, 
\text{ with } \, \, \,
a_l \leq b_l. 
\tag 4.2.1
$$

First by the assumption (A-2) 
with Corollary 2.9, 
$$
R^2 g_* \Cal O_X(2K_X) = 0. 
$$
This combined with the exact sequence; 
$$
0 \longrightarrow 
\Cal O_X (2K_X) \longrightarrow 
\Cal O_X (K_X) \longrightarrow 
\Cal O_D (K_X) \longrightarrow 0
$$
we get 
$$
R^1 (g|_D)_* \Cal O_D(K_X) = 0. 
\tag 4.2.2
$$
This, the formal function theorem, and the exact sequence
$$
\align
0 \longrightarrow 
J_l^{(m)} / J_l^{(m+1)} \otimes \Cal O_X(K_X) 
\longrightarrow 
\Cal O_D &/ J_l^{(m+1)} \otimes \Cal O_X(K_X) \\
&\longrightarrow 
\Cal O_D / J_l^{(m)} \otimes \Cal O_X(K_X) 
\longrightarrow 0
\endalign
$$
yield 
$$
H^1 (l, \, J_l / J_l^{(2)} \otimes 
\Cal O_{\Bbb P^1}(-1)) = 0. 
$$
(Recall $\Cal O_l(K_X) \simeq \Cal O_{\Bbb P^1}(-1)$.) 
Hence 
$$
\cases
a_l \geq 0 \, \, \, \text{ in (4.2.1), \, \, \, and} \\
\text{in particular } \, \, \, c_1(J_l/J_l^{(2)}) 
\geq 0. 
\endcases
\tag 4.2.3
$$

Finally, 

\flushpar
(4.2.4) \, \, \, There is a natural injection 
$$
J_l / J_l^{(2)} \hookrightarrow 
(N_{E/X})^{\vee} \otimes \Cal O_l 
$$
(which is an isomorphism if 
$l \cap \text{Sing } X = \emptyset$). 

In fact, let $\delta_{E/X} : I_E / I_E^2 \longrightarrow 
\Omega_X^1 \otimes \Cal O_E$, 
$\delta_{l/D} : J_l / J_l^2 \longrightarrow 
\Omega_D^1 \otimes \Cal O_l$ be the natural 
homomorphisms, then 
$\text{Im } \delta_{l/D} \simeq 
(\text{Im } \delta_{E/X}) \otimes \Cal O_l$, and so 
$$
\align
J_l / J_l^{(2)} &\simeq (\text{Im } \delta_{l/D})^{\vee \vee} 
\simeq ((\text{Im } \delta_{E/X}) \otimes \Cal O_l)^{\vee \vee} \\
&\hookrightarrow 
((\text{Im } \delta_{E/X})^{\vee}  \otimes \Cal O_l)^{\vee} 
\simeq (N_{E/X} \otimes \Cal O_l)^{\vee} 
\simeq N_{E/X}^{\vee} \otimes \Cal O_l.
\endalign
$$
Hence (4.2.4). By (4.0.1), (4.2.3) and (4.2.4), 
$$
(N_{E/X})^{\vee} \otimes \Cal O_l 
\simeq \Cal O_{\Bbb P^1}(1)^{\oplus 2} 
\, \, \text{ or } \, \, 
\Cal O_{\Bbb P^1} \oplus \Cal O_{\Bbb P^1}(2), 
$$
{\it i.e.\/} 
$$
N_{E/X} \otimes \Cal O_l 
\simeq \Cal O_{\Bbb P^1}(-1)^{\oplus 2} 
\, \, \text{ or } \, \, 
\Cal O_{\Bbb P^1} \oplus \Cal O_{\Bbb P^1}(-2). 
\quad \qed 
$$
\enddemo

\proclaim{Corollary 4.3}\ \ 
({\it cf.\/} \S Appendix.)

Let $X \supset E \simeq \Bbb P^2 
\overset{g}\to{\longrightarrow}Y \ni Q$ 
be of Type (R). Then 

\flushpar
(1) \, \, \, For each singular 
point $P \in \operatorname{Sing} X$, 
$$
\operatorname{emb.codim} (X, P) \leq 2. 
$$

\flushpar
(2) \, \, \, If $\operatorname{emb.codim} (X, P) = 2$, 
then there is no other singular point of $X$ 
on $E$: 
$$
\operatorname{Sing} X = \{P\}. 
$$
\endproclaim

\demo{Proof}\ \ (1) \, \, \, 
Take a general $D \in |-K_X|$ passing through 
$P$, and $l := D \cap E$, 
then since $X$ is Gorenstein {\it i.e.\/} $D$ 
is Cartier, 

\flushpar
(4.3.1) \qquad \qquad \quad \, \, \, \, 
$\text{emb.codim }(X, P) 
= \text{emb.codim }(D, P)$
$$
\align
& \qquad \qquad \qquad \qquad \qquad \qquad \qquad \quad 
\leq \varepsilon_P (D \supset l) 
\qquad \qquad \qquad 
\text{(Corollary 3.4)} \\
& \qquad \qquad \qquad \qquad \qquad \qquad \qquad \quad 
\leq \sum\limits_{P' \in \, 
l \, \cap \, \operatorname{Sing} D} 
\varepsilon_{P'} (D \supset l) \\
& \qquad \qquad \qquad \qquad \qquad \qquad \qquad \quad 
= 2 - c_1(J_l/J_l^{(2)}) 
\qquad \qquad \qquad \quad \, \, 
\text{(4.2.0)} \\
& \qquad \qquad \qquad \qquad \qquad \qquad \qquad \quad 
\leq 2.
\qquad \qquad \qquad \qquad \qquad \qquad \quad 
\text{(4.2.3)} 
\endalign
$$

\flushpar
(2) \, \, \, Assume \, $\text{emb.codim }(X, P) = 2$. Then 

\flushpar
(4.3.2) \, \, \, All the inequalities in (4.3.1) must be 
the equality. 

If there exists another singular point 
$P' \in \text{Sing }X$, then by taking 
a general $D \in |-K_X|$ so that $D \ni P, P'$, 
we have $\varepsilon_{P'} (D \supset l) = 0$ (4.3.2). 
This means that $(D, P')$ is smooth (Corollary 3.4 (2)), 
hence so is $(X, P')$, a contradiction.
\quad \qed 
\enddemo

The second step is:

\proclaim{Proposition 4.4}\ \ (Step 2)

In Proposition 4.2, 
$$
N_{E/X} \otimes \Cal O_l \simeq 
\Cal O_{\Bbb P^1} (-1)^{\oplus 2} 
$$
is impossible. 
\endproclaim

The proof requires the following theorem 
which is a generalization of Yo. Nami- 
\allowlinebreak
kawa's {\it local moduli\/} [Nam3] 
(see also Remark 4.9 below):

\proclaim{Theorem 4.5}\ \ (see also [loc.cit] \S 1)

Let $U \overset{\varphi}\to{\longrightarrow} V = 
\{U_t \overset{\varphi_t}\to{\longrightarrow} V_t\}_{t 
\in \Delta (t)}$ be a proper bimeromorphic morphism over 
$\Delta (t) = \{t \in \Bbb C \, | \, |t|<1\}$ 
between 4-dimensional normal analytic spaces 
$U$, $V$. Let 
$$
C_t := \operatorname{Exc} \varphi_t. 
$$
Assume the following conditions: 

\flushpar
(1) \, \, \, $C_t \simeq \Bbb P^1$ \, \, \, 
$(\forall t \in \Delta (t))$, 

\flushpar
(2) \, \, \, $(K_{U_t} \, . \, C_t) = 0$ \, \, \, 
$(\forall t \in \Delta (t))$, 

\flushpar
(3) \, \, \, $U$ has only isolated terminal complete intersection 
singularities such that 

\quad 
$\emptyset \not= \operatorname{Sing} U \subset C_0$. 

\flushpar
Then for $t \not= 0$, 
$$
N_{C_t/U_t} \simeq \Cal O_{\Bbb P^1} \oplus \Cal O_{\Bbb P^1}(-2)
\, \, \text{ or } \, \, 
\Cal O_{\Bbb P^1}(1) \oplus \Cal O_{\Bbb P^1}(-3). 
$$
\endproclaim

\demo{Proof} 
Assume 
$N_{C_t/U_t} \simeq \Cal O_{\Bbb P^1}(-1)^{\oplus 2}$ 
to derive a contradiction. (See (5.1.2).)

Since $R^i \varphi_* \Cal O_U = 0$ $(i \geq 1)$, by 
Koll\'ar--Mori [KoMo] we can proceed the deformation theory 
of the morphism $\varphi : U \to V$ formulated by Z.~Ran 
[Ra1]. (For details see \S 6, especially 6.2 below.) 
Since $R^2 \varphi_* T_U = 
\text{Ext}^2_{\Cal O_U} (\Omega^1_U, \, \Cal O_U) = 0$, 

\flushpar
(4.5.0) \, \, \, The natural holomorphic map 
$$
\operatorname{Def} U \longrightarrow 
\prod\limits_{P \in \operatorname{Sing} U} 
\operatorname{Def} (U, P)
$$
is surjective (see 6.2). 

Let 
$$
\overline{C^0} := \bigcup\limits_{t \in \Delta (t)} C_t 
= \operatorname{Exc} \varphi.
\tag 4.5.1
$$
(4.5.2) \, \, \, 
Pick $P \in \text{Sing } U$ arbitrarily, and write locally 
$$
(U, P) \simeq \{f_1(x_1, ... , x_N) = 
... = f_{N-4}(x_1, ... , x_N) = 0\} \supset 
\overline{C^0} = \{x_3 = ... = x_N = 0\} 
$$
in irredundant form {\it i.e.\/} 
$\text{emb.codim }(U, P) = N-4$, 
where $N = 5$ or $6$ (Corollary 4.3). 

\flushpar
(4.5.3) \, \, \, 
Consider a local deformation $\Delta (s) \to \text{Def }(U, P)$
of $(U, P)$ by 
$$
s \mapsto \{f_1(x) + s = f_2(x) = ... = f_N(x) = 0\}.
$$ 

\flushpar
(4.5.4) \, \, \, Then by (4.5.0) this yields a deformation 
of $U$: \, $\Delta (s) \to \text{Def }U$. 
This furthermore produces a deformation of $\varphi$: 
$$
\Cal U \overset{\psi}\to{\longrightarrow} \Cal V
= \{\Cal U^s \overset{\psi^s}\to{\longrightarrow} 
\Cal V^s\}_{s \in \Delta(s)} 
$$
($\psi^0 = \varphi$, $\Cal U^0 = U$, $\Cal V^0 = V$). 

Let us consider the {\it relative Hilbert scheme\/}: 
$$
H := \operatorname{Hilb}_{\Cal U/\Cal V/\Delta(s), \, 
\big[\overline{C^0}\big]} 
\overset{\lambda}\to{\longrightarrow} \Delta(s) 
$$
(see Koll\'ar--Miyaoka--Mori [KoMiMo2,3]).

\vskip 2mm 

\flushpar
{\bf Claim.}\ \ $H$ is irreducible and 
$\lambda$ is an isomorphism.

\vskip 2mm

\flushpar
{\it Proof.\/}\ \ Let us consider first 
$$
H_1 := \operatorname{Hilb}_{\Cal U/\Cal V/\Delta(s), \, 
[C_t]} \overset{\mu}\to{\longrightarrow} \Delta(s). 
$$
$\mu^{-1}(0)$ just parametrizes 
$\{C_t\}_{t \in \Delta(t)}$. Since 
$(K_U \, . \, C_t) = 0$ $(\forall t \in \Delta(t))$ 
(assumption (2) of the theorem), 
$$
\dim H_{1, \, [C_t]} \geq 
\dim U + (-K_U \, . \, C_t) 
- \dim \operatorname{Aut} \Bbb P^1 
+ \dim \Delta (s) = 2 
\tag 4.5.5
$$
(Theorem 1.2). Moreover since we assumed that 
$N_{C_t/U} \simeq \Cal O_{\Bbb P^1} (-1)^{\oplus 2} 
\oplus \Cal O_{\Bbb P^1}$ $(t \not= 0)$, 
$H^1(N_{C_t/U}) = 0$. Hence the equality holds in (4.5.5) 
for those $[C_t]$'s ($t \not= 0$). In more precise

\flushpar
(4.5.6) \, \, \, There exists a closed analytic 
subset $A \subsetneq H_1$ (may contain a whole irreducible 
component of $H_1$) such that 
$$
\cases
A \cap \mu^{-1}(0) \subset \{[C_0]\}, \\
H_1 - A \text{ is an irreducible smooth surface, and} \\
\mu|_{H_1 - A} : H_1 - A \to \Delta (s) 
\text{ has connected fibers.} 
\endcases
$$

We are going to prove that $H_1$ is irreducible. 
So suppose on the contrary that $H_1$ is reducible: 
$$
H_1 = H_1^{\operatorname{pr}} \cup H_1' \cup ...
$$
where $H_1^{\text{pr}}$ is the irreducible 
component containing $H_1 - A$ (4.5.6) as an open 
dense subset, and $H_1'$ one another. 
Then $H_1' \cap \mu^{-1}(0) = \{[C_0]\}$ 
by (4.5.6), so each $H_1'$ must dominate $\Delta (s)$. 
By the upper-semi-continuity of fiber dimensions, 
$H_1$ is purely of dimension 1. This is however impossible, 
since for $[C^s] \in H_1'$ lying over 
$s \in \Delta (s) - \{0\}$, we have 
$\dim H_{1, \, [C^s]} \geq 2$ by the same reason as in 
(4.5.5), a contradiction. 

Hence $H_1$ must be irreducible. By this, together with 
(4.5.6), we get

$$
\cases
H_1 \text{ is an irreducible surface,} \\ 
\operatorname{Sing} H_1 \cap \mu^{-1}(0) 
\subset \{[C_0]\}, \text{ and} \\
H_1 \overset{\mu}\to{\longrightarrow} \Delta (s) 
\text{ has connected fibers.} 
\endcases
\tag 4.5.7
$$

Let $h$ and $p$ be the universal family over $H_1$ and 
the natural projection, respectively;
$$
\CD
\Cal H_1 @>{p}>> \Cal U \\
@V{h}VV @VVV \\
H_1 @>>{\mu}> \Delta(s)
\endCD
$$
and let 
$$
\overline{C^s} := p(h^{-1} \mu^{-1}(s)). 
$$
(Note that $\overline{C^s}$ coincides with the given 
$\overline{C^0}$ (4.5.1) when $s = 0$.)
Then $[\overline{C^s}] \in H = 
\operatorname{Hilb}_{\Cal U/\Cal V/\Delta(s), \, 
\big[\overline{C^0}\big]}$ and we can define a morphism 
$$
H_1 \longrightarrow H \quad (\text{over } \Delta(s))
$$
by sending the class of a curve ($\simeq \Bbb P^1$)
lying over $\mu^{-1}(s)$ to $[\overline{C^s}]$. 
Since $H_1$ was an irreducible surface (4.5.7), 
$H$ must be an irreducible curve. 
Moreover, $H$ birationally dominates 
$\Delta (s)$ through $\lambda$ 
(4.5.7). Hence $H$ is smooth and $\lambda$ is 
an isomorphism. 

\vskip 2mm

\qquad \qquad \qquad \qquad \qquad \qquad \qquad 
\qquad \qquad \qquad \qquad 
{\it (Proof of the Claim completed.)\/}

\vskip 2mm 

\flushpar
(4.5.8) \, \, \, 
{\it Proof of Theorem 4.5 continued.\/}

By the above Claim, 
$\{\overline{C_s}\}_{s \in \Delta (s)}$ 
forms a flat family where $\overline{C^0}$
is also reduced. Recall that we gave a deformation $\Cal U$ 
of $(U, P)$ (4.5.3) as follows:
$$
\align
(\Cal U, P) &\simeq \{f_1(x_1, ... , x_N) + s 
= f_2(x_1, ... , x_N) = \, ... \, = f_{N-4}(x_1, ... , x_N) = 0\} \\
&\supset \overline{C^0} = \{x_3 = \, ... \, = x_N = s = 0\}. 
\endalign
$$
Let $I_s$ be the ideal of $\overline{C^s}$ in 
$\Cal O_{(\Cal U, \, P)} 
= \Bbb C\{x_1, ... , x_N, s\} / (f_1 + s, f_2, ... , f_{N-4})$ 
(or rather its lift to $\Bbb C\{x_1, ... , x_N, s\}$). 
Then $I_0 = (x_3, \, ... \, , x_N)$ by the Claim, and hence 
$$
I_s = (x_3 + s \cdot e_3(x, s), \,  ... \, , x_N + s \cdot e_N(x, s)) 
\quad (\exists e_i(x, s) \in \Bbb C\{x, s\}). 
$$
By the condition $f_1(x) + s \in I_s$; 
$$
\align
f_1(x) + s = \xi_3(x, s) (x_3 + s \cdot e_3 (x, s)) 
+ ... &+ \xi_{N-1}(x, s) (x_{N-1} + s \cdot e_{N-1} (x, s)) \\
&+ \xi_N(x, s) (x_N + s \cdot e_N(x, s)) \\
& \qquad \qquad \qquad \quad 
(\exists \xi_i(x, s) \in \Bbb C \{x, s\})
\endalign
$$
for at least one $i$, say $i = N$, 
$e_N \cdot \xi_N \in \Bbb C\{x,s\}^{\times}$. 
In particular $\{f_1(x)= 0\} \subset (\Bbb C^N, 0)$ defines  
a smooth germ, which contradicts the asssumption that 
$(U, P) = \{f_1(x) = ... = f_{N-4}(x) = 0\}$ 
is an irredundant expression. Hence we are done.
\quad \qed
\enddemo

\flushpar
{\bf 4.6.}\ \ {\it Proof of Theorem 4.5\/} $\Longrightarrow$
{\it Proposition 4.4.\/}

Let $D_0 \in |-K_X|$ be a member passing through a 
singular point, say $P$, of $X$, and let $l_0 := D_0 \cap E$. 
Now those $l$'s such that 
$N_{E/X} \otimes \Cal O_l \simeq 
\Cal O_{\Bbb P^1}(-1)^{\oplus 2}$ consist 
of a Zariski open subset, say $W$, 
of $(\Bbb P^2)^{\vee}$ (see [OSS]). Assume 
$W \not= \emptyset$, to derive a contradiction.

If we take a smooth conic $B$ in $E$ tangent to 
$l_0$ general enough (so that $B \not\ni P$), then 

\flushpar
(4.6.1) \, \, \, 
We can take a family $\{D_t\}_{t \in \Delta}$ of 
members of $|-K_X|$ tangent to $B$ such that 
for each $\l_t := D_t \cap E$ \, $(t \not= 0)$, 
$$
l_t \cap \operatorname{Sing} X = \emptyset, \, \, \, 
\text{ and } \, \, \, [l_t] \in W.
$$

This provides a simultaneous contraction 
$$
U \overset{\varphi}\to{\longrightarrow} V 
= \{D_t \overset{g|_{D_t}}\to{\longrightarrow} 
g(D_t)\}_{t \in \Delta}
$$
satisfying the assumptions of Theorem 4.4. 
Hence
$N_{l_t/D_t} \simeq N_{E/X} \otimes \Cal O_{l_t}$
should not be $\Cal O_{\Bbb P^1}(-1)^{\oplus 2}$
by Theorem 4.4, which contradicts 
$[l_t] \in W$ (4.6.1). 
\quad \qed

\vskip 2mm

\flushpar
{\bf 4.7.}\ \ {\it Proof of Theorem 4.1.\/} 

Now the vector bundle $N_{E/X}$ of rank 2 
on $E \simeq \Bbb P^2$ has the property: 
$$
N_{E/X} \otimes \Cal O_l \simeq 
\Cal O_{\Bbb P^1} \oplus 
\Cal O_{\Bbb P^1}(-2) \quad 
(\forall l: \text{ line in } E). 
$$
Hence by Van de Ven's characterization theorem 
[V] of {\it uniform bundles\/}, 
$$
N_{E/X} \simeq 
\Cal O_{\Bbb P^2} \oplus 
\Cal O_{\Bbb P^2}(-2). 
\quad \qed
$$

\vskip 2mm 

Finally we shall give some examples of families 
as in Theorem 4.5, in the case 
$$
N_{C_t/U_t} \simeq 
\Cal O_{\Bbb P^1} \oplus 
\Cal O_{\Bbb P^1}(-2).
$$

\definition{Example 4.8}\ \ ({\it cf.\/}  5.2 below.)

\flushpar
(1) \, \, \, Case $\# \text{Sing }U = 1$ \, (Yo.~Namikawa)
$$
\align
V &:= \{x_1 x_2 + x_3^3 + x_3 x_4^2 + t \cdot x_3^2 = 0\} \\
(&\subset (\Bbb C^5, 0) = \{(x_1, ... , x_4, t)\}), 
\endalign
$$
$\varphi : U \to V$ is the blow-up with the codimension 1 center: 
$$
F := \{x_1 = x_3 = 0\} \subset V.
$$
$\text{Sing }U_0$ is one point which is 
an ordinary double point. 

\vskip 2mm

\flushpar
(2) \, \, \, Case $\# \text{Sing }U = 2$ 
\, (jointly with Yo.~Namikawa)
$$
V := \{x_1 x_2 + x_3^2(t + x_3)^2 - x_4^{2n} = 0\}. 
$$
$\varphi : U \to V$ is the blow-up with the center: 
$$
F := \{x_1 = x_3(t + x_3) + x_4^n = 0\} \subset V.
$$
$\text{Sing }U_0$ consists of two points both of which are 
isomorphic to: 
$$
\{XY + Z^2 + W^n = 0\}. \quad \text{---}
$$
\enddefinition

\definition{Remark 4.9}\ \ 
Note that in Theorem 4.5 we did not assume anything about the 
singularities of the 3-fold $U_0$. What we assumed is 
just that of the 4-fold $U$ (assumption (3) there).  

One of the disadvantages in 4-dimensional birational geometry 
is that one cannot tell much about the definite structures 
of terminal singularities, even for hypersurface ones, as opposed 
to dimension 3. 
This is considered to be one of the principal reasons that 
a 4-dimensional algebraic variety seems too complicated
to handle with. As an evidence, there in fact exists an example of 
4-dimensional terminal hypersurface singularity whose general 
hyperplane-section is a `$K3${\it-singularity\/}' (M.Reid [Re4]).
(See also \S Appendix.)

Suppose thus in Theorem 4.5 that $U_0$ has only $cDV$-singularities. 
Then the conclusion of the theorem has been known as a special 
case of Namikawa's local moduli [Nam3], whose proof had been 
based much on the description of the versal deformations of 
Du Val singularities (Brieskorn [Br], {\it cf.\/} 
[Pi1]). In the general case however his method could not be 
applicable, and what we proceeded instead is the deformation 
theory for contractions $\varphi$ (see the proof). 
Actually, this methodology does not require any particular kind of 
assumptions on the defining equations of given singularities, 
but rather makes us possible to advance under enough generality. 

This idea indicates a prospect to overcome difficulties 
arising from the complexity of 4-dimensional singularities. 
\enddefinition

\vskip 5mm

\head \S 5. Widths (after M.~Reid).
\endhead

In this section we extend the notion of widths 
introduced by M.~Reid [Re1] for 3-fold 
{\it flopping contractions\/} ((5.1.5) for the definition) 
to our 4-fold flipping contractions (Definition 5.5). 

First we recall M.~Reid's {\it Pagoda\/}:

\vskip 2mm

\flushpar
{\bf 5.1.}\ \ (Reid's ``Pagoda'' [loc.cit])

Let $U_0 \to V_0$ be a proper bimeromorphic morphism 
from a 3-fold $U_0$ with only terminal singularities 
to a germ of a normal 3-fold singularity $(V_0, Q_0)$, 
satisfying:
$$
\cases
C_0 := \operatorname{Exc} \varphi_0 \simeq \Bbb P^1 
\, \, (\text{so that } \varphi_0(C_0) = \{Q_0\}), 
\, \, \text{ and} \\
(K_{U_0} \, . \, C_0) = 0.
\endcases
\tag 5.1.0
$$
This is called a {\it 3-fold flopping contraction\/}. 
Necessarily $(V_0, Q_0)$ is a terminal 
singularity. Write it simply 
$$
U_0 \supset C_0 \simeq \Bbb P^1 
\overset{\varphi_0}\to{\longrightarrow}
V_0 \ni Q_0.
\tag 5.1.1
$$

\flushpar
(5.1.2) \, \, \, If we further assume that 
$U_0$ is smooth, then 
$$
N_{C_0/U_0} \simeq 
\cases
\Cal O_{\Bbb P^1} (-1)^{\oplus 2}, \\
\Cal O_{\Bbb P^1} \oplus \Cal O_{\Bbb P^1} (-2), 
\, \, \text{ or} \\
\Cal O_{\Bbb P^1} (1) \oplus \Cal O_{\Bbb P^1} (-3).
\endcases
$$
(This is just because 
$c_1(N_{C_0/U_0}) = 2, \, \, \text{ and } \, \, 
H^1(N_{C_0/U_0}^{\vee}) = 0$.)

We call then $C_0$ a $(-1, -1)${\it-curve\/}, 
$(0, -2)${\it-curve\/}, $(1, -3)${\it-curve in\/} 
$U_0$, respectively. 

\flushpar
(5.1.3) \, \, \, From now on we specifically consider 
the case 
$$
N_{C_0/U_0} \simeq \Cal O_{\Bbb P^1} \oplus 
\Cal O_{\Bbb P^1} (-2). 
$$
({\it cf.\/} Laufer [L], Pinkham [Pi1], 
Katz--Morrison [KaMo] and Kawamata [Kaw9] for 
$(1,-3)$-curve case.) 

\flushpar
(5.1.4) \, \, \, (First step, upwards)

Blow up $U_0$ with the center $C_0$; 
$$
\varphi^{(1)} : U^{(1)} \to U_0.
$$
Then $F^{(1)} := 
\text{Exc }\varphi^{(1)} \simeq \Sigma_2$ 
(the second Hirzebruch surface 
$\Bbb P(\Cal O_{\Bbb P^1} \oplus \Cal O_{\Bbb P^1}(-2))$), 
and the negative section $C^{(1)} \subset F^{(1)}$ 
has the normal bundle either 
$$
N_{C^{(1)}/U^{(1)}} \simeq 
\Cal O_{\Bbb P^1} (-1)^{\oplus 2}, 
\, \, \text{ or } \, \, 
\Cal O_{\Bbb P^1} \oplus \Cal O_{\Bbb P^1} (-2). 
$$
In the former case we shall stop the procedure, 
while in the latter case, blow $U^{(1)}$ up again 
with the center $C^{(1)}$: 
$$
\varphi^{(2)} : U^{(2)} \to  U^{(1)}. 
$$
Then $F^{(2)} := \text{Exc }\varphi^{(2)} 
\simeq \Sigma_2$, and let $C^{(2)}$ be the negative section 
of it. Then $N_{C^{(2)}/U^{(2)}}$ should again be 
either $\Cal O_{\Bbb P^1} (-1)^{\oplus 2}$ 
or $\Cal O_{\Bbb P^1} \oplus \Cal O_{\Bbb P^1} (-2)$. 
In the latter case repeat the process, then after a 
finitely many steps; 
$$
\varphi^{(m-1)} : U^{(m-1)} \to U^{(m-2)}
\, \, \, (\text{with } F^{(m-1)} 
:= \operatorname{Exc} \varphi^{(m-1)})
$$
we eventually arrive 
$$
N_{C^{(m-1)}/U^{(m-1)}} \simeq 
\Cal O_{\Bbb P^1} (-1)^{\oplus 2}. 
$$ 

\flushpar
(5.1.5) \, \, \, ({\it Definition of widths\/} 
[loc.cit])

Define: 
$$
\operatorname{width} \varphi_0 \, \, \, 
(\text{or } \, \operatorname{width}_{U_0} C_0) 
:= m \, \, (\geq 2). 
$$

\flushpar
(5.1.6) \, \, \, (Second step, the roof)

Blow $U^{(m-1)}$ up once again with the center 
$C^{(m-1)}$: $\varphi^{(m)} : U^{(m)} \to U^{(m-1)}$. 
Then $F^{(m)} 
:= \text{Exc }\varphi^{(m)} 
\simeq \Bbb P^1 \times \Bbb P^1$ 
({\it i.e.\/} no more $\Sigma_2$). 
Since $N_{C^{(m-1)}/U^{(m-1)}}^{\vee} 
\simeq \Cal O_{\Bbb P^1} (1)^{\oplus 2}$ 
is generated by global sections, 
$U^{(m)}$ has the contraction 
to the opposite direction, namely, there exists 
a proper bimeromorphic morphism 
$$
\varphi^{(m) \, +} : U^{(m)} \to U^{(m-1) \, +}
$$
with $F^{(m) \, +} := \text{Exc }
\varphi^{(m) \, +} = F^{(m)}$, 
$\varphi^{(m) \, +} (F^{(m) \, +}) \simeq \Bbb P^1$, 
and the fibers of $\varphi^{(m) \, +}$ 
is not linearly equivalent in $F^{(m)}$ 
to those of $\varphi^{(m)}$. 

\vskip 2mm

\flushpar
(5.1.7) \, \, \, (Third step, downwards)

Then the proper transform $F^{(m-1) \, +}$ 
of $\text{Exc } \varphi^{(m-1)}$ to $U^{(m-1) \, +}$ 
is still isomorphic to $\Sigma_2$, and is contractible 
to $\Bbb P^1$ (induction on $m$). Contract it down: 
$$
\varphi^{(m-1) \, +} : U^{(m-1) \, +} \to  
U^{(m-2) \, +}. 
$$
Consider the proper transform $F^{(m-2) \, +}$ 
of $F^{(m-2)}$ to $U^{(m-2) \, +}$, and so on. 

After $m$ steps of these contraction procedures, 
we arrive 
$$
\varphi^{(1) \, +} : U^{(1) \, +} \to U_0^+,  
$$
where $U_0^+$ is an neighborhood of 
$C_0^+ \simeq \Bbb P^1$, with $N_{C_0^+/U_0^+} 
\simeq \Cal O_{\Bbb P^1} \oplus \Cal O_{\Bbb P^1} (-2)$. 

Finally this is contractible to the original basement 
$(V_0, Q_0)$, since everything in the above procedures is 
defined over $(V_0, Q_0)$: 

\flushpar
(5.1.1)$^+$ \qquad \qquad \qquad \qquad \, 
$U_0^+ \supset C_0^+ \simeq \Bbb P^1 
\overset{\varphi_0^+}\to{\longrightarrow}
V_0 \ni Q_0$. 

We call this $\varphi_0^+$ the {\it flop\/} 
of $\varphi_0$. Also we call this whole operation 
producing $\varphi_0^+$ out of $\varphi_0$ 
the {\it flop\/}. \quad ---

\vskip 2mm

The following observation is due also to J.~Koll\'ar 
[Kol1]: 

\proclaim{Fact 5.2}\ \ 
Let $U_0 \supset C_0 \simeq \Bbb P^1 
\overset{\varphi_0}\to{\longrightarrow}
V_0 \ni Q_0$ be as above, with the assumption (5.1.3). 
Let $m = \text{width } \varphi_0$. 
Then 
$$
(V_0, Q_0) \simeq \{xy + z^2 - w^{2m} = 0\} (\subset 
(\Bbb C^4, 0) = \{(x, y, z, w)\}), 
$$
and $\varphi_0$, $\varphi_0^+$ is the blow-up of 
$U_0$, $U_0^+$ with the codimension 1 center 
$$
\{x = z + w^m = 0\}, \quad 
\{x = z - w^m = 0\}, 
$$
respectively. 
\endproclaim

\definition{Remark 5.3}\ \ 
It is convenient to define the width also for 
$(-1, -1)$-curves. Namely, for a contraction 
$U_0 \supset C_0 \simeq \Bbb P^1 
\overset{\varphi_0}\to{\longrightarrow}
V_0 \ni Q_0$ with $N_{C_0/U_0} \simeq 
\Cal O_{\Bbb P^1}(-1)^{\oplus 2}$, 
$$
\operatorname{width} \varphi_0 := 1. 
$$
\enddefinition

This invariant width has an interpretation 
from quite different point of view: 

\proclaim{Theorem 5.4}\ \ 
(Yo.~Namikawa [Nam3])

\flushpar
(1) \, \, \, 
Let $U_0 \supset C_0 \simeq \Bbb P^1 
\overset{\varphi_0}\to{\longrightarrow}
V_0 \ni Q_0$ be a contraction as in 
5.1. Then there exists an unique number $m$, 
and a 1-parameter deformation  
$$
U \overset{\varphi}\to{\longrightarrow} V 
= \{U_t \overset{\varphi_t}\to{\longrightarrow}
V_t\}_{t \in \Delta}
$$
of $U_0 \overset{\varphi_0}\to{\longrightarrow} V_0$
such that for $t \not= 0$, $\varphi_t$ is a 
contraction of $m$ disjoint $(-1, -1)$-curves. 

\flushpar
(2) \, \, \, Moreover, if 
$N_{C_0/U_0} \simeq 
\Cal O_{\Bbb P^1}(-1)^{\oplus 2}$ or 
$\Cal O_{\Bbb P^1} \oplus \Cal O_{\Bbb P^1}(-2)$, 
then this number $m$ coincides with the width 
of $\varphi_0$ (5.1.5):
$$
m = \operatorname{width} \varphi_0. 
\quad \qed 
$$
\endproclaim

Now we are ready to define the width also for our 
original 4-fold flipping contractions. 

\definition{Definition 5.5}\ \ 
(widths for dimension 4)

Let $X \supset E \simeq \Bbb P^2 
\overset{g}\to{\longrightarrow} Y \ni Q$ 
be a flipping contraction with the Assumption A, 
{\it i.e.\/} of Type (R). 
Take a general smooth $D \in |-K_X|$, 
then by Theorem 4.1 (or Proposition 4.2 and 4.4) 
$$
g|_D : D \supset l \to g(D) \ni Q 
$$ 
(where $l := D \cap E$) gives a contraction 
of the $(0, -2)$-curve $l$ (or $(-1, -1)$-curve 
only if $X$ is smooth [Kaw4]). 
Define 
$$
\operatorname{width} g := \operatorname{width} (g|_D) 
\, \text{ in the sense of } (5.1.5). 
$$
\enddefinition

The following is a paraphrase of Proposition 4.4 
in terms of widths:

\proclaim{Corollary 5.6}\ \ (Characterization of 
the Kawamata flip [Kaw4])

Let  $X \supset E \simeq \Bbb P^2 
\overset{g}\to{\longrightarrow} Y \ni Q$ 
be a flipping contraction of Type (R). 
Then 
$$
X \text{ is smooth} \, \, \, 
\Longleftrightarrow \, \, \, 
\operatorname{width} g = 1. 
\quad \qed
$$
\endproclaim

In \S 4 we determined the normal bundle 
$N_{E/X}$ (Theorem 4.1). This provides 
us the information on the embedding 
$E \subset X$ in the first order level, while the 
width actually measures that in higher orders, and 
these supplement to each other toward the investigation 
of the contraction $g$.

\vskip 5mm 

\head \S 6.\ Deformations of contractions $g : X \to Y$. 
\endhead

This section is devoted to the study of deformations of 
$X \supset E \simeq \Bbb P^2 \overset{g}\to{\longrightarrow} 
Y \ni Q$ of Type (R). 
The main result of this section is: 

\proclaim{Theorem 6.1}\ \ (Globalizability of local deformations 
of $\operatorname{Sing} X$ for Type (R)) 

Assume that $X \supset E \simeq \Bbb P^2 
\overset{g}\to{\longrightarrow} Y \ni Q$ is a flipping contraction 
of Type (R). Let 
$\{P_1 , \dots , P_r\} := 
\operatorname{Sing} X (\not= \emptyset)$, 
let $U_i$ be a sufficiently small Stein open neighborhood of 
$P_i$ in $X$, and let 
$\{U_{i, \, t}\}_{t \in \Delta}$ be any small 
deformation of $(U_i, P_i)$ $(i = 1 , \dots , r)$. 
Then there exists a 1-parameter deformation 
$$
\Cal X \to \Cal Y = \{X_t 
\overset{g_t}\to{\longrightarrow} Y_t\}_{t \in \Delta} 
$$ 
($(X_0 \overset{g_0}\to{\longrightarrow} Y_0) \, = \, 
(X \overset{g}\to{\longrightarrow} Y)$) 
satisfying the following conditions:

\flushpar
(1) \, \, \, For sufficiently small Stein open 
neighborhood $\Cal U_i$ of $P_i$ in $\Cal X$, 
$\Cal U_i \to \Delta$ coincides with 
the given deformation $\{U_{i, \, t}\}_{t \in \Delta}$ 
$(i=1 , \dots , r)$, and 

\flushpar
(2) \, \, \, $g_t$ is again a flipping contraction 
of Type (R) $(\forall t)$ (where $E_t := 
\operatorname{Exc} g_t$ maynot be irreducible).
\endproclaim

\flushpar
{\bf 6.2.} \quad 
Since $R^i g_* \Cal O_X = 0$ $(i=1,2)$, 
if we give a deformation $\Cal X \to \Delta$ of $X$, 
then by Koll\'ar--Mori ([KoMo] 11.4), there exists 
a deformation $\Cal Y \to \Delta$ of $Y$ and a proper bimeromorphic 
morphism $\Cal X \to \Cal Y$ over $\Delta$ which extends 
$g : X \to Y$. 

Let $\Cal T_X := \Cal H \! om_X (\Omega^1_X, \Cal O_X)$ be the 
tangent sheaf of $X$. Note that since 
$X$ is a sufficiently small analytic neighborhood of $E$, 
$$
\text{Ext}^i_X (\Omega^1_X, \Cal O_X) = 
H^i (\Bbb R g_* \Bbb R \Cal H \! om_X(\Omega^1_X, \Cal O_X)), 
$$
which is the abutment of the spectral sequence:

\flushpar
(6.2.1) \qquad \, \, \, 
$E_2^{p,q} = R^p g_* \Cal E \! xt^q_X (\Omega^1_X, \Cal O_X) 
\Longrightarrow 
E^{p+q} = \text{Ext}^{p+q}_X (\Omega^1_X, \Cal O_X)$. 

\flushpar
Let us look at the edge sequence:

\flushpar
(6.2.2) \, \, \, 
$0 \longrightarrow R^1 g_* \Cal T_X 
\longrightarrow \text{Ext}^1_X (\Omega^1_X, \Cal O_X) 
\overset{\alpha}\to{\longrightarrow} 
g_* \Cal E \! xt^1_X (\Omega^1_X, \Cal O_X)$ 

\flushpar
\qquad \qquad \qquad \qquad \qquad 
\qquad \qquad \qquad \qquad \qquad 
$\longrightarrow R^2 g_* \Cal T_X 
\longrightarrow \text{Ext}^2_X (\Omega^1_X, \Cal O_X)$ 

By Z.~Ran [Ra1] (after the smooth case 
of Kodaira [Kod] and Horikawa [Ho1,2,3]), 
(6.2.2) describes the deformation of $g : X \to Y$ 
({\it cf.\/} [Sc], [Fri], [Nam1,2,3], [G], [NS]). 

Let $\{P_1 , \dots , P_r\}$ be as in Theorem 6.1. 
Toward the existence of a global deformation $\Cal X$ of 
$X$ as in Theorem 6.1, it is sufficient to prove: 

\flushpar
(6.2.3) \, \, \, The vanishing of the 
{\it obstruction space\/}: 
$\text{Ext}^2_X (\Omega^1_X, \Cal O_X) =0$, 
which implies the smoothness of the {\it global deformation 
space\/} $\text{Def }X$, and 

\flushpar
(6.2.4) \, \, \, The surjectivity of 
$$
\text{Ext}^1_X (\Omega^1_X, \Cal O_X) 
\overset{\alpha}\to{\longrightarrow} 
g_* \Cal E \! xt^1_X (\Omega^1_X, \Cal O_X),  
$$
which (together with (6.2.3)) implies the surjectivity of 
the natural morphism:
$$
\operatorname{Def} X \to 
\prod\limits_{i=1}^r \operatorname{Def} (X,P_i). 
$$ 

First we shall give the following: 

\proclaim{Lemma 6.3}\ \ 
The exact sequence (6.2.1) can be extended as 
\endproclaim

\flushpar
\qquad \qquad \quad \, \, \, 
$R^2 g_* \Cal T_X 
\longrightarrow \text{Ext}^2_X (\Omega^1_X, \Cal O_X) 
\longrightarrow g_* \Cal E \! xt^2_X (\Omega^1_X, \Cal O_X) 
\longrightarrow 0$. 

\demo{Proof}\ \ 
Let us consider the spectral sequence (6.2.2). 
If $q \geq 1$, then 
$\text{Supp } \Cal E \! xt^q_X (\Omega^1_X,$ 
$\Cal O_X) \subset \text{Sing }X$, 
so by the assumption $\dim (\text{Sing }X \cap E) \leq 0$, 

\flushpar
(6.3.1) \qquad \qquad \qquad \qquad \quad \, \, 
$E_2^{p,q} =0$ \, \, \, $(\forall p \geq 1, \forall q \geq 1)$. 

\flushpar
On the other hand, since $g$ has at most a 2-dimensional fiber, 

\flushpar
(6.3.2) \qquad \qquad \qquad \qquad \qquad \quad \, 
$E_2^{p,0} =0 \, \, \, (\forall p \geq 3)$. 

\flushpar
Let $F^{\cdot}(E^r)$ be the associated filtration 
of $E^r$. Recall that the homomorphism 
$E_2^{2,0} \to E^2$ in the edge sequence is constructed as 
$$
E_2^{2,0} \twoheadrightarrow E_{\infty}^{2,0} = F^2(E^2) 
\subset E^2. 
$$
Since $E_2^{1,1}=0$ (6.3.1), 
$E_{\infty}^{1,1} = \text{Gr}^1(E^2) =0$, that is, 
$F^2(E^2) = F^1(E^2)$. Thus 
$$
\text{Coker } [E_2^{2,0} \to E^2] \simeq 
\text{Gr}^0(E^2) = E_{\infty}^{0,2} \simeq E_3^{0,2}. 
$$
Moreover $E_2^{2,1} =0$ (6.3.1), 
so $E_3^{0,2} \simeq E_2^{0,2}$: 

\flushpar
(6.3.3) \qquad \qquad \qquad \qquad \quad \, \, 
$\text{Coker } [E_2^{2,0} \to E^2] \simeq E_2^{0,2}$, 

\flushpar
which proves our lemma. 
\quad \qed 
\enddemo

\proclaim{Lemma 6.4} \qquad \qquad \qquad \quad \, \, \, \, 
$\Cal E \! xt^2_X (\Omega^1_X, \Cal O_X) =0$. 

\flushpar
Hence by Lemma 6.3, $R^2 g_* \Cal T_X 
\to \operatorname{Ext}^2_X (\Omega^1_X, \Cal O_X)$ 
is surjective. 
\endproclaim

\demo{Proof}\ \ (See also Koll\'ar [Kol2])

Let $P_i \in \text{Sing }X \cap E$ be any point. 
By the assumption, there exists a local immersion 
$(X,P_i) \subset (\Bbb C^N, 0) =: V$ 
which is of complete intersection. 
Let $J$ be the ideal of $X$ in $V$, then 
both $\Omega^1_V \otimes \Cal O_X$ and 
$J / J^2$ are free $\Cal O_{(X,P)}\text{-modules}$, and 
hence $\Cal E \! xt^1_X(J/J^2, \Cal O_X) = 
\Cal E \! xt^2_X(\Omega^1_V \otimes \Cal O_X, \Cal O_X) =0$. 
By the exact sequence 
$$  
0 \longrightarrow J / J^2 
\longrightarrow \Omega^1_V \otimes \Cal O_X 
\longrightarrow \Omega^1_X 
\longrightarrow 0
$$
we have 
$\Cal E \! xt^2_X (\Omega^1_X, \Cal O_X)=0$. 
\quad \qed
\enddemo

\proclaim{Lemma 6.5}\ \ 
For $g$ of Type (R), (6.2.3) and (6.2.4) hold. 
\endproclaim

\demo{Proof}\ \ 
By Lemma 6.3 and 6.4, (6.2.3) and (6.2.4) are reduced to 

\flushpar
(6.5.1) \qquad \qquad \qquad \qquad \qquad \qquad \, \, 
$R^2 g_* \Cal T_X =0$. 

\flushpar
By the formal function theorem, 
$$
(R^2 g_* \Cal T_X)^{\wedge} =
\varprojlim H^2(\Cal O_X / I_E^n \otimes \Cal T_X). 
$$
This, combined with; 
$$
0 \longrightarrow I_E^n / I_E^{n+1} \otimes \Cal T_X 
\longrightarrow \Cal O_X / I_E^{n+1} \otimes \Cal T_X 
\longrightarrow \Cal O_X / I_E^n \otimes \Cal T_X 
\longrightarrow 0
$$
it is enough to show

\flushpar
(6.5.2) \qquad \qquad \qquad \quad \, \, \, \, 
$H^2(I_E^n / I_E^{n+1} \otimes \Cal T_X) =0$ \, \, \, 
$(\forall n \geq 0)$. 

\flushpar
Since there are homomorphisms 
$$
I_E^n / I_E^{n+1} \to (I_E^n / I_E^{n+1})^{\vee \vee} 
\, \, \, \text{ and } \, \, \, 
S^n (N_{E/X}^{\vee}) \hookrightarrow 
(I_E^n / I_E^{n+1})^{\vee \vee} 
$$
which are isomorphisms outside the finite set of points 
$\{P_1, \dots , P_m\}= \text{Sing }X$, 
(6.5.2) is further reduced to 

\flushpar
(6.5.3) \qquad \qquad \qquad \, \, \, 
$H^2 \bigl( E, \, S^n (N_{E/X}^{\vee}) 
\otimes \Cal T_X \bigr) =0$ \, \, \, 
$(\forall n \geq 0)$. 

\flushpar
The left-hand side is Serre dual to 

\flushpar
(6.5.4) \qquad \qquad \qquad \, \, 
$H^0 \bigl( E, \, S^n (N_{E/X}) 
\otimes (\Cal T_X \otimes \Cal O_E)^{\vee} 
\otimes \omega_E \bigr)$. 

\flushpar
We claim that this certainly vanishes. 

Take a general line $l \subset E$. Since 
for any locally free $\Cal O_E\text{-module}$ 
$\Cal F$, $H^0(l, \Cal F \otimes \Cal O_l)=0$ implies 
$H^0(E, \Cal F)=0$, and since 
$$
(\Cal T_X \otimes \Cal O_E)^{\vee} \otimes \Cal O_l 
\simeq (\Omega_X^1 \otimes \Cal O_E)^{\vee \vee} 
\otimes \Cal O_l, 
$$
it is enough to prove: 

\flushpar
(6.5.5) \qquad \qquad 
$
H^0 \bigl(l, \, S^n (N_{E/X}) 
\otimes (\Omega_X^1 \otimes \Cal O_E)^{\vee \vee} 
\otimes \omega_E \otimes \Cal O_l \bigr) =0$. 

\flushpar
Since $N_{E/X} \simeq \Cal O_{\Bbb P^2} \oplus 
\Cal O_{\Bbb P^2} (-2)$ (Theorem 4.1), 
this in fact follows from the exact sequence; 
$$
0 \longrightarrow N_{E/X}^{\vee} 
\longrightarrow (\Omega_X^1 \otimes \Cal O_E)^{\vee \vee} 
\longrightarrow \Omega_E^1 
\longrightarrow 0
$$
tensorized with $S^n (N_{E/X}) \otimes 
\omega_E \otimes \Cal O_l \simeq 
S^n \bigl( \Cal O_{\Bbb P^1} \oplus 
\Cal O_{\Bbb P^1} (-2) \bigr) \otimes 
\Cal O_{\Bbb P^1} (-3)$. 
\quad \qed 
\enddemo

\proclaim{6.6}\ \ Proof of Theorem 6.1. 
\endproclaim

The remaining thing we have to prove is that 
$g_t : X_t \to Y_t$ again gives a flipping 
contraction of Type (R). 

\flushpar
(6.6.0) \, \, \, To begin with, $Y_t$ has only isolated 
rational singularities by Elkik [E], and 
hence is normal. 

Let 
$$
H := 
\operatorname{Hilb}_{\Cal X/\Cal Y/\Delta, \, [E]}
\overset{\lambda}\to{\longrightarrow} \Delta.
\tag 6.6.1
$$
We claim that $H$ is purely 1-dimensional. 
First take any line $l \subset E$ 
and consider 
$$
H_1 := \operatorname{Hilb}_{\Cal X/\Cal Y/\Delta, \, [l]}
\overset{\mu}\to{\longrightarrow} \Delta.
$$
$$
\mu^{-1}(0) \simeq (\Bbb P^2)^{\vee}, 
\tag 6.6.2
$$ 
and
$$
\dim H_1 \geq \dim X + (-K_X \, . \, l) - 
\dim \operatorname{Aut} \Bbb P^1 + \dim \Delta = 3.
\tag 6.6.3
$$
Hence for any $t \in \Delta$, $\dim \mu^{-1}(t) \geq 2$. 
By the upper-semi-continuity of fiber dimensions of $\mu$ 
and (6.6.2) we must have 
$$
\dim \mu^{-1}(t) = 2 \quad (\forall t \in \Delta). 
\tag 6.6.4
$$
In particular, $g_t$ is not an isomorphism and 
$\dim E_t \geq 2$, where $E_t := \text{Exc }g_t$. 
Again by the upper-semi-continuity of fiber dimensions applied 
in turn to $\Cal X \to \Cal Y$, 
$$
\dim E_t = 2, \quad \dim g_t(E_t) = 0. 
\tag 6.6.5
$$
(6.6.4) and (6.6.5) imply that 
$-K_{X_t}$ is $g_t\text{-ample}$. 
Finally by Corollary 2.9 and the upper-semi-continuity,  
$$
R^2 g_{t \, *} \Cal O_{X_t}(2K_{X_t}) = 0. 
\tag 6.6.6
$$ 
Since $X_t$ has again only terminal complete 
intersection singularities, $g_t$ is a flipping contraction 
of Type (R), and we are done. 
\quad \qed

For later arguments we have to look at the 
ideal structure of $E$ associated to a global 
deformation of $g$, determined by the degeneration of 
$\{E_t\}_{t \in \Delta - \{0\}}$ as $t$ tends to 0.

\definition{Definition 6.7}\ \ 
($E_0$: The ideal structure of $E \simeq \Bbb P^2$ 
under a given deformation $\Cal X \to \Cal Y$)

Let $\Cal X \to \Cal Y = 
\{X_t \overset{g_t}\to{\longrightarrow} Y_t\}_{t \in \Delta}$ 
be as in Theorem 6.1. Let $E_t := \text{Exc }g_t$ $(t \not= 0)$. 
$E_t$ is reduced and is a disjoint union of several 
$\Bbb P^2$'s. Naturally 
$$
\operatorname{Hilb}_{\Cal X/\Cal Y/\Delta, \, [E_t]} 
\overset{\sim}\to{\longrightarrow} \Delta 
\tag 6.7.1
$$
and so we may regard 
$\lambda : H \to \Delta$ (6.6.1) as 
$$
\lambda([\text{a connected component of }E_t]) 
= [E_t]. 
\tag 6.7.2
$$
In particular 
$$
\operatorname{deg} \lambda 
= \# (\text{connected components of } E_t) \, \, \, 
(t \not= 0). 
\tag 6.7.3
$$
(6.7.4) \, \, \, Let $E_0$ be a closed analytic subspace 
corresponding to $0 \in 
\operatorname{Hilb}_{\Cal X/\Cal Y/\Delta, \, [E_t]}$ 
(6.7.1). Then $\{E_t\}_{t \in \Delta}$ forms a flat family 
of closed analytic subspaces of $\Cal X$. 
Clearly 
$$
\operatorname{red} E_0 = E (\simeq \Bbb P^2). 
\tag 6.7.5
$$
 
\flushpar
(6.7.6) \, \, \, From now on we always distinguish 
$E$ and $E_0$. Note that $E_0 = E_0^{\rho}$ 
depends on the choice of deformations $\Cal X \to \Cal Y$ 
of $g : X \to Y$ {\it i.e.\/} $\rho : \Delta \to \text{Def }X$. 
(For example for the trivial deformation; 
$\Cal X = X \times \Delta$, $E_0 = E$.) Let 
$I_{E_t} = I_{E_t}^{\rho}$ be the associated ideal sheaf 
in $\Cal O_{\Cal X}$ (or in $\Cal O_{X_t}$) (also for $t = 0$). 
Define {\it the universal ideal of\/} $E$ by 
$$
I_E^* = I_{E^*} := 
\bigcap\limits_{\rho : \Delta \to \operatorname{Def} X} 
I_{E_0}^{\rho}, 
\tag 6.7.7
$$
and $E^*$ the associated closed subspace. 
Denote 

\flushpar
(6.7.8) \qquad \quad \, \, \, \,  
$\text{mult }E_0 := 
\text{ the multiplicity of } E_0 
\text{ at a general point}$ 

\qquad \qquad \qquad \qquad \qquad \qquad 
\qquad \qquad \qquad \quad \, \, 
(or of $E_0$ as a 2-cycle). 
\enddefinition

Then the following is the key observation 
toward the local classification of singularities
on $X$ (\S 7): 

\proclaim{Proposition 6.8}\ \ 

\flushpar
(1) \qquad \qquad \qquad \qquad \qquad \qquad 
\, \, \, \, \, 
$\operatorname{mult} E_0 = \operatorname{deg} \lambda$. 

\flushpar
(2) \, \, \, $E_0$ has no embedded components. 
\endproclaim

\demo{Proof}\ \ (1) \, \, \, 
$\text{mult }E_0 \geq \text{deg }\lambda$ is clear. 
Let $r := \text{deg }\lambda$. 
Take a sufficiently large integer $d$, and a general 
$\Cal L \in |-d K_{\Cal X}|$. Let $L_t := \Cal L|_{X_t}$. 
Then since for any connected component 
$E_{t, i} \simeq \Bbb P^2$ of $E_t$, 
$\Cal O_{E_{t, i}}(-K_{X_t}) \simeq \Cal O_{\Bbb P^2}(1)$ 
$(i = 1, ... , r)$ (Theorem 1.1), 
$C_{t, i} := L_t \cap E_{t, i}$ is a smooth curve of 
degree $d$ in $E_{t, i} \simeq \Bbb P^2$. 
So on $E = \text{red }E_0 \simeq \Bbb P^2$, 
$L_0 \cap E$ must be a curve of degree 
$$
d' := \dfrac{rd}{\operatorname{mult} E_0}. 
$$
If $\text{mult }E_0 > r$, then 
$d < d'$, so on $E$ we have 
$$
0 < \operatorname{deg} (-K_X|_E) = \dfrac{d}{d'} < 1, 
$$
a contradiction. Hence (1). 

\flushpar
(2) comes from (1). 
\quad \qed 
\enddemo

The following is an immediate corollary of Theorem 6.1:

\proclaim{Corollary 6.9}\ \ (Existence of a smoothing)

Let $X \supset E \simeq \Bbb P^2 
\overset{g}\to{\longrightarrow} Y \ni Q$ 
be of Type (R), with 
$$
\operatorname{width} g = m (\geq 2). 
$$
Then 

\flushpar
(1) \, \, \, There exists a 1-parameter deformation 
$\Cal X \to \Cal Y = \{X_t 
\overset{g_t}\to{\longrightarrow} Y_t\}_{t \in \Delta}$ 
such that for $t \not= 0$, 
$X_t$ is smooth. 

\flushpar
(2) \, \, \, Moreover, for any such $\Cal X \to \Cal Y$, 
$\operatorname{Exc} g_t$ consists of an $m$ 
disjoint union of $\Bbb P^2$'s each of which has the 
normal bundle $N_{E_{t, i}/X_t} \simeq 
\Cal O_{\Bbb P^2}(-1)^{\oplus 2}$. 

In particular, 
$$
\operatorname{mult} E_0 = m. 
$$
\endproclaim

\demo{Proof}\ \ 
Since $X$ is assumed to have only isolated complete intersection 
singularities, which are locally smoothable, 
this follows immediately from Theorem 6.1 and the theorem of 
Kawamata (Theorem 0.5). 
\quad \qed 
\enddemo

\vskip 5mm

\head \S 7.\ $\operatorname{Sing} X =\{P\}, \quad 
\varepsilon_P (X \supset E) =1$. 
\endhead

In this section, we shall prove the following: 

\proclaim{Theorem 7.1}\ \ 
Let $X \supset E \simeq \Bbb P^2 
\overset{g}\to{\longrightarrow} Y \ni Q$ 
be a flipping contraction of Type (R). 
Then 
$$
\# \operatorname{Sing} X =1. 
$$
So let $\{P\} = \operatorname{Sing} X$. Then 
$$
\varepsilon_P (X \supset E) =1. 
$$
\endproclaim

The proof will be divided into several steps, and will be completed 
by Proposition 7.5. 

To start with: 

\proclaim{Proposition 7.2}\ \ 
Let $P \in \operatorname{Sing} X$, then $(X, P)$ is 
a hypersurface singularity: 
$$
\operatorname{emb.codim} (X, P) =1. 
$$
\endproclaim

\demo{Proof}\ \ 
Assume:
$$
\operatorname{emb.codim} (X, P) =2, \, \, 
\text{ and hence } \, \, \operatorname{Sing} X = \{P\}, 
\tag 7.2.0
$$
to derive a contradiction (Corollary 4.3). 
Then as in (3.1.1) 
$$
(X, P) \simeq 
\{f_1(x_1, ... , x_6) = f_2(x_1, ... , x_6) = 0\} 
\tag 7.2.1
$$
\qquad \qquad \qquad \qquad \qquad \qquad 
$\supset E = \{x_3= x_4 = x_5 = x_6 = 0\}$, 
$$
\cases
f_1(x_1, ... , x_6) = 
\sum\limits_{j = 3}^6 g_{1,j}(x_1, x_2) \cdot x_j 
+ h_1(x_1, ... , x_6), \\
f_2(x_1, ... , x_6) = 
\sum\limits_{j = 3}^6 g_{2,j}(x_1, x_2) \cdot x_j 
+ h_2(x_1, ... , x_6) \\
\qquad \qquad \qquad \qquad \qquad 
(h_1, h_2 \in (x_3, x_4, x_5, x_6)^2).
\endcases
$$
(7.2.2) \, \, \, By abuse of notations, we shall denote the lifts 
of the ideals $I_{E_0, P}$, $I^*_{E, P} 
\subset \Cal O_X \simeq \Bbb C\{x_1, ... , x_6\} / (f_1, f_2)$ 
(Definition 6.7) to $\Bbb C\{x_1, ... , x_6\}$ 
by the same notations $I_{E_0,P}$ and $I^*_{E, P}$, 
respectively:
$$
f_1, f_2 \in I^*_{E, P} \subset I_{E_0,P}. 
$$

Let $A := \Bbb C\{x_1, x_2\}$ and $M := A^{\oplus 2}$ 
regarding as a free $A$-module of rank 2. 
First we claim: 

\vskip 2mm

\flushpar
{\bf Claim 1.}\ \ After a suitable permutation of 
$\{x_3, x_4, x_5, x_6\}$, 
$$
\left( \matrix
g_{1,6} \\ 
g_{2,6} \endmatrix \right) 
\in 
\left( \matrix
g_{1,3} \\ 
g_{2,3} \endmatrix \right) 
A + 
\left( \matrix
g_{1,4} \\ 
g_{2,4} \endmatrix \right) 
A + 
\left( \matrix
g_{1,5} \\ 
g_{2,5} \endmatrix \right) 
A. 
$$

\flushpar
{\it Proof.\/}\ \ If the Claim 1 is false, 
then since $I^*_{E, P}$ 
does not have embedded primary ideals (Proposition 6.8 (2)), 
and since the radical of $I^*_{E, P}$ is 
$(x_3, x_4, x_5, x_6)$, 
$$
I^*_{E, P} = (x_3, x_4, x_5, x_6),  
$$
a contradiction to Corollary 6.9 (2). Hence the Claim 1. 

By the Claim 1, after a suitable biholomorphic change of 
coordinates we may rewrite (7.2.1) as: 
$$
\cases
f_1(x_1, ... , x_6) = 
\sum\limits_{j = 3}^5 (g_{1,j}(x_1, x_2) + h_{1,j}(x_1, ... , x_6)) 
\! \cdot \! x_j +  c_1 \, x_6^{n_1}, \\
f_2(x_1, ... , x_6) = 
\sum\limits_{j = 3}^5 (g_{2,j}(x_1, x_2) + h_{2,j}(x_1, ... , x_6)) 
\! \cdot \! x_j + c_2 \, x_6^{n_2} \\
\quad  
(x_j \! \cdot \! h_{i,j} \in 
(x_3, x_4, x_5, x_6)^2, \, \, 
n_1, n_2 \geq 2, \, \, c_1, c_2 \in \{0\} \cup 
\Bbb C\{x, s\}^{\times}).
\endcases
\tag 7.2.3
$$

\vskip 2mm 

\flushpar
{\bf Claim 2.}\ \ In (7.2.3), for any permutation 
$\{x_i, x_j, x_k\}$ of $\{x_3, x_4, x_5\}$, 
$$
\left( \matrix
g_{1,k} \\ 
g_{2,k} \endmatrix \right) 
\not\in 
\left( \matrix
g_{1,i} \\ 
g_{2,i} \endmatrix \right) 
A + 
\left( \matrix
g_{1,j} \\ 
g_{2,j} \endmatrix \right) 
A. 
$$

\flushpar
{\it Proof.\/}\ \ 
Assume to the contrary that 
$$
\left( \matrix
g_{1,5} \\ 
g_{2,5} \endmatrix \right) 
\in 
\left( \matrix
g_{1,3} \\ 
g_{2,3} \endmatrix \right) 
A + 
\left( \matrix
g_{1,4} \\ 
g_{2,4} \endmatrix \right) 
A,  
$$
say. Then 
$$
\left( \dfrac{\partial f_i}{\partial x_j} 
\right)\Bigg|_{x_3 = ... = x_6 = 0} 
\sim \left( \matrix
0 & 0 & g_{1,3} & g_{1,4} & 0 & 0 \\
0 & 0 & g_{2,3} & g_{2,4} & 0 & 0
\endmatrix \right),
$$
where $\sim$ means the equivalence of matrices 
by the fundamental linear operations of rows, in the matrix ring 
$M_{2,6}(A)$. So 
$$
\operatorname{Sing} X \supset \left\{
\left| \matrix 
g_{1,3}(x_1, x_2) & g_{1,4}(x_1, x_2) \\
g_{2,3}(x_1, x_2) & g_{2,4}(x_1, x_2) 
\endmatrix \right| 
= x_3 = x_4 = x_5 = x_6 = 0 \right\}, 
$$
which is of dimension 1, a contradiction to the assumption (A-1). 
Hence the Claim 2. 

\vskip 2mm

By the Claim 2 with Proposition 6.8 (2), 
$x_3, x_4, x_5 \in I_{E, P}^*$, and 
hence 
$$
I_{E, P}^* = (x_3, x_4, x_5, x_6^k) 
\quad (\exists k \geq 2). 
\tag 7.2.4
$$

\vskip 2mm 

\flushpar
{\bf Claim 3.}\ \ $c_1, c_2 \not= 0$ and $n_1 = n_2$ 
in (7.2.3). 

\vskip 2mm 

\flushpar
{\it Proof.\/}\ \ Let 
$$
\{f_1(x) + t = f_2(x) + \alpha t = 0\}
\quad (\alpha =0, 1)
\tag 7.2.5
$$
be a deformation of $(X, P)$, $\Cal X^{(\alpha)}
\to \Cal Y^{(\alpha)}$ 
its globalization (Theorem 6.1), and 
$E_0^{(\alpha)}$ the associated ideal structure 
(Definition 6.7). Then by 
$I_E^* \subset I_{E_0}^{(\alpha)}$ 
and Proposition 6.8 (2), 
$$
I_{E_0}^{(\alpha)} = (x_3, x_4, x_5, x_6^{k(\alpha)}) 
\quad (2 \leq \exists k(\alpha) \leq k), 
\tag 7.2.6
$$
and hence 
$$
I_{E_t}^{(\alpha)} = 
(x_3 +t \cdot e_3^{(\alpha)}(x, t), \, \, \, 
x_4 +t \cdot e_4^{(\alpha)}(x, t), \, \, \, 
x_5 +t \cdot e_5^{(\alpha)}(x, t), 
\tag 7.2.7
$$
\qquad \qquad \qquad \qquad \qquad \qquad \qquad 
\qquad \qquad \qquad \qquad \qquad \quad 
$x_6^{k(\alpha)} + t \cdot e_6^{(\alpha)}(x, t))$

\flushpar
\qquad \qquad \qquad \qquad \qquad \qquad \qquad 
\qquad \qquad \qquad \qquad \qquad \, \, \,  
$(\exists e_i^{(\alpha)}(x, t) \in \Bbb C\{x, t\})$. 

\flushpar
Since $f_1(x) + t \in I_{E_t}^{(\alpha)}$, 
$$
f_1(x) + t = \sum\limits_{j = 3}^5 
\xi_j(x, t) \! \cdot \! (x_j + t \! \cdot \! e_j^{(\alpha)}(x, t))
\, \, \, + \, \, 
\xi_6(x, t) \! \cdot \! (x_6^{k(\alpha)} 
+ t \! \cdot \! e_6^{(\alpha)}(x, t))
$$
\qquad \qquad \qquad \qquad \qquad \qquad \qquad 
\qquad \qquad \qquad \qquad \qquad \qquad 
$(\exists \xi_j(x, t) \in \Bbb C\{x, t\})$, 

\flushpar
in particular, 
$$
\xi_3, \xi_4, \xi_5 \in (x, t) \Bbb C\{x, t\} 
\, \text{ and } \, 
\xi_6 \! \cdot \! e_6^{(\alpha)} \in \Bbb C\{x, t\}^{\times}. 
\tag 7.2.8
$$
Hence 
$$
c_1 \not= 0, \, \text{ and } \, 
k(\alpha) = n_1 \, \, \, 
(\alpha = 0, 1) \, \text{ in (7.2.3), 
(7.2.6)}. 
\tag 7.2.9
$$
By running the same argument through $f_2$ 
with $\alpha = 1$, we get $c_2 \not= 0$, $k(1) = n_2$, 
and hence the Claim 3. 

\vskip 2mm 

Let $n := n_1 = n_2$. 
By (7.2.7), (7.2.8) and (7.2.9), 
$$
I_{E_t}^{(\alpha)} = (x_3 +t \! \cdot \! e_3^{(\alpha)}(x, t), 
\, \, x_4 +t \! \cdot \! e_4^{(\alpha)}(x, t), \, \, 
x_5 +t \! \cdot \! e_5^{(\alpha)}(x, t), \, \, 
x_6^n + ct)
\tag 7.2.10
$$
\qquad \qquad \qquad \qquad \qquad \qquad \qquad 
\qquad \qquad \qquad \qquad \qquad \qquad \qquad \quad 
$(c \in \Bbb C\{x, t\}^{\times})$. 

\flushpar
Let us put $\alpha = 0$ in (7.2.5) in turn, then 
$f_2(x) \in I_{E_t}^{(0)}$, {\it i.e.\/}
there exist $\eta_j(x, t) \in \Bbb C\{x, t\}$ 
$(j = 3, ... , 6)$ such that 
$$
f_2(x) = \sum\limits_{j = 3}^5 
\eta_j(x, t) \! \cdot \! (x_j + t \cdot e_j^{(0)}(x, t))
\, \, \, + \, \, 
\eta_6(x, t) \! \cdot \! (x_6^n + ct) 
\, \, \, (\forall t \in \Delta), 
\tag 7.2.11
$$
where $c := e_6^{(0)} \in \Bbb C\{x, t\}^{\times}$ (7.2.8). 
As in (7.2.8), we have 
$$
\eta_3, \eta_4, \eta_5 \in (x, t)\Bbb C\{x, t\}, 
\, \text{ and } \, 
\eta_6 \in \Bbb C\{x, t\}^{\times},  
\tag 7.2.12
$$
by the expression of $f_2$ (7.2.3), with $m = n$ (Claim 3). 
This is however absurd, because the left-hand side of (7.2.11) 
does not depend on $t$. 
Hence the assumption (7.2.0) is false and we are done.
\quad \qed
\enddemo

\proclaim{Lemma 7.3}\ \ 

\flushpar
(1) \qquad \qquad \, $(X, P) \simeq 
\{f(x_1, ... , x_5) = 0\} \supset E 
= \{x_3= x_4 = x_5 = 0\}$, 
$$
\align
f(x_1, ... , x_5) = 
(g_3(x_1, x_2) &+ h_3(x_1, ... , x_5)) \cdot x_3 \\ 
&+ (g_4(x_1, x_2) + h_4(x_1, ... , x_5)) \cdot x_4 + x_5^k
\endalign
$$
for some $k \geq 2$, $g_i \in (x_1, x_2) \Bbb C\{x_1, x_2\}$, 
and $h_i \in \Bbb C\{x_1, ... , x_5\}$ with 
$x_i \! \cdot \! h_i \in (x_3, x_4, x_5)^2$.

\flushpar
(2) \, \, \, Moreover, 
$$
I_E^* = (x_3, x_4, x_5^k). 
$$ 
\endproclaim

\demo{Proof}\ \ 
By Proposition 7.2, $(X, P)$ is a hypersurface singularity, 
and we may write down the defining equation as: 
$$
f(x_1, ... , x_5) = g_3 (x_1, x_2) \cdot x_3 + 
g_4 (x_1, x_2) \cdot x_4 + g_5 (x_1, x_2) \cdot x_5
+ h(x_1, ... , x_5)
$$
$(h \in (x_3, x_4, x_5)^2)$. Let 
$$
\Cal G : = (g_3, g_4, g_5) \subset \Bbb C\{x_1, x_2\}. 
\tag 7.3.1
$$
Then in a similar way to Claim 1 in the proof of Proposition 7.2, 
$$
g_5 \in (g_3, g_4), \, \, \, i.e. \, \, \, 
\Cal G = (g_3, g_4) 
\tag 7.3.2
$$
(after a suitable permutation of 
$\{x_3, x_4, x_5\}$). Moreover, since 
$$
\varepsilon_P (X \supset E) = 
\dim_P \Cal E \! xt^1_E 
(\Omega^1_X \otimes \Cal O_E, \Cal O_E) = 
\operatorname{length} \Bbb C \{x_1, x_2\}/\Cal G 
< \infty 
\tag 7.3.3
$$
(Theorem 3.1), $\Cal G$ is not a principal ideal. 
Hence $I_E^* = (x_3, x_4, x_5^k)$ $(k \geq 2)$, and 
the rest is clear from the information $f \in I_E^*$ and 
(7.3.2). \quad \qed
\enddemo

\proclaim{Proposition 7.4} 
\qquad \qquad \qquad $\varepsilon_P (X \supset E) =1$. 
\endproclaim

\demo{Proof}\ \ 

\flushpar
(7.4.0) \, \, \, 
Let $\Cal G = (g_3, g_4)$, where $g_3$, $g_4$ are in the 
expression of $f$ in Lemma 7.3 (1). 

Assume that $\varepsilon_P (X \supset E) \geq 2$. 
Then $(x_1, x_2) \subsetneq \Cal G$ (7.3.3), so 
we may assume 
$$
x_2 \not\in \Cal G, \text{ say.} 
\tag 7.4.1
$$
Consider the deformation 
$$
\{f(x) + t x_2 = 0\}
\tag 7.4.2
$$
of $(X, P)$ and let $E_0$ be as usual. Then by Lemma 7.3 
(with Proposition 6.8 (2)), 
$$
I_{E_0} = (x_3, x_4, x_5^n) \quad (\exists n \geq 2). 
$$
Write the condition $f(x) + t x_2 \in I_{E_t}$: 
$$
f(x) + t x_2 = \xi_3 (x, t) \! \cdot \! 
(x_3 + t \! \cdot \! e_3(x, t)) 
+ \xi_4 (x, t) \! \cdot \! 
(x_4 + t \! \cdot \! e_4(x, t)) 
\tag 7.4.3
$$
\qquad \qquad \qquad \qquad \qquad \qquad \qquad 
\qquad \qquad \qquad \qquad \, \, \, 
$+ \xi_5 (x, t) \! \cdot \! 
(x_5^n + t \! \cdot \! e_5(x, t))$ 

\flushpar
($\exists \xi_i(x, s) \in \Bbb C\{x, t\}$), where 
$$
\xi_3, \xi_4 \in (x, t) \Bbb C\{x, t\} 
\tag 7.4.4
$$

\flushpar
as in (7.2.8), (7.2.12). 
Clearly, 

\flushpar
(7.4.5) \, \, \, 
$\xi_i(x, t) \! \cdot \! e_i(x, t)$ contains 
$c \! \cdot \! x_2$ as a monomial $(c \in \Bbb C\{x, t\}^{\times})$
for at least one $i \in \{3, 4, 5\}$. 

\flushpar
(7.4.6) \, \, \, Assume in (7.4.5) $i = 3$ or $4$, 
say $i = 4$. Then $e_4 \in \Bbb C\{x, t\}^{\times}$. 

\flushpar
Rewrite (7.4.3): 
$$
f(x) + t x_2 = \xi_3(x, t) \! \cdot \! 
(x_3 + t \! \cdot \! e_3(x, t)) 
+ (c_1 x_2 + \xi_4'(x, t)) \! \cdot \! 
(x_4 + c_2 t) 
\tag 7.4.7
$$
\qquad \qquad \qquad \qquad \qquad \qquad \qquad \qquad 
\qquad \qquad \quad \, \, \, \, 
$+ \xi_5 (x, t) \! \cdot \! 
(x_5^n + t \! \cdot \! e_5(x, t))$

\flushpar
\qquad \qquad \qquad \qquad \qquad \qquad \qquad \qquad 
\qquad \qquad \qquad \qquad \qquad 
($c_1, c_2 \in \Bbb C\{x, t\}^{\times}$). 

\flushpar
Put $t = 0$, then $f(x)$ contains 
$c_1 c_2 \! \cdot \! x_2 x_4$ 
as a monomial, which contradicts the 
assumption (7.4.1) and the expression of $f$ 
(Lemma 7.3 (1)).  

So in (7.4.5) we must have $i = 5$: 
$$
\xi_5(x, t) \! \cdot \! e_5(x, t) = c x_2. 
\tag 7.4.8
$$
Again by the expression of $f$ (Lemma 7.3 (1)) and 
(7.4.3), $\xi_5 \in 
\Bbb C\{x, t\}^{\times}$, $e_5 = c' x_2$ $(c' \in 
\Bbb C\{x, t\}^{\times})$.  
Hence 
$$
I_{E_t} = (x_3 + t \! \cdot \! e_3(x, t), \, \, 
x_4 + t \! \cdot \! e_4(x, t), \, \, 
x_5^n + c' \! \cdot \! t x_2). 
$$
In particular, 

\flushpar
(7.4.9) \, \, \, 
$\text{mult }E_0 = n \geq 2$, while 
$E_t$ is irreducible $(t \not= 0)$, which 
contradicts Proposition 6.8 (1). 

Hence $\varepsilon_P (X \supset E)$ must be 1. 
\quad \qed
\enddemo

\proclaim{Proposition 7.5}\ \ 
\qquad \qquad \qquad 
$\# \operatorname{Sing} X = 1$. 
\endproclaim 

\demo{Proof}\ \ 
Assume $P, P' \in \text{Sing }X$ $(P \not= P')$, say. 
By Proposition 7.4 with Corollary 3.5, 

\flushpar
(7.5.1) \qquad \quad \, 
$(X, P) \simeq \{f(x_1, ... , x_5) = 0\} 
\supset E = \{x_3 = x_4 = x_5 = 0\}$, 

\qquad \qquad \, \, \, \, \, 
$(X, P') \simeq \{f'(y_1, ... , y_5) = 0\} \supset 
E = \{y_3 = y_4 = y_5 = 0\}$,

\qquad \qquad \qquad \qquad \qquad 
$f(x_1, ... , x_5) = x_1 x_3 + x_2 x_4 + x_5^n$, 

\qquad \qquad \qquad \qquad \quad \, \, \, \, \, 
$f'(y_1, ... , y_5) = y_1 y_3 + y_2 y_4 + y_5^{n'}$. 

\flushpar
Consider deformations of $(X, P)$ and of $(X, P')$: 
$$
\{f(x) + t = 0\}, \quad \{f'(y) + t^{n'} = 0\}. 
\tag 7.5.2
$$
These are patched together to give a global deformation 
$\Cal X \to \Cal Y = \{X_t 
\overset{g_t}\to{\longrightarrow} Y_t\}_{t \in \Delta}$ 
(Theorem 6.1), and let $\Cal E$ be the exceptional locus of 
it: 
$$
\Cal E = \bigcup\limits_{t \in \Delta} E_t. 
$$
In a similar way to the proofs of Propositions 7.2 and 7.3, 
$$
(E_t, P) = \{x_3 + t \! \cdot \! e_3(x, t) = 
x_4 + t \! \cdot \! e_4(x, t) = x_5^n + t = 0\}, 
\tag 7.5.2
$$
in particular 

\flushpar
(7.5.3) \, \, \, $\#$(connected components of $E_t$) $= n$ 
$(t \not= 0)$, and 
$$
(E_t, P') = \{y_3 + t \! \cdot \! e'_3(x, t) = 
y_4 + t \! \cdot \! e'_4(x, t) = \prod\limits_{i \in I} 
(y_5 + \zeta^i_{n'} \! \cdot \! t) = 0\} 
\tag 7.5.4
$$
where $\zeta_{n'} \in \Bbb C$ is a primitive $n'$-th roots 
of unity, and $I \subset \{1, ... , n'\}$ with $\# I = n$ 
(7.5.3). 

\flushpar
So 

\flushpar
(7.5.5) \, \, \, 
$\Cal E$ is smooth in an neighborhood of $P$ (7.5.2), 

\flushpar
while $\Cal E$ is a union of $n \geq 2$ irreducible 
components meeting at the whole 
$(E, P') = \{y_3 = y_4 = y_5 = 0\}$ (7.5.3), in particular 
$$
\operatorname{Sing} \Cal E \cap U' = E \cap U' 
\quad (\exists U' \ni P': \text{ analytic open set of } 
\Cal X).
\tag 7.5.6
$$ 
(7.5.5) and (7.5.6) contradict to each other, since 
$\text{Sing }X$ is a Zariski closed subset of $E$. Hence 
$\# \text{Sing }X = 1$. 
\quad \qed
\enddemo

Now Theorem 7.1 follows from Proposition 7.4 and 7.5. 

\proclaim{Corollary 7.6}\ \ 
Let $X \supset E \simeq \Bbb P^2 
\overset{g}\to{\longrightarrow} Y \ni Q$ 
be of Type (R). By Theorem 7.1 (with Corollary 3.5), 
$\operatorname{Sing} X = \{P\}$, and 
$$
(X, P) \simeq \{x_1 x_3 + x_2 x_4 + x_5^m = 0\}. 
$$
Then 
$$
\operatorname{width} g = m. 
$$
\endproclaim 

\demo{Proof}\ \ 
Take a smoothing 
$\{x_1 x_3 + x_2 x_4 + x_5^m + t = 0\}$, 
then 
$$
\# (\text{connected components of } E_t) = m \, \, \, 
(t \not= 0). 
$$
The result follows from this, Proposition 6.8 (1), 
and Corollary 6.9 (2). 
\quad \qed
\enddemo

\head \S 8.\ Existence of flips 
--- La Torre Pendente. 
\endhead

In this section we shall prove the existence of flips 
for Type (R) contractions. 

\proclaim{Theorem 8.1}\ \ 
Let $X \supset E \simeq \Bbb P^2 
\overset{g}\to{\longrightarrow} Y \ni Q$ 
be a flipping contraction of Type (R). 
Then the flip $g^+$ exists. 
\endproclaim 

\vskip 2mm

\flushpar
{\bf 8.2.}\ \ According to the results of \S 7 
(Theorem 7.1 and Corollary 7.6), 
$$
\operatorname{Sing} X = \{P\}, 
\tag 8.2.1
$$
$$
(X, P) \simeq \{x_1 x_3 + x_2 x_4 + x_5^m = 0\} 
\supset E = \{x_3 = x_4 = x_5 = 0\}, 
\tag 8.2.2
$$
where
$$
m = \operatorname{width} g. 
\tag 8.2.3
$$

\flushpar
(8.2.4) \, \, \, Let $\varphi : \overline{X} \to X$ 
be the blow-up of $X$ with the center $E$. 
Let 
$$
F := \operatorname{Exc} \varphi, 
$$
which is a projective 3-fold.  

Sine $N_{E/X} \simeq \Cal O_{\Bbb P^2} \oplus 
\Cal O_{\Bbb P^2} (-2)$, 
there exists a birational map between 
$F$ and $\Bbb P := \Bbb P_{\Bbb P^2} (\Cal O_{\Bbb P^2} \oplus 
\Cal O_{\Bbb P^2} (-2))$ compatible with their projections: 
$$
\CD
\quad F \quad @. \dashrightarrow 
@. \quad \Bbb P \\
@V{\varphi|_F}VV \quad @. \quad @VV{\pi}V \\
\quad E \quad @. \overset{\sim}\to\longrightarrow 
@. \quad \Bbb P^2  
\endCD
\tag 8.2.5
$$
(8.2.6) \, \, \, Let $\overline{E} \subset F$ be the 
proper transform of the negative section of $\Bbb P$, 
{\it i.e.\/} the section corresponding to 
$\Cal O_{\Bbb P} (-1_{\Cal O_{\Bbb P^2} \oplus 
\Cal O_{\Bbb P^2} (-2)})$. 

\vskip 2mm 

\proclaim{Lemma 8.3}\ \ 
(Local description of $\varphi$ along fibers)

\flushpar
(1) \, \, \, Let $S := \varphi^{-1}(P)$ 
be the fiber of $\varphi$ at $P$, and $f$ any other 
fiber of $\varphi$ over $E - \{P\}$. Then 
$$
\varphi|_{F - S} : F - S \to E - \{P\}
$$ 
is a $\Bbb P^1$-bundle: $f \simeq \Bbb P^1$, while 
$$
S \simeq \Bbb P^2. 
$$
(A {\it jumping fiber\/} in the sense of 
Andreatta--Wi\'sniewski [AW2].)

\flushpar
(2) \, \, \, $F$ is a Cartier divisor of $\overline{X}$, 
and \, $\# \operatorname{Sing} F = 1$. 

\flushpar
Let $\{\overline{P}\} := \operatorname{Sing} F$. 
Then $\overline{P} \in S$, and 
$$
\operatorname{Sing} \overline{X} = 
\cases
\{\overline{P}\} \qquad \qquad (\text{if } m \geq 3), \\
\emptyset \qquad \qquad \quad \, \, (\text{if } m = 2).  
\endcases
$$
Moreover, 
$$
(\overline{X}, \overline{P}) \simeq 
\{y_1 y_3 + y_2 y_4 + y_5^{m-1} = 0\} 
\supset F = \{y_5 = 0\}. 
$$
\endproclaim

\demo{Proof}\ \ Straightforward. See also 
[loc.cit], and Beltrametti [Be]. \quad \qed 
\enddemo

\proclaim{Lemma 8.4}\ \ (Global description of $\varphi$ 
along $F$ --1)

\flushpar
(1) \, \, \, For any line $l \subset E$ 
such that $l \not\ni P$, 
$$
Z_l := \varphi^{-1}(l) \simeq \Sigma_2 
\overset{\varphi|_{Z_l}}\to{\longrightarrow}
l \simeq \Bbb P^1. 
$$

\flushpar
(2) \, \, \, Let $C_l := \overline{E} \cap Z_l$ 
for such $l$. Then 
$$
(-K_{\overline{X}} \, . \, C_l) = 1. 
$$
\endproclaim

\demo{Proof}\ \ 
(1) is clear, since $N_{E/X} \simeq 
\Cal O_{\Bbb P^2} \oplus \Cal O_{\Bbb P^2}(-2)$ 
(Theorem 4.1). 

\flushpar
(2) \, \, \, Denote simply 
$C_l$, $Z_l$ by $C$, $Z$, respectively. 
Then $\overline{X}$ and $F$ are 
smooth along $C$ (Lemma 8.3), and 
$$
C = Z \cap \overline{E} \quad 
\text{(meeting transversally).}
$$
Hence we have the exact sequence: 
$$
0 \longrightarrow 
N_{C/F} \longrightarrow 
N_{C/\overline{X}} \longrightarrow 
N_{F/\overline{X}} \otimes \Cal O_C 
\longrightarrow 0
$$
with 
$$
N_{C/F} \simeq \Cal O_{\Bbb P^1}(-2) 
\oplus \Cal O_{\Bbb P^1}(1). 
$$
From this, and 
$$
\align
\qquad \qquad \qquad \quad 
N_{F/\overline{X}} \otimes \Cal O_C &\simeq 
\Cal O_F (-1_{N_{E/X}^{\vee}}) \otimes \Cal O_C 
\qquad \qquad (8.2.5) \qquad \qquad \qquad \quad 
\, \, \, \, \, \\
&\simeq \Cal O_{Z} 
(-1_{N_{E/X}^{\vee} \otimes \Cal O_{l}}) 
\otimes \Cal O_C \\
&\simeq \Cal O_{\Bbb P(\Cal O_{\Bbb P^1} \oplus 
\Cal O_{\Bbb P^1}(-2))} 
(-1_{\Cal O_{\Bbb P^1} \oplus \Cal O_{\Bbb P^1}(2)}) 
\otimes \Cal O_C \\
&\simeq \Cal O_{\Sigma_2}(-C-2f) \otimes 
\Cal O_C \\
&\simeq \Cal O_{\Bbb P^1}, 
\endalign
$$
we have 
$$
c_1(N_{C/\overline{X}}) = c_1(N_{C/F}) + 
c_1(N_{F/\overline{X}} \otimes \Cal O_C)
= -1, 
$$
{\it i.e.\/} $(-K_{\overline{X}} \, . \, C) = 1$. 
Hence (2). 
\quad \qed 
\enddemo

\proclaim{Proposition 8.5}\ \ (Global description of $\varphi$ 
along $F$ --2)
$$
\overline{E} \cap S = \{P\}, 
\, \, \, \text{ and } \, \, \, 
\varphi|_{\overline{E}} : \overline{E} \to E \simeq 
\Bbb P^2 \, \, \, \text{ is an isomorphism}. 
$$
$$
C_l \text{ is a line in } \overline{E} \simeq \Bbb P^2. 
$$
\endproclaim

\demo{Proof}\ \ 
\flushpar
(8.5.0) \, \, \, 
Assume $\dim (\overline{E} \cap S) = 1$, 
to get a contradiction. 

\vskip 2mm

\flushpar
{\bf Step 1.}\ \ 
First by the construction (8.2.5) and (8.2.6), 
$$
\varphi|_{\overline{E}} : \overline{E} \to E 
\simeq \Bbb P^2
$$
 is a birational morphism and 
$$
(\overline{E} \, . \, f)_F = 1. 
\tag 8.5.1
$$
Since $f$'s and lines $n$ on $S \simeq \Bbb P^2$ 
belong to the same irreducible component of 
$\text{Hilb}_{\overline{X}}$ (and also of 
$\text{Hilb}_F$) (see Andreatta--Wi\'sniewski 
[AW2]), 
$$
(\overline{E} \, . \, n)_F = 1 \quad 
(\forall n \subset S: \text{ line}). 
\tag 8.5.2
$$
Thus 

\flushpar
(8.5.3) \, \, \, 
$n_0 := \overline{E} \cap S$ is purely 1-dimensional and 
is a reduced line in $S$. 

\flushpar
Since $F$ is smooth outside $\{\overline{P}\}$ 
(Lemma 8.3 (1)), 

\flushpar
(8.5.4) \, \, \, $\overline{E}$ is smooth 
({\it resp.\/} $\overline{E} - \{\overline{P}\}$ is smooth) 
if $\overline{E} \not\ni \overline{P}$ ({\it resp.\/} 
if $\overline{E} \ni \overline{P}$). 

\vskip 2mm

\flushpar
{\bf Step 2.}\ \ 
Let $l_{\lambda} \subset E$ be any line 
passing through $P$ $(\lambda \in \Bbb P^1)$, and 
let 
$$
m_{\lambda} := (\overline{E} \cap 
\varphi^{-1} (l_{\lambda} - \{P\}))^{\text{---}}
\tag 8.5.5
$$
where ${}^{\text{---}}$ denotes the closure in 
$F$ (or in $\overline{X}$). 
Then 
$$
\varphi|_{m_{\lambda}} : m_{\lambda} \to l_{\lambda} 
\simeq \Bbb P^1
$$
is a birational morphism, and is hence an isomorphism: 
$$
m_{\lambda} \simeq \Bbb P^1. 
\tag 8.5.6
$$

\vskip 2mm

\flushpar
{\bf Step 3.}\ \ 

\flushpar
(8.5.7) \, \, \, Let $\mu : H \simeq \Bbb P^1 
\to \operatorname{Hilb}_{\overline{X}}$ 
be a morphism defined by $\lambda \mapsto 
[m_{\lambda}]$, then this is set theoretically an injection.

\flushpar
(8.5.8) \, \, \, Let $h : \Cal H \to H$ be
the family over $H$ induced from $\text{Hilb}_{\overline{X}}$, 
with $\Cal H$ being a normal surface. 

\flushpar
Since every fiber is irreducible and $h$ is a 
$\Bbb P^1$-bundle over a general point on $H$, 

\flushpar
(8.5.9) \, \, \, $h$ is in fact a $\Bbb P^1$-bundle. 
(see [Mo1], {\it cf.\/} the argument of \S 1.)

\flushpar
Let 
$$
\nu : \Cal H \to \overline{E} \subset \overline{X}
\tag 8.5.10
$$
be the projection.

\flushpar
By construction, 

\flushpar
(8.5.11) \, \, \, The composition 
$$
\Cal H \overset{\nu}\to{\longrightarrow} 
\overline{E} 
\overset{\varphi|_{\overline{E}}}\to{\longrightarrow} 
E \simeq \Bbb P^2
$$
is a birational morphism, which is an isomorphism 
over $E - \{P\}$. 

\flushpar
In particular,  
$$
\Cal H \simeq \Sigma_1. 
\tag 8.5.12
$$
Moreover, $(\varphi|_{\overline{E}})^{-1}(P) = 
n_0$ (8.5.3), so $\nu$ is finite birational, 
{\it i.e.\/} 

\flushpar
(8.5.13) \, \, \, $\nu$ is the normalization morphism, 
which is an isomorphism over $E - \{P\}$ 

\flushpar
(8.5.4). Let $C_l$ be as in Lemma 8.4 (1), then 

\flushpar
(8.5.14) \, \, \, $\nu^{-1}(C_l)$ is a $(+1)$-section, 
and $\nu^{-1}(n_0)$ is the $(-1)$-section (negative section), 
of $\Cal H \simeq \Sigma_1$. 

Summing up, we have: 

\flushpar
(8.5.15) \, \, \, Inside $F$, we have a 1-parameter family 
$\{m_{\lambda}\}_{\lambda \in \Lambda}$ $(\Lambda \simeq 
\Bbb P^1)$ such that 
$$
\align
&m_{\lambda} \simeq \Bbb P^1 \, \, \, 
(\forall \lambda \in \Lambda), \, \, \, 
(-K_{\overline{X}} \, . \, m_{\lambda}) = 0, 
\, \, \, \text{ and } \\
&\Lambda \text{ forms a whole connected component of } 
\operatorname{Hilb}_{\overline{X}}. 
\endalign
$$
({\it cf.\/} Matsuki [Ma].)

\vskip 2mm

\flushpar
{\bf Step 4.}\ \ Finally, 

\flushpar
(8.5.16) \, \, \, let us consider 
a smoothing 
$$
\{x_1 x_3 + x_2 x_4 + x_5^m + t^m = 0\} 
$$
of $(X, P)$. Then the globalization $\Cal X \to 
\Cal Y$ (Theorem 6.1) satisfies 
$$
E_t = \{x_3 = x_4 = \prod\limits_{i \in I} 
(x_5 + \zeta_m^i t) = 0\}
\tag 8.5.17
$$
$(\exists I \subset \{1, ... , m\})$ 
(after a suitable biholomorphic change of 
coordinates $\{(x_1, ... ,$ $x_5 , t)\}$), 
as in the argument of \S 7. Moreover by Corollary 6.9, 
$I = \{1 , ... , m\}$: 
$$
E_t = \{x_3 = x_4 = x_5^m + t^m = 0\}. 
\tag 8.5.18
$$

\flushpar
(8.5.19) \, \, \,  Take 
$$
E_{t, 1} := \{x_3 = x_4 = x_5 + t = 0\}, \, \, \, 
\text{ and } \, \, \, 
\Cal E_1 := \bigcup\limits_{t \in \Delta} E_{t, 1}. 
$$

\flushpar
(8.5.20) \, \, \, $\Cal E_1 \simeq \Bbb P^2 \times \Delta$, 
and $\{E_{t, 1}\}_{t \in \Delta}$ 
forms a flat family, with 
$$
E_{t, 1} = E \, \, (\text{reduced}) \text{ when } t = 0.
$$

\flushpar
(8.5.21) \, \, \, Blow up $\Cal X$ with the center 
$\Cal E_1$: 
$$
\overline{\Cal X} \to \Cal X \to \Delta
$$
By (8.5.20), 

\flushpar
(8.5.22) \, \, \, The fiber over $0 \in \Delta$ 
just coincides with $\varphi : \overline{X} \to X$. 

\flushpar
(8.5.23) \, \, \, By Corollary 6.9, 
$\text{Exc }(\overline{X_t} 
\overset{\varphi_t}\to{\longrightarrow} X_t)$ 
consists of a disjoint union of 
$(m-1)$ $\Bbb P^2$, plus one $\Bbb P^2 \times \Bbb P^1$. 
As in Kawamata [Kaw4], 

\flushpar
(8.5.24) \, \, \, $-K_{\overline{X_t}}$ is 
$(g_t \circ \varphi_t)$-ample.

\flushpar
Hence $m_{\lambda}$ cannot move outside 
$\overline{X_0}$, {\it i.e.\/} 

\flushpar
(8.5.25) \, \, \, $\Lambda \simeq \Bbb P^1$ (8.5.15) 
is a connected component also of 
$\operatorname{Hilb}_{\overline{\Cal X}/\Cal Y/\Delta}$, 

\flushpar
while by Theorem 1.2, 
$$
\align
\dim \operatorname{Hilb}_{\overline{\Cal X}/\Cal Y/\Delta, 
\, [m_{\lambda}]} 
&\geq \dim \overline{X}
+ (-K_{\overline{X}} \, . \, m_{\lambda}) - 
\dim \operatorname{Aut} \Bbb P^1 + \dim \Delta \\
&= 2, \qquad \qquad \qquad (8.5.15)
\endalign
$$
a contradiction to each other.

Hence $\dim (\overline{E} \cap S) = 0$. 
Since $\overline{E}$ and $S$ are surfaces 
inside the 3-fold $F$, 
$$
\overline{E} \cap S \subset \operatorname{Sing} F = \{P\} 
$$
(Lemma 8.3), and we have done.
\quad \qed 
\enddemo

\proclaim{Corollary 8.6} \qquad
$\rho(\overline{X}/Y) = 2$, 
\, \, \, $-K_{\overline{X}}$ is 
$(g \circ \varphi)\text{-ample}$, and 
$$
\overline{NE}(\overline{X}/Y) 
= \Bbb R_{\geq 0} [f] + \Bbb R_{\geq 0} [C_l]. 
\quad \qed
$$
\endproclaim

\proclaim{Proposition 8.7}\ \ 
The extremal ray $\Bbb R_{\geq 0} [C_l]$ 
determines a flipping contraction 
$\overline{g} : \overline{X} \to \overline{Y}$ 
of Type (R), with 
$$
\operatorname{Exc} \overline{g} = \overline{E}, 
\, \, \text{ and } \, \, 
\operatorname{width} \overline{g} = m-1. 
$$
\endproclaim

\demo{Proof}\ \ 

\flushpar
{\bf Claim 1.}\ \ 
$\Bbb R_{\geq 0}[C_l]$ gives a flipping contraction. 

In fact, if $\Bbb R_{\geq 0}[C_l]$ defines a 
divisorial contraction $\overline{g} : \overline{X} \to 
\overline{Y}$, then clearly 
$$
\operatorname{Exc} \overline{g} = F. 
$$
Since $\varphi(S) = \{P\}$, $\overline{g}$ does not 
contract $S$, and so 
$$
\dim \overline{g}(F) = 2. 
$$
Take a general fiber 
$\overline{C} \simeq \Bbb P^1$ of $\overline{g}$ over 
$\overline{g}(F)$, then  
$$
(F \, . \, \overline{C}) = -1, \, \, \text{ and } \, \, 
\dim \varphi(\overline{C}) = 1. 
$$
Hence 
$$
\align
0 < (-K_{\overline{X}} \, . \, \varphi_*(\overline{C})) 
&= (-\varphi^* K_X \, . \, \overline{C}) \\
&= (-K_{\overline{X}} + F \, . \, \overline{C}) \\
&= 1 - 1 = 0, 
\endalign
$$
a contradiction. 
So $\overline{g}$ must be a flipping contraction. 

Since $\text{Sing }\overline{X} \subset 
\{\overline{P}\}$, and $(\overline{X}, \overline{P})$ 
is again an isolated hypersurface singularity (Lemma 8.3), 
$\text{Exc }\overline{g}$ is a union of $\Bbb P^2$'s 
(Corollary 2.8). Moreover, 

\flushpar
{\bf Claim 2.} \qquad \qquad \qquad \qquad \qquad \quad \, 
$\text{Exc }\overline{g} = \overline{E}$. 

In fact, let $Z_l \simeq \Sigma_2$ 
be as in Lemma 8.3. This is a surface 
contained in the smooth locus of the 
3-fold $F$. So if we choose $l$ general enough, 
and if we assume that $\text{Exc }\overline{g}$ is reducible, 
then $Z_l \cap \text{Exc }\overline{g}$ is a reducible 
curve. Consider the birational morphism 
$\overline{g}|_{Z_l} : \Sigma_2 \simeq Z_l \to \overline{g}(Z_l)$. 
Since $C_l$ (Lemma 8.4) is the only contractible irreducible 
curve on $Z_l$, we get a contradiction. Hence 
$\text{Exc }\overline{g}$ must be irreducible and 
$\text{Exc }\overline{g} \, \cap Z_l = C_l$, thus 
necessarily $\operatorname{Exc} \overline{g} = \overline{E}$, 
{\it i.e.\/} the Claim 2. 

Finally, 

\flushpar
{\bf Claim 3.}\ \ $\overline{g}$ is of Type (R), with 
\, $\text{width }\overline{g} = \text{width }g - 1$. 

In fact, take a general $D = D_l \in |-K_X|$ with 
$l := D \cap E$ (so that $D \not\ni P$), 
and let $\overline{D}$ be its proper transform in 
$\overline{X}$, then clearly 
$\overline{D} \cap \overline{E} = C_l$ (Lemma 8.4), and 
$\overline{D} 
\overset{\varphi|_{\overline{D}}}\to{\longrightarrow} D$ 
is the blow-up of the $(0, -2)$-curve $l$ in $D$. 
So by Reid [Re1] (see also \S 5, particularly 5.1), 
$$
N_{\overline{E}/\overline{X}}|_{C_l} 
\simeq N_{C_l/\overline{D}} 
\simeq 
\cases
\Cal O_{\Bbb P^1} (-1)^{\oplus 2} \qquad \quad \, \, \, \, 
(m = 2), \\
\Cal O_{\Bbb P^1} \oplus \Cal O_{\Bbb P^1}(-2) 
\qquad (m \geq 3). 
\endcases
$$
In both cases 
$R^1 (\overline{g}|_{\overline{D}})_* 
\Cal O_{\overline{D}} (K_{\overline{X}}) = 0$, {\it i.e.\/} 
$$
R^2 \overline{g}_* \Cal O_{\overline{X}} (2K_{\overline{X}}) = 0. 
$$
Hence by Corollary 2.9, $\overline{g}$ is again of Type (R). 
Also 
$$
\operatorname{width} \overline{g} 
= \operatorname{width} (\overline{g}|_{\overline{D}})
= m - 1. \quad \qed
$$
\enddemo

\vskip 2mm

\flushpar
{\bf 8.8.}\ \ {\it Set-up of the induction argument 
toward Theorem 8.1.\/} 

We shall prove the following 9 set of statements 
(8.8.1)$_m$ through (8.8.9)$_m$ by the 
induction on $m$: 

\flushpar
(8.8.1)$_m$ \, \, \, 
The flip $g^+$ of $g$ exists for any Type (R) 
contraction $X \supset E \simeq \Bbb P^2 
\overset{g}\to{\longrightarrow} Y \ni Q$ 
with 
$$
\operatorname{width} g \leq m. 
$$

Let $g : X \to Y$ 
be any such one, with 
$$
2 \leq \operatorname{width} g \leq m. 
$$
Let 
$\varphi : \overline{X} \to X$ be the blow-up and 
$F = \text{Exc }\varphi$, as in 8.2, then as proved in 
Corollary 8.7, there exists a unique 
flipping contraction $\overline{g} : \overline{X} \to 
\overline{Y}$ from $\overline{X}$, which is of Type (R), with 
$\text{width }\overline{g} =  \text{width }g-1$. 

Let $\overline{g}^+ : \overline{X}^+ \to \overline{Y}$ 
be the flip of $\overline{g}$ (where the existence is assured 
in (8.8.1)$_m$), and let $\overline{E}^+ := \text{Exc }\overline{g}^+$. 
Let $F^+$ be the proper transform of $F$ in $\overline{X}^+$. 
Then 

\flushpar
(8.8.2)$_m$ \, \, \, 
$\overline{X}^+$ is smooth, 

\flushpar
(8.8.3)$_m$ \, \, \, $\overline{E}^+ \simeq 
\Bbb P^1$,  

\flushpar
(8.8.4)$_m$ \, \, \, There exists a birational morphism  
$\varphi^+ : \overline{X}^+ \to X^+$ with 
$\operatorname{Exc} \varphi^+ = F^+$, 
$E^+ := \varphi^+(\operatorname{Exc} \varphi^+) \simeq \Bbb P^1$, 

\flushpar
(8.8.5)$_m$ \, \, \, $X^+$ is also smoooth, 

\flushpar
(8.8.6)$_m$ \, \, \, 
$\varphi^+$ is the blow-up of $X^+$ 
with the center $E^+$, 

\flushpar
(8.8.7)$_m$ \, \, \, $\overline{E}^+$ is a section 
of the $\Bbb P^2$-bundle $\varphi^+|_{F^+} : F^+ \to E^+$, 
and 

\flushpar
(8.8.8)$_m$ \qquad \qquad \qquad 
$N_{E^+/X^+} \simeq 
\Cal O_{\Bbb P^1} \oplus \Cal O_{\Bbb P^1}(-1) 
\oplus \Cal O_{\Bbb P^1}(-2)$. 

(In particular, $X^+ \supset E^+ \to Y \ni Q$ gives 
the flip $g^+$ of $g$.) 

\flushpar
Finally, let 
$W_{\lambda} := \varphi^{-1}(l_{\lambda} - \{P\})^{\text{---}} 
\subset F$, where $\{l_{\lambda}\}_{\lambda \in \Lambda}$ 
$(\Lambda \simeq \Bbb P^1)$ is the complete family of 
lines in $E$ passing through $P$. Then 

\flushpar
(8.8.9)$_m$ \, \, \, 
The fibers of $\varphi^+|_{F^+}$ 
$(\simeq \Bbb P^2)$ are exactly the proper transforms 
$W_{\lambda}^+$ of $W_{\lambda}$'s in $\overline{X}^+$. 

\vskip 2mm

To prove Theorem 8.1, it is enough to show the following couple of 
statements: 

\flushpar
(i) \, \, \, (8.8.1)$_2$ through (8.8.9)$_2$, and 

\flushpar
(ii) \, \, \, (8.8.1)$_{m-1}$ through (8.8.9)$_{m-1}$ 
imply (8.8.1)$_m$ through (8.8.9)$_m$. 

\vskip 2mm

\flushpar
{\bf 8.9.}\ \ [{\it Proof of\/} (8.8.1)$_2$ 
--- (8.8.9)$_2$]

First, (8.8.1)$_2$ is nothing but Kawamata [Kaw4]'s 
result (Theorem 0.5). 

\flushpar
(8.9.0) \, \, \, Let $X \supset E \simeq \Bbb P^2
\overset{g}\to{\longrightarrow} Y \ni Q$ 
be any flipping contraction of Type (R), with 
$$
\operatorname{width} g = 2, 
$$
and $\overline{g} : \overline{X} 
\to \overline{Y}$, $\overline{E}$ be as in 8.8. 
Then since $\text{width } \overline{g} = 1$, 
$\overline{g}$ is a flipping contraction
of Kawamata Type [loc.cit] (Corollary 5.6). 
Thus by [loc.cit], 
$$
N_{\overline{E}/\overline{X}} \simeq 
\Cal O_{\Bbb P^2}(-1)^{\oplus 2}, \, \, 
\overline{X}^+ \text{ is smooth}, \, \, 
\text{ and } \, \, \overline{E}^+ \simeq \Bbb P^1, 
$$
{\it i.e.\/} we get (8.8.2)$_2$, (8.8.3)$_2$. 

\vskip 2mm 

\flushpar
{\bf Claim 1.}\ \ Let $\overline{\varphi} : 
\overline{\overline{X}} \to \overline{X}$ 
be the blow-up of $\overline{X}$ with the 
center $\overline{E}$, let $\overline{F} := \text{Exc } 
\overline{\varphi}$, and $F' \subset \overline{\overline{X}}$ 
the proper transform of $F \subset \overline{X}$ in 
$\overline{\overline{X}}$. Then 
$\overline{\overline{E}} := \overline{F} \cap F'$ 
is a smooth $(1,1)$-divisor in 
$\overline{F} \simeq \Bbb P^2 \times \Bbb P^1$. 
In particular, 
$\overline{\overline{E}} \simeq \Sigma_1$.

\demo{Proof}\ \ 
Consider $\varphi \circ \overline{\varphi} : 
\overline{\overline{X}} \to X$. 
The restriction to $\overline{F}$: 
$$
\overline{F} \simeq \Bbb P^2 \times \Bbb P^1 
\overset{\overline{\varphi}|_{\overline{F}}}\to{\longrightarrow}
\overline{E} 
\overset{\varphi|_{\overline{E}}}\to{\overset{\sim}\to{\longrightarrow}}
E \simeq \Bbb P^2
\tag 8.9.1
$$
coincides with the projection. Moreover, 
take a general smooth $D \in |-K_X|$ 
(so that $D \not\ni P$), and let 
$l := D \cap E$. Then 
$$
(\varphi \circ \overline{\varphi})^{-1}(l) 
= \overline{\varphi}^{-1}(C_l) \cup Z_l' 
\subset \overline{D}'
\tag 8.9.2
$$
where $Z_l'$, $\overline{D}'$ is the proper transform of 
$Z_l (\subset \overline{X})$ (Lemma 8.4 (1)), 
$D (\subset X)$ in $\overline{\overline{X}}$, 
respectively, and $C_l$ is as in Lemma 8.4 (2). 
By Reid [Re1], this forms the Pagoda of width 2 
(see 5.1), {\it i.e.\/} 
$$
\overline{\varphi}^{-1}(C_l) \simeq \Bbb P^1 \times \Bbb P^1, 
\, \, Z_l' \simeq \Sigma_2, \, \, 
\text{ and } \, \, 
\tag 8.9.3
$$
$\overline{\varphi}^{-1}(C_l) \cap Z_l'$ is an irreducible 
$(1,1)$-divisor of $\overline{\varphi}^{-1}(C_l)$. 
The Claim 1 follows immediately from this. 
\enddemo

\vskip 2mm

\flushpar
{\bf Claim 2.}\ \ 
$F'$ is smooth, and 
$$
F' 
\overset{\overline{\varphi}|_{F'}}\to{\longrightarrow} F
$$ 
is a contraction of the $(-1, -1)$-curve 
$\overline{\varphi}^{-1}(\overline{P}) \simeq \Bbb P^1$. 
$S' := \overline{\varphi}^{-1}(S) \simeq \Sigma_1$. 

\demo{Proof}\ \ 
$\overline{\varphi}|_{F'}$ is the blow-up of $F$ with 
the codimension 1 center $\overline{E}$. Since 

\flushpar
(8.9.4) \qquad $\text{Sing }F = \{\overline{P}\}$, \, \, 
$(F, \overline{P}) \simeq \{y_1 y_3 + y_2 y_4 = 0\}$ 

\qquad \qquad \qquad \qquad \qquad \qquad \qquad 
$\supset \overline{E} = \{y_3 = y_4 = 0\}$, \, \, 
$\supset S = \{y_1 = y_2 = 0\}$

\flushpar
(Lemma 8.3, Proposition 8.5), this is clear. 
\enddemo

\vskip 2mm

\flushpar
{\bf Claim 3.}\ \ 
$(\varphi \circ \overline{\varphi})^{-1}|_{F'} 
: F' \to E \simeq \Bbb P^2$ is factored as 
$$
F' \overset{h}\to{\longrightarrow} 
\operatorname{Bl}_P E \simeq \Sigma_1 
\overset{\operatorname{bl}_P}\to{\longrightarrow} 
E \simeq \Bbb P^2,
$$
where $\text{bl}_P$ 
is the blow-up of $E$ with the center $P$. 
Moreover, $h$ is a $\Bbb P^1$-bundle. 

\demo{Proof}\ \ 
Let $f \simeq \Bbb P^1$ be the fiber of the 
$\Bbb P^1$-bundle 
$\varphi|_{F - S} : F - S \to E - \{P\}$ (Lemma 8.3), 
and $f'$ its proper transform in $\overline{\overline{X}}$. 
Let $H$, $H'$ be the irreducible component of 
$\text{Hilb}_{\overline{X}}$, 
$\text{Hilb}_{\overline{\overline{X}}}$, containing the 
point $[f]$, $[f']$, respectively. 
Then as in Andreatta--Wi\'sniewski [AW2], 
$H \simeq \text{Bl}_P E \simeq \Sigma_1$. Hence by Claim 2, 
$$
H' \simeq H \simeq \Sigma_1. 
\tag 8.9.5
$$
Let 
$$
\CD
\Cal H' @>{p}>> \overline{\overline{X}} \\
@V{h}VV @. \\
H' @. @. 
\endCD
\tag 8.9.6
$$
be the standard diagram of the universal family 
$h$ over $H'$, and the projection $p$. 
$$
p(\Cal H') = F', 
\tag 8.9.7
$$
and $p$ is set theoretically an injection 
[loc.cit]. Since $F'$ is smooth (Claim 2), 
$p$ gives an isomorphism onto $F'$: 
$$
p : \Cal H' \overset{\sim}\to{\longrightarrow} 
F'.
\tag 8.9.8
$$
Thus 

\flushpar
(8.9.9) \, \, \, $h$ gives a $\Bbb P^1$-bundle 
over $\Sigma_1$: 
$$
F' \overset{\sim}\to{\longrightarrow} \Cal H 
\overset{h}\to{\longrightarrow} \Sigma_1. 
$$
Hence the Claim 3. 
\enddemo

\vskip 2mm 

\flushpar
(8.9.10) \, \, \, Let us consider 
the composite surjective morphism:
$$
F' \longrightarrow \Sigma_1 \longrightarrow B := \Bbb P^1. 
$$
Clearly $\rho(F'/B) = 2$. 

\flushpar
(8.9.11) \, \, \, Let $\overline{\varphi}^+ : 
\overline{\overline{X}} \to \overline{X}^+$ 
be the contraction of 
$\overline{F} \simeq \Bbb P^2 \times \Bbb P^1$ to 
$\Bbb P^1$ so that $\overline{X}^+$ gives the flip 
of $\overline{g} : \overline{X} \to \overline{Y}$ 
(8.9.0) [Kaw4]. Let 
$$
\overline{E}^+ := \overline{\varphi}^+(\overline{F}) 
\simeq \Bbb P^1, \, \, \, F^+ := \overline{\varphi}^+(F'). 
$$

\vskip 2mm

\flushpar
{\bf Claim 4.}\ \ $-K_{F'}$ is relatively ample over 
$B$, and the extremal ray, other than that determines 
$h$, coincides with 
$$
\overline{\varphi}^+|_{F'} 
: F' \to F^+,  
$$
where $\overline{\varphi}^+$, $F^+$ are as in (8.9.11). 
Moreover, this is a divisorial contraction, with the exceptional 
divisor 
$$
\operatorname{Exc} (\overline{\varphi}^+|_{F'}) = 
\overline{\overline{E}} (\simeq \Sigma_1), 
\, \, \, 
\overline{\varphi}^+(\overline{\overline{E}}) 
= \overline{E}^+ \simeq \Bbb P^1. 
$$
$F^+$ is smooth, and $F^+ \to B = \Bbb P^1$ is 
a $\Bbb P^2$-bundle.

\demo{Proof}\ \ 
Let $W_{\lambda} \subset F$ be as in 8.8, and 
$W_{\lambda}'$ its proper transform in $F'$. Then 
$W_{\lambda} \simeq W_{\lambda}' \simeq \Sigma_1$. 

\flushpar
(8.9.12) \, \, \, Let $m_{\lambda}$, 
$m_{\lambda}'$ be the negative section of 
$W_{\lambda}$, $W_{\lambda}'$, respectively. 
Clearly 

\flushpar
(8.9.13) \, \, \, $\{m_{\lambda}\}_{\lambda \in \Lambda}$ 
forms the complete family of lines in $\overline{E} 
\simeq \Bbb P^2$ passing through $P$. Hence by the Claim 2, 

\flushpar
(8.9.14) \, \, \, $\{m_{\lambda}'\}_{\lambda \in \Lambda}$ 
forms the complete family of rulings in $\overline{\overline{E}}
\simeq \Sigma_1$ (Claim 1). 

Since $\overline{\varphi}^+ : \overline{\overline{X}} \to 
\overline{X}^+$ maps $\overline{F} \simeq \Bbb P^2 \times 
\Bbb P^1$ to $\Bbb P^1$, 

\flushpar
(8.9.15) \, \, \, $\overline{\varphi}^+|_{F'} : 
F' \to F^+$ contracts exactly $m_{\lambda}'$'s 
(Claim 1). Hence 

\flushpar
(8.9.16) \, \, \, $\overline{\varphi}^+|_{F'}$ 
factors through $F' \to B = \Bbb P^1$. 

Moreover, 
$$
(-K_{F'} \, . \, m_{\lambda}') 
= (-K_{\overline{\overline{X}}} - F' \, . \, m_{\lambda}') 
= 2 - 1 = 1, 
\tag 8.9.17
$$
in particular, $-K_{F'}$ is 
$\overline{\varphi}^+|_{F'}$-ample. Since 
$\rho(F'/B) = 2$ (8.9.10), necessarily 

\flushpar
(8.9.18) \, \, \, $\overline{\varphi}^+|_{F'}$ is the contraction 
of the extremal ray 
$\Bbb R_{\geq 0} [m_{\lambda}']$ of 
$\overline{NE} (F'/B)$. 

\flushpar
(8.9.19) \, \, \, 
This is furthermore a divisorial contraction from a smooth 3-fold 
$F'$ (Claim 2), with 
$$
\operatorname{Exc} (\overline{\varphi}^+|_{F'}) = 
\overline{\overline{E}} \simeq \Sigma_1, 
\, \, \, \overline{\varphi}^+(\overline{\overline{E}}) 
\simeq \Bbb P^1. 
$$
Hence by Mori [Mo2], $F^+ = \overline{\varphi}^+(F')$ 
is smooth. Since $\overline{\varphi}^+|_{F'}$ is 
defined over $B = \Bbb P^1$ (8.9.18), 

\flushpar
(8.9.20) \, \, \, There is a surjective morphism 
$F^+ \to B$. 

\flushpar
By construction, the fibers are exactly 
$\overline{\varphi}^+(W_{\lambda}') \simeq \Bbb P^2$ (8.9.15), 
so 

\flushpar
(8.9.21) \, \, \, $F^+$ is a $\Bbb P^2$-bundle 
over $B = \Bbb P^1$. 

Hence the Claim 4. 
\enddemo

\vskip 2mm

\flushpar
{\bf Claim 5.} \qquad \qquad \qquad \qquad \qquad \, \, \, 
$\rho(\overline{X}^+/Y) = 2$. 
$$
\overline{NE} (\overline{X}^+/Y) = 
\Bbb R_{\geq 0}[\overline{E}^+] + R^+, 
$$
where $R^+$ is the unique extremal ray. 
Let $\varphi : \overline{X}^+ \to X^+$ be the associated 
contraction. Let $\varphi^+(F^+) = : E^+$. 
Then $X^+$ is smooth, $E^+ \simeq \Bbb P^1$, and 
$$
\varphi^+|_{F^+} : F^+ \to E^+ \simeq \Bbb P^1
$$
coincides with the $\Bbb P^2$-bundle in the Claim 4. 

\demo{Proof}\ \ 
$\rho(\overline{X}^+/Y) = 2$ is clear, and
$$
(K_{\overline{X}^+} \, . \, \overline{E}^+) = 1
\tag 8.9.22
$$
[Kaw4]. So

\flushpar
(8.9.23) \, \, \, $\overline{NE}(\overline{X}^+/Y)$ 
has at most one extremal ray. 

Let $f^+$ be any line in a fiber of the $\Bbb P^2$-bundle 
$F^+ \to B$ (Claim 4). Then 
$$
(-K_{\overline{X}^+} \, . \, f^+) = 
(-K_{\overline{\overline{X}}} \, . \, m_{\lambda}' + f') 
= (-K_{\overline{X}} \, . \, m_{\lambda} + f) = 2.
\tag 8.9.24
$$
Hence 

\flushpar
(8.9.25) \, \, \, $\overline{NE}(\overline{X}^+/Y)$ 
has exactly one extremal ray (8.9.22), (8.9.24). 
Let $\varphi^+ : \overline{X}^+ \to X^+$ be the associated 
contraction. Then $\varphi|_{F^+}$ coincides with 
the $\Bbb P^2$-bundle $F^+ \to B = \Bbb P^1$. Moreover, 
$X^+$ is smooth along $E^+ := \varphi^+(F^+) = B$, 
since $N_{W_{\lambda}^+/\overline{X}^+} \simeq 
\Cal O_{\Bbb P^2} \oplus \Cal O_{\Bbb P^2}(-1)$ 
(8.9.24). Finally, $\overline{E}^+$ is clearly the section of 
$\varphi^+|_{F^+}$. Hence the Claim 5. 
\enddemo

\vskip 2mm

\flushpar
{\bf Claim 6.} \qquad \qquad \quad \, \, 
$N_{E^+/X^+} \simeq \Cal O_{\Bbb P^1} \oplus 
\Cal O_{\Bbb P^1}(-1) \oplus \Cal O_{\Bbb P^1}(-2)$. 

$X^+ \to Y$ gives the flip of $g : X \to Y$ in (8.9.0).

\demo{Proof}\ \ Take $D \in |-K_X|$ general enough. 

\flushpar
(8.9.26) \, \, \, In the above procedure 
$$
X \overset{\varphi}\to{\longleftarrow}
\overline{X} \overset{\overline{\varphi}}\to{\longleftarrow}
\overline{\overline{X}} 
\overset{\overline{\varphi}^+}\to{\longrightarrow}
\overline{X}^+ \overset{\varphi^+}\to{\longrightarrow}
X^+,
$$
let $\overline{D}$, $\overline{\overline{D}}$, 
$\overline{D}^+$, $D^+$ be the proper transforms of $D$ 
in $\overline{X}$, $\overline{\overline{X}}$, 
$\overline{X}^+$, $X^+$, respectively. 
Then 
$$
D \longleftarrow \overline{D} \longleftarrow
\overline{\overline{D}} \longrightarrow
\overline{D}^+ \longrightarrow D^+
$$
forms Reid [Re1]'s Pagoda (5.1) of width 2, in particular 
they are all smooth. 
So 
$$
N_{E^+/D^+} \simeq 
\Cal O_{\Bbb P^1} \oplus \Cal O_{\Bbb P^1}(-2).  
$$
By this, and the exact sequence: 
$$
0 \longrightarrow N_{E^+/D^+} 
\longrightarrow N_{E^+/X^+} 
\longrightarrow N_{D^+/X^+} \otimes 
\Cal O_{E^+} \longrightarrow 0
\tag 8.9.27
$$
with 
$$
(D^+ \, . \, E^+)_{X^+} = 
(-K_{\overline{\overline{X}}} + 2F' + 4 \overline{F}
\, . \, \overline{f})_{\overline{\overline{X}}} 
= 1 + 2 - 4 = -1 
\tag 8.9.28
$$
(where $\overline{f} \simeq \Bbb P^1$ is a fiber of 
$\overline{\varphi}|_{\overline{F}} 
: \overline{F} \to \overline{E}$), we get 
$$
N_{E^+/X^+} \simeq \Cal O_{\Bbb P^1} \oplus 
\Cal O_{\Bbb P^1}(-1) \oplus \Cal O_{\Bbb P^1}(-2). 
$$
Hence
$$
(K_{X^+} \, . \, E^+) = -2 - c_1(N_{E^+/X^+}) 
= 1, 
$$ 
and so $X^+ \to Y$ certainly gives the flip of $g$, 
{\it i.e. \/} the Claim 6. 
\enddemo

Now (8.8.4)$_2$ through (8.8.9)$_2$ are consequences of 
the Claim 5 and 6. Also (8.8.1)$_2$ follows from 
the existence of the flip $X^+ \to Y$ (Claim 6). 

\vskip 2mm

\flushpar
{\bf 8.10.}\ \ [Proof of (8.8.1)$_{m-1}$ --- 
(8.8.9)$_{m-1}$ $\Longrightarrow$ 
(8.8.1)$_{m}$ --- (8.8.9)$_{m}$]

Assume (8.8.1)$_{m-1}$ --- (8.8.9)$_{m-1}$. 

\flushpar
(8.10.0) \, \, \, Let $X \supset E \simeq \Bbb P^2 
\overset{g}\to{\longrightarrow} Y \ni Q$ 
be any flipping contraction of Type (R) with 
$$
\operatorname{width} g = m \geq 3. 
$$

\flushpar
(8.10.1) \, \, \, Let $\varphi : \overline{X} \to X$, 
$F$, $\overline{E}$, $\overline{g} : \overline{X}
\to \overline{Y}$ be as in 8.8. 

\flushpar
(8.10.2) \, \, \, Let $\overline{\varphi} : 
\overline{\overline{X}} \to \overline{X}$ 
be the blow-up of $\overline{X}$ with the center 
$\overline{E}$, let $\overline{F} := 
\text{Exc }\overline{\varphi}$, 
$F'$ the proper transform of $F$ in 
$\overline{\overline{X}}$, 
$\overline{E}' := \overline{F} \cap F'$, and 
$\overline{S} := \overline{\varphi}^{-1}(\overline{P}) 
\simeq \Bbb P^2$. 

\flushpar
(8.10.3) \, \, \, Since 
$\text{width }\overline{g} = m-1 \geq 2$, 
there exists a flipping contraction 
$\overline{\overline{g}} : \overline{\overline{X}} 
\to \overline{\overline{Y}}$ of Type (R) with 
$\text{width }\overline{\overline{g}} = m-2$ 
(Proposition 8.7). 

\flushpar
Let 
$\overline{\overline{g}}^+ : \overline{\overline{X}}^+
\to \overline{\overline{Y}}$ be the flip, 
where the existence is assured by the induction 
hypothesis (8.8.1)$_{m-1}$. 
Let $\overline{\overline{E}} := 
\text{Exc } \overline{\overline{g}}$, 
$\overline{\overline{\psi}} := 
\overline{\overline{g}}^{+ \, -1} \circ 
\overline{\overline{g}} : 
\overline{\overline{X}} \dashrightarrow 
\overline{\overline{X}}^+$, $\overline{F}^+ := 
\overline{\overline{\psi}}_*(\overline{F})$, and 
${F'}^+ := \overline{\overline{\psi}}_*(F')$. 

\vskip 2mm

\flushpar
{\bf Claim 1.} \qquad \qquad \qquad \quad \, \, \, 
$\overline{\overline{E}}
\cap \overline{E}' = \emptyset$, \, 
{\it i.e.\/} \, $\overline{\overline{E}}
\cap F' = \emptyset$. 

\flushpar
Hence $\overline{\overline{\psi}}|_{F'} : F' \to 
{F'}^+$ is an isomorphism: $F' \simeq {F'}^+$. 

\demo{Proof}\ \ 
For any line $l \not\ni P$ in $E$, 

$$
(\varphi \circ \overline{\varphi})^{-1}(l) 
\cap \overline{\overline{E}}
\cap \overline{E}' = \emptyset
\tag 8.10.4
$$
[Re1]. On the other hand, 
$$
\align
(\overline{X}, \overline{P}) &\simeq 
\{y_1 y_3 + y_2 y_4 + y_5^{m-1} = 0\} \\
\supset \overline{E} &= \{y_3 = y_4 = y_5 = 0\}, \, \, \, 
\supset S = \{y_1 = y_2 = y_5 = 0\}
\endalign
$$
(Lemma 8.3, Proposition 8.5). 
Hence on the blown-up $\overline{\overline{X}}$ (with the center 
$\overline{E}$), 
the proper transform $S'$ of $S$ in 
$\overline{\overline{X}}$ is away from 
$\text{Sing }\overline{\overline{X}}$, while again by 
Proposition 8.5, 
$$
\overline{\overline{E}} \cap 
\overline{\varphi}^{-1}(\overline{P}) = 
\{\overline{\overline{P}}\}. 
$$
So 
$$
(\varphi \circ \overline{\varphi})^{-1}(P) 
\cap \overline{\overline{E}} \cap \overline{E}' 
= \emptyset.
\tag 8.10.5
$$
(8.10.4) and (8.10.5) show the Claim 1. 
\enddemo

\vskip 2mm

\flushpar
(8.10.6) \, \, \, By the induction hypothesis 
(8.8.4)$_{m-1}$ through (8.8.7)$_{m-1}$, 

\flushpar
(8.10.7) \, \, \, There exists a contraction 
$$
\overline{\varphi}^+ : \overline{\overline{X}}^+ 
\to \overline{X}^+
$$
with $\text{Exc }\overline{\varphi}^+ = \overline{F}^+$, 
$\overline{E}^+ := \overline{\varphi}^+(\overline{F}^+) 
\simeq \Bbb P^1$, and $\overline{\varphi}^+$ 
is the blow-up of $\overline{X}^+$ with the center 
$\overline{E}^+$. 

Let $\overline{\varphi}^+({F'}^+) =: F^+$. 
By the Claim 1, we have the following claim in a similar way to 
the proof in 8.9:

\vskip 2mm

\flushpar
{\bf Claim 2.}\ \ There exist surjective morphisms 
$F' \to \Sigma_1$ and $F^+ \to \Bbb P^1$ which make 
the following diagram commutative:
$$
\CD
F' @. \simeq @. {F'}^+ 
@>{\overline{\varphi}}^+|_{{F'}^+}>> F^+ \\
@VVV @. @. @VVV \\
\quad \Sigma_1 @. @. 
@>>{\Bbb P^1\text{-bundle}}> \Bbb P^1
\endCD
\tag 8.10.8
$$

Let $\overline{W_{\lambda}} := \overline{\varphi}^{-1} 
(\overline{l_{\lambda}} - \{\overline{P}\})^{\text{---}}$ 
where $\{\overline{l_{\lambda}}\}_{\lambda \in \Lambda}$ 
$(\Lambda \simeq \Bbb P^1)$ is the complete 
family of lines in $\overline{E}$ passing through 
$\overline{P}$. Let $\overline{W_{\lambda}}^+$ 
be the proper transform of it in $\overline{\overline{X}}^+$. 
By (8.8.9)$_{m-1}$, 

\flushpar
(8.10.9) \, \, \, 
$\overline{W_{\lambda}}^+ \simeq \Bbb P^2$ 
and $\overline{\varphi}^+|_{{F'}^+}$ contracts exactly with 
$$
\overline{W_{\lambda}}^+ \cap {W_{\lambda}'}^+ 
=: m_{\lambda}^+ \simeq \Bbb P^1
$$
(where ${W_{\lambda}'}^+ \subset {F'}^+$ is the proper 
transform of $W_{\lambda}$ (8.8.9)$_m$ in 
$\overline{\overline{X}}^+$). 
This is the negative section of 
${W_{\lambda}'}^+ \simeq \Sigma_1$ (Claim 1). 

In particular, the morphism 
$F^+ \to \Bbb P^1$ in the Claim 2, (8.10.8) is a 
$\Bbb P^2$-bundle. 

The rest is exactly the same as in 8.9. 
\quad \qed

\vskip 2mm

8.8, 8.9 and 8.10 show Theorem 8.1. 

\proclaim{Corollary 8.11}\ \ 
(The description of the flip $g^+$)

Let $X \supset E \simeq \Bbb P^2
\overset{g}\to{\longrightarrow} Y \ni Q$ 
be a flipping contraction of Type (R), 
with 
$$
\operatorname{width} g = m (\geq 2). 
$$
Let $X^+ \supset E^+ 
\overset{g^+}\to{\longrightarrow} Y \ni Q$ 
be its flip (Theorem 8.1), where $E^+ := \operatorname{Exc} g^+$.
Then 

\flushpar
(1) \, \, \, $X^+$ is smooth, 

\flushpar
(2) \, \, \, $E^+ \simeq \Bbb P^1$, \, \, \, 

\flushpar
(3) \, \, \, $\operatorname{Bs} |-K_{X^+}| = E^+$, 

\flushpar
(4) \, \, \, 
$N_{E^+/X^+} \simeq \Cal O_{\Bbb P^1} \oplus 
\Cal O_{\Bbb P^1}(-1) \oplus \Cal O_{\Bbb P^1}(-2)$,  

\flushpar
(5) \, \, \, 
a general member $D^+ \in |-K_{X^+}|$ is smooth, and 

\flushpar
(6) \, \, \, $E^+ \subset D^+$ is a $(0, -2)$-curve, 
with 
$\operatorname{width} = \operatorname{width} (g^+|_{D^+}) 
= m$. \quad \qed
\endproclaim

\vskip 2mm

\flushpar
{\bf 8.12.}\ \ {\it (La Torre Pendente)\/}

Follow the induction procedure 8.8, with 
the proofs 8.9, 8.10, then we get the complete 
picture of the whole flip operation, starting from 
the given flipping contraction $X \supset E \simeq \Bbb P^2
\overset{g}\to{\longrightarrow} Y \ni Q$. Namely, 
the flip $X^+ \supset E^+ \overset{g^+}\to{\longrightarrow} Y \ni Q$ 
is obtained by blowing $X$ up $m$ times, and then 
down $m$ times, like: 
$$
X \leftarrow \overline{X}
\leftarrow ... \leftarrow \overline{X}^{(m-1)} 
\leftarrow \overline{X}^{(m)} \to {\overline{X}^{(m-1)}}^+ 
\to ... \to \overline{X}^+ \to X^+. 
$$
We name this {\it La Torre Pendente\/}, after M.~Reid's 
Pagoda. 

Let $D$ be a general smooth member of $|-K_X|$, and 
$\overline{D}^{(i)}$,  ${\overline{D}^{(i)}}^+$ be its 
proper transforms in $\overline{X}^{(i)}$,   
${\overline{X}^{(i)}}^+$, respectively, 
then 
$$
D \leftarrow \overline{D}
\leftarrow ... \leftarrow \overline{D}^{(m-1)} \leftarrow 
\overline{D}^{(m)} \to {\overline{D}^{(m-1)}}^+ \to ... \to 
\overline{D}^+ \to D^+ 
$$
is nothing but the whole Pagoda diagram of width $m$ (5.1). 
So our La Torre Pendente contains Pagoda as an anti-canonical 
divisor (and its proper transforms). 

Here we note that though Pagoda was symmetric with respect to 
the flop, ours is no more.

\vskip 2mm

At the end we shall give some examples, 
particularly related with {\it fiber type\/} 
contractions.  

\definition{Example 8.13}\ \ (S.~Mukai)

There exists a projective surjective morphism 
$h : X \to Z$ from a terminal Gorenstein 4-fold 
$X$ to a germ of a normal 3-fold singularity $(Z,Q)$ 
such that $h|_{X - h^{-1}(Q)} : X - h^{-1}(Q) \to Z - \{Q\}$ 
is a $\Bbb P^1$-bundle, while the fiber 
$F := h^{-1}(Q)$ at $Q$ is a union of two singular 
quadric surfaces (in $\Bbb P^3$): $F = F_1 \cup F_2$, 
meeting at a point $P$ which is the vertices of both. 
$X$ has a Gorenstein quotient singularity 
of type $\dfrac{1}{2}(1,1,1,1)$ at $P$ 
({\it i.e.\/} the quotient by the involution 
($\Bbb C^4, 0) \ni z \mapsto -z \in (\Bbb C^4, 0)$), and is smooth 
elsewhere. 

$h$ is factored by a flipping contraction 
$g : X \to Y$ which contracts either one of $F_1$, $F_2$. 

So the fact $\text{Exc }g \simeq \Bbb P^2$ is no longer true 
when $X$ has worse than complete intersection singularities. 
\enddefinition

\definition{Example 8.14}\ \ (A singular analogue of 
Mukai--Shepherd-Barron--Wi\'sniewski contraction (MSW))

\flushpar
An example of the following type of contraction is known 
(Mukai--Shepherd-Barron--Wi\'sniewski contraction (MSW), 
see [Kac1]):

A projective surjective morphism $h : X \to Z$ from a 
{\it smooth\/} 4-fold $X$ to a germ of a normal 3-fold 
singularity $(Z,Q)$, which is a $\Bbb P^1$-bundle over $Z - \{Q\}$, 
while the central fiber $F := h^{-1}(Q)$ is a union of 
two $\Bbb P^2$'s: $F = F_1 \cup F_2$, meeting at a point $P$. 

$h$ is factored by a flipping contraction 
$g : X \to Y$ contracting either one of $F_1$, $F_2$, 
which is of Kawamata type [Kaw4]. 
In this case, 
$$
(Z, Q) \simeq \{z_1 z_2 + z_3^2 + z_4^2 = 0\}. 
$$

\vskip 2mm

\flushpar
{\bf Fact.}\ \ ([Kac1])

This is the only contraction among those $h : X \to (Z, Q)$ 
from a smooth 4-fold $X$ to a 3-fold germ $(Z,Q)$, 
with $-K_X$ being $h$-ample, such that 

\flushpar
(0) \, \, \, $F := \dim h^{-1}(Q) = 2$, 

\flushpar
(1) \, \, \,there exists a rational curve $l \subset F$ with 
$(-K_X \, . \, l) = 1$, and 

\flushpar
(2) \, \, \, $h$ is a $\Bbb P^1$-bundle over $Z-\{Q\}$. \quad 
---

\vskip 2mm 

What we are going to construct is a similar $h$ from a 
{\it singular\/} $X$ which is factored by a flipping contraction 
of Type (R). 

The following construction is a generalization of Mukai who 
did in the case $m=1$, that is, the case that $X$ is smooth.

\vskip 2mm

\flushpar
{\bf Construction.}\ \ 
Let us start with the 3-fold flopping contraction 
$U \supset C \simeq \Bbb P^1 \to (Z, Q)$ 
of a $(0, -2)$-curve $C$, with width $m (\geq 2)$. 
Let $L$ be an irreducible smooth divisor on $U$ 
with $(L \, . \, C) = 1$, let 
$\Cal E := \Cal O_U \oplus \Cal O_U (-L)$, and 
let 
$$
X^+ := \Bbb P (\Cal E) \overset{\psi}\to{\longrightarrow} U. 
$$
This is a $\Bbb P^1$-bundle, with the section $D^+$ corresponding 
to $\Cal O_U (-L)$. 
Then ${E'}^+ := \psi^{-1}(C) \simeq \Sigma_1$.  
Let $E^+ := {E'}^+ \cap D^+$, then 
since $D^+ \supset E^+$ is mapped isomorphically onto 
$U \supset C$, $E^+$ is a $(0, -2)$-curve of 
width $m$ inside $D^+$.  
Moreover $E^+ \simeq \Bbb P^1$ is the negative section of 
${E'}^+$, and is a contractible curve in $X^+$: 
$X^+ \to Y$, 
with 
$$
N_{E^+/X^+} \simeq \Cal O_{\Bbb P^1} \oplus 
\Cal O_{\Bbb P^1}(-1) \oplus \Cal O_{\Bbb P^1}(-2). 
$$
Follow the induction argument 8.8 (with 8.9, 8.10) 
in the reverse direction starting from this 
$X^+ \supset E^+ \simeq \Bbb P^1$, then 
we get a flipping contraction of Type (R), 
with width $m$:
$$
X \supset E \simeq \Bbb P^2 
\to Y \ni Q 
$$
(the {\it inverse flip\/}). 
In particular there is exactly one singular point 
$P$ on $X$, which is of form: 
$$
(X,P) \simeq \{x_1 x_3 + x_2 x_4 + x_5^m = 0\} 
\supset E = \{x_3 = x_4 = x_5 = 0\}. 
$$
Let $E'$ be the proper transform of ${E'}^+ \subset X^+$ 
in $X$, then $E' \simeq \Bbb P^2$, and 
$$
E \cap E' = \{P\}, \, \, \, 
E' = \{x_1 = x_2 = x_5 = 0\}
$$
in the above description. 

Since everything in the above procedure is defined over 
$(Z,Q)$, and since we never touched the part 
$X^+ - {E'}^+$, 
there exists a contraction $h : X \to (Z,Q)$, 
which is a $\Bbb P^1$-bundle over $Z - \{Q\}$. 
As is easily seen the central fiber is a union 
of two $\Bbb P^2$'s meeting just at $P$: 
$h^{-1}(Q) = E \cup E'$.
Moreover 
$$
(Z, Q) \simeq \{z_1 z_2 + z_3^2 + z_4^{2m} = 0\}. 
$$
This $h$ gives the desired example. 
\enddefinition

\vskip 5mm

\head \S Appendix.\ On terminal complete intersection 
singularities \\
--- S.~Ishii's theorem \endhead

\flushpar
{\bf A.0.}\ \ 
The following implications hold in arbitrary dimensions: 
$$
\align
\text{Hypersurface singularity} \, \, \, &\Longrightarrow \, \, \, 
\text{Complete intersection singularity} \\ 
&\Longrightarrow \, \, \, \text{Gorenstein singularity}. 
\endalign
$$

\vskip 2mm

\flushpar
{\bf A.1.}\ \ (Reid [Re1,4])

For 3-dimensional terminal singularities, all of these 
are equivalent, and general hyperplane sections of those 
are Du Val singularities.  (These are called 
$cDV${\it-singularities\/}.)

\vskip 2mm

\flushpar
{\bf A.2.}\ \ In dimension 4, however, the reverse directions 
of both arrows in A.0 do not hold, even for terminal singularities 
(Mori [Mo5], Kac--Watanabe [KW]). 

\vskip 2mm

Here is a collection of references on terminal singularities in 
dimension 3 and 4: [D], [Kol1], [KS], [Mo3,4], [MoMoMo], [MS], 
[Re1,2,4]. 

\vskip 2mm

After the completion of this paper, the author learned 
the following theorem, which 
was originally conjectured independently by 
Yu.~Prokhorov and K.-i.~Watanabe, and has been solved by 
S.~Ishii:

\proclaim{Theorem A.3}\ \ (S.~Ishii)

Let $(Z, P)$ be a germ of an isolated terminal 
complete intersection singularity, of dimension $n$. 
Then 

\flushpar
(1) \qquad \qquad \qquad \qquad \qquad \, \, \, \,  
$\operatorname{emb.codim} (Z, P) \leq n-2$. 

\flushpar
(2) \, \, \, Let 
$$
(Z, P) \simeq \{f_1(x_1, ... , x_N) = \, ... \, = 
f_r(x_1, ... , x_N) = 0\} \subset (\Bbb C^N, 0), 
$$
where $r := \operatorname{emb.codim} (Z, P)$ and $N = n + r$ 
(so that this is an irredundant expression). 
Then
$$
\sum\limits_{i = 1}^r \operatorname{mult} f_i 
\leq n + r - 2. 
$$

\flushpar
(3) \, \, \, In particular if the equality holds in (1) 
{\it i.e.\/} $r = n - 2$, then 
$$
\operatorname{mult} f_i = 2 \quad (\forall i = 1, ... , n-2). 
$$
\endproclaim

This reproduces the following, which is a weaker version of 
the above mentioned theorem (A.1) of Reid [Re1] in the case $n = 3$: 

\proclaim{Corollary A.4}\ \ 
A 3-dimensional terminal complete intersection singularity 
is a hypersurface double point.  
\endproclaim

The following is the case $n = 4$ in Theorem A.3:

\proclaim{Corollary A.5}\ \ (S.~Ishii)

A 4-dimensional isolated terminal complete intersection singularity 
is either 

\flushpar
(a) \, \, \, a hypersurface double or a triple point in 
$(\Bbb C^5, 0)$, or 

\flushpar
(b) \, \, \, an intersection of two double hypersurfaces 
in $(\Bbb C^6, 0)$. 
\endproclaim

(See [Mo5] for some example.)

By virtue of this Corollary A.5, our Corollary 4.3 became 
unnecessary.


\newpage

\Refs
\widestnumber\key{KoMiMo3}
\ref
\key ABW
 \by M.~Andreatta, E.~Ballico and J.~Wi\'sniewski
 \paper Two theorems on elementary contractions
 \jour Math\. Ann\.
 \vol 297
 \pages 191--198
 \yr 1993
 \endref
\ref
\key AW1
 \by M.~Andreatta and J.~Wi\'sniewski
 \paper A note on nonvanishing and applications
 \jour Duke Math\. J\.
 \vol 72
 \pages 739--755
 \yr 1993
 \endref
\ref
\key AW2
 \by M.~Andreatta and J.~Wi\'sniewski
 \paper On good contractions of smooth varieties
 \jour Preprint
 \yr 1996
\endref
\ref
\key Be
 \by M.~Beltrametti
 \paper On $d$-folds whose canonical bundle is not numerically effective,
        according to Mori and Kawamata
 \jour Ann\. Mat\. Pura\. Appl\.
 \vol 147
 \pages 151--172
 \yr 1987
 \endref
\ref
\key Br
 \by E.~Brieskorn
 \paper Singular elements of semi-simple algebraic groups
 \jour Proc\. Int\. Cong\. Math\. Nice
 \vol 2
 \pages 279--284
 \yr 1970
\endref
\ref
\key Co
 \by A.~Corti
 \paper Semi-stable 3-fold flips
 \jour Preprint
 \yr 1994
\endref
\ref
\key D
 \by V.I.~Danilov
 \paper Birational geometry of toric 3-folds
  \jour Math\. USSR\. Izv\.
 \vol 21
 \pages 269--279
 \yr 1983
 \endref
\ref
\key E
 \by R.~Elkik
 \paper Singularit\'es rationelles et d\'eformations
 \jour Invent\. Math\.
 \vol 47
 \pages 139--147 
 \yr 1978
 \endref
\ref
\key Fra
 \by P.~Francia
 \paper Some remarks on minimal models
 \jour Compos\. Math\.
 \vol 40
 \yr 1980
 \pages 301--313
 \endref
\ref
\key Fri
 \by R.~Friedman 
 \paper \, Simultaneous resolution of threefold double points
 \jour Math\. Ann\.
 \vol 274
 \pages 671--689
 \yr 1986
 \endref
\ref
\key G1
 \by M.~Gross
 \paper Deforming Calabi-Yau threefolds
 \jour Preprint
 \yr 1995
 \endref
\ref
\key G2
 \by M.~Gross
 \paper Primitive Calabi-Yau threefolds
 \jour Preprint
 \yr 1995
\endref
\ref
\key H1
 \by E.~Horikawa
 \paper Deformations of holomorphic maps, I
 \jour J\. of Math\. Soc\. Japan
 \yr 1973
 \vol 25
 \pages 372--396
 \endref
\ref
\key H2
 \by E.~Horikawa
 \paper ibid, II
 \jour J\. of Math\. Soc\. Japan
 \yr 1974
 \vol 26
 \pages 647-667
 \endref
\ref
\key H3
 \by E.~Horikawa
 \paper ibid, III
 \jour Math\. Ann\.
 \yr 1976
 \vol 222
 \pages 275--282
 \endref
\ref
\key Io
 \by P.~Ionescu
 \paper Generalized adjunction and applications
 \jour Math\. Proc\. Camb\. Phil\. Soc\.
 \vol 99
 \pages 457--472
 \yr 1986
 \endref
\ref
\key Is
 \by S.~Ishii
 \paper On isolated Gorenstein singularities
 \jour Math\. Ann\.
 \vol 270
 \pages 541--554
 \yr 1985
 \endref
\ref
\key KW
 \by V.~Kac, K.-i.~Watanabe
 \paper Finite linear groups whose ring of invariants 
        is a complete intersection
 \jour Bull\. Amer\. Math\. Soc\.
 \yr 1982
 \vol 6, No2
 \pages 221--223
 \endref
\ref
\key Kac1
 \by Y.~Kachi
 \paper Extremal contractions from 4-dimensional manifolds to 3-folds
 \jour Ann\. di Pisa
 \yr To appear
 \endref 
\ref
\key Kac2
 \by Y.~Kachi
 \paper Flips from semi-stable 4-folds whose degenerate fibers are 
        unions of Cartier divisors which are terminal factorial 3-folds
 \jour Math\. Ann\.
 \yr To appear
 \endref
\ref
\key Kaw1
 \by Y.~Kawamata
 \paper Elementary contractions of algebraic 3-folds
 \jour Ann\. of Math\.
 \vol 119
 \pages 95--110
 \yr 1984
 \endref
\ref
\key Kaw2
 \by Y.~Kawamata
 \paper The cone of curves of algebraic varieties
 \jour Ann\. of Math\.
 \vol 119
 \pages 603--633
 \yr 1984
 \endref
\ref
\key Kaw3
 \by Y.~Kawamata
 \paper Crepant blowing-up of 3-dimensional canonical 
        singularities and its application to degenerations of surfaces
 \jour Ann\. of Math\.
 \vol 127
 \pages 93--163
 \yr 1988
 \endref
\ref
\key Kaw4
 \by Y.~Kawamata
 \paper Small contractions of four dimensional algebraic manifolds
 \jour Math\. Ann\.
 \vol 284
 \pages 595--600
 \yr 1989
 \endref
\ref
\key Kaw5
 \by Y.~Kawamata
 \paper On the length of an extremal rational curve
 \jour Invent\. math\.
 \vol 105
 \pages 609--611
 \yr 1991
 \endref
\ref
\key Kaw6
 \by Y.~Kawamata
 \paper Unobstructedness deformations -- A remark to a paper of Z.~Ran
 \jour J\. of Alg\. Geom\. 
 \vol 1
 \pages 183--190
 \yr 1992
 \endref
\ref
\key Kaw7
 \by Y.~Kawamata
 \paper Termination of log flips for algebraic 3-folds
 \jour Int\. J\. of Math\.
 \vol 3
 \pages 653--659
 \yr 1992
 \endref
\ref
\key Kaw8
 \by Y.~Kawamata
 \paper Semistable minimal models of threefolds in 
        positive or mixed characteristic
 \jour J\. of Alg\. Geom\.
 \vol 3
 \pages 463--491
 \yr 1994
 \endref
\ref 
\key Kaw9
 \by Y.~Kawamata
 \paper General hyperplane sections of nonsingular flops in dimension 3
 \jour Math\. Res\. Lett\.
 \vol 1
 \pages 49--52
 \yr 1994
 \endref
\ref
\key KaMaMa
 \by Y.~Kawamata, K.~Matsuda and K.~Matsuki
 \paper Introduction to the minimal model problem
 \jour Adv\. St\. Pure Math\.
 \vol 10
 \pages 283--360
 \yr 1987
 \endref
\ref 
\key KaMo
 \by S.~Katz and D.~Morrison
 \paper Gorenstein threefold singularities with small resolutions 
        via invariant theory of Weyl groups
 \jour J\. of Alg\. Geom\. 
 \vol 1
 \pages 449--530
 \yr 1992
 \endref
\ref
\key Kod
 \by K.~Kodaira
 \paper On stability of compact submanifolds of complex manifolds
 \jour Amer\. J\. of Math\.
 \yr 1963
 \vol 85
 \pages 79--84
\endref
\ref
\key Kol1
 \by J.~Koll\'ar
 \paper Flops
 \jour Nagoya Math\. J\.
 \vol 113
 \pages 15--36
 \yr 1989
\endref
\ref
\key Kol2
 \by J.~Koll\'ar
 \book Rational curves on algebraic varieties
 \publ Univ\. of Utah 
 \yr 1994
 \endref
\ref
\key Kol3
 \by J.~Koll\'ar
 \paper Flatness criteria
 \jour Preprint
 \yr 1994
 \endref
\ref
\key KoMiMo1
 \by J.~Koll\'ar, Y.~Miyaoka and S.~Mori
 \paper Rational curves on Fano varieties
 \jour Preprint, RIMS Kyoto
 \yr 1991
 \endref
\ref
\key KoMiMo2
 \by J.~Koll\'ar, Y.~Miyaoka and S.~Mori
 \paper Rationally connected varieties
 \jour J\. of Alg\. Geom\.
 \vol 1
 \pages 429--448
 \yr 1992
 \endref
\ref
\key KoMiMo3
 \by J.~Koll\'ar, Y.~Miyaoka and S.~Mori
 \paper Rational connectedness and boundedness of Fano manifolds
 \jour J\. of Diff\. Geom\.
 \vol 36
 \pages 765--779
 \yr 1992
 \endref
\ref
\key KoMo
 \by J.~Koll\'ar and S.~Mori
 \paper Classification of three dimensional flips
 \jour J\. of Amer\. Math\. Soc\.
 \vol 5
 \pages 533--703
 \yr 1992
 \endref
\ref
\key KS
 \by J.~Koll\'ar and N.~I.~Shepherd-Barron
 \paper Threefolds and deformations of surface singularities
 \jour Invent\. Math\.
 \vol 91
 \pages 299--338
 \yr 1988
 \endref
\ref
\key L
 \by H.~Laufer
 \book On $\Bbb C \Bbb P^1$ as an exceptional set
 \bookinfo in Ann\. of Math\. Studies
 \vol 100
 \pages 261--275
 \publ Princeton Univ\. Press
 \yr 1981
 \endref
\ref
\key Ma
 \by K.~Matsuki
 \paper Weyl groups and ...
 \jour ??
 \vol ??
 \pages ??
 \yr 1995
 \endref
\ref
\key Mo1
 \by S.~Mori
 \paper Projective manifolds with ample tangent bundles
 \jour Ann\. of Math\.
 \vol 110
 \pages 593--606
 \yr 1979
 \endref
\ref
\key Mo2
 \by S.~Mori
 \paper Threefolds whose canonical bundles are not numerically effective
 \jour Ann\. of Math\.
 \vol 116
 \pages 133--176
 \yr 1982
 \endref
\ref
\key Mo3
 \by S.~Mori
 \paper On 3-dimensional terminal singularities
 \jour Nagoya Math\. J\.
 \vol 98
 \pages 43--66
 \yr 1985
\endref
\ref
\key Mo4
 \by S.~Mori
 \paper Flip theorem and the existence of minimal models for 3-folds
 \jour J\. of Amer\. Math\. Soc\.
 \vol 1
 \pages 117--253
 \yr 1988
 \endref
\ref
\key Mo5
 \by S.~Mori
 \paper Dear Miles
 \jour Appendix to [Re5]
 \endref
\ref 
\key MoMoMo
 \by S.~Mori, D.~Morrison and I.~Morrison
 \paper On four dimensional terminal quotient singu- \linebreak
        larities
 \jour Preprint
 \endref
\ref
\key MS
 \by D.~Morrison and G.~Stevens
 \paper Terminal quotient singularities in dimension 3 and 4
 \jour Proc\. Amer\. Math\. Soc\.
 \vol 90
 \pages 15--20
 \yr 1984
 \endref
\ref 
\key Nak
 \by N.~Nakayama
 \paper On smooth exceptional curves in threefolds
 \jour Thesis, Univ\. of Tokyo
 \yr 1989
 \endref
\ref
\key Nam1
 \by Yo.~Namikawa
 \paper Smoothing Fano 3-folds
 \jour Preprint
 \yr 1994
 \endref
\ref
\key Nam2
 \by Yo.~Namikawa
 \paper Deformation theory of Calabi-Yau threefolds and 
        certain invariants of singularities
 \jour Preprint, Max-Planck-Inst\. 
 \yr 1995
 \endref
\ref
\key Nam3
 \by Yo.~Namikawa
 \paper Stratified local moduli of Calabi-Yau threefolds
 \jour Preprint, Max-Planck-Inst\. 
 \yr 1995
 \endref
\ref
\key NS
 \by Yo.~Namikawa and J.H.M.~Steenbrink
 \paper Global smoothing of Calabi-Yau threefolds
 \jour Pre- print, Max-Planck-Inst\. 
 \yr 1995
 \endref
\ref
\key OSS
 \by C.~Okonek, M.~Schneider and H.~Spindler
 \book Vector bundles on complex projective spaces
 \bookinfo Progress math\.
 \vol 3
 \publ Birkh\"auser
 \yr 1980
\endref
\ref
\key Pi1
 \by H.~Pinkham
 \book R\'esolution simultan\'ee de points doubles rationels
 \bookinfo Lect\. Notes Math\.
 \vol 777
 \publ Springer 
 \yr 1980
 \endref
\ref
\key Pi2
 \by H.~Pinkham
 \paper Factorization of birational maps in dimension 3
 \jour Proc\. Symp\. Pure\. Math\.
 \vol 40
 \yr 1983
 \pages 343--371
 \endref
\ref
 \key Pr1
 \by Y.~Prokhorov
 \paper On Q-Fano fiber spaces with two-dimensional base
 \jour Preprint, Max-Planck-Inst\. 
 \yr 1995
 \endref
\ref
 \key Pr2
 \by Y.~Prokhorov
 \paper On extremal contractions from threefolds to surfaces
 \jour Preprint, Univ\. of Warwick\.
 \yr 1995
 \endref
\ref
\key Ra1
 \by Z.~Ran
 \book Deformation of maps
 \bookinfo Algebraic curves and projective geometry, LNM 
 \linebreak
 \vol 1389
 \pages 246--253
 \publ Springer
 \yr 1989
 \endref
\ref
\key Ra2
 \by Z.~Ran
 \paper Deformations of manifolds with torsion or negative 
        canonical bundle
 \jour J\. of Alg\. Geom\.
 \vol 1
 \pages 279--291
 \yr 1992
 \endref
\ref
\key Re0
 \by M.~Reid
 \book Canonical 3-folds
 \bookinfo in Journ\'ee de G\'eom\'etrie Alg\'ebrique d'Angers
 \ed A.~Beauville
 \publ Sijthoff and Noordhoff
 \yr 1980
 \pages 273--310
 \endref
\ref
\key Re1
 \by M.~Reid
 \paper Minimal models of canonical 3-folds
 \jour Adv\. St\. Pure Math\.
 \vol 1
 \pages 131--180
 \yr 1983
 \endref
\ref
\key Re2
 \by M.~Reid
 \book Decomposition of toric morphisms
 \bookinfo in Arithmetic and Geometry II, Progress Math\.
 \vol 36
 \publ Birkha\"user
 \pages 395--418
 \yr 1983
\endref
\ref
\key Re3
 \by M.~Reid
 \paper Projective morphisms according to Kawamata
 \jour Preprint
 \yr 1983
 \endref
\ref
\key Re4
 \by M.~Reid
 \paper Young person's guide to canonical singularities
 \jour Proc\. Symp\. Pure Math\.
 \vol 46, \, Vol\. 1
 \pages 345--414
 \yr 1987
 \endref
\ref
\key Re5
 \by M.~Reid
 \paper What is a flip?
 \jour Preprint
 \yr 1993
 \endref
\ref
\key Sc
 \by M.~Schlessinger
 \paper Rigidity of quotient singularities
 \jour Invent\. Math\.
 \vol 14
 \pages 17--26
 \yr 1971
\endref
\ref
\key Sh1
 \by V.V.~Shokurov
 \paper The nonvanishing theorem
 \jour Math\. USSR\. Izv\.
 \vol 26, No.3
 \pages 591--604
 \yr 1986
 \endref
\ref
\key Sh2
 \by V.V.~Shokurov
 \paper 3-fold log flips
 \jour Math\. USSR\. Izv\.
 \vol 40, No.1
 \pages 95--202
 \yr 1993
 \endref
\ref
\key Utah
 \by J.~Koll\'ar et.al.
 \book Flips and abundance for algebraic threefolds
 \bookinfo Ast\'erisque
 \vol 211
 \publ Soc\. Math\. de France
 \yr 1992
 \endref
\ref
\key V
 \by A.~Van de Ven
 \paper On uniform vector bundles
 \jour Math\. Ann\.
 \vol 195
 \pages 245--248
 \yr 1972
 \endref
\ref
\key W
 \by J.~Wi\'sniewski
 \paper On contraction of extremal rays of Fano manifolds
 \jour J\. Reine\. Angew\. Math\.
 \vol 417
 \pages 141--157
 \yr 1991
 \endref
\endRefs
\enddocument